\documentclass[reqno]{amsart}

\usepackage{amssymb,epic,eepic}

\usepackage[cmtip,matrix,arrow]{xy}
\CompileMatrices
\newdir{((}{{\lhook\kern-.1em}}
\newdir{))}{{\rhook\kern-.1em}}
\newdir{ >}{{}*!/-5pt/@{>}}

\numberwithin{equation}{section}
\newcommand\note[1]%
{$^\dagger$\marginpar{\footnotesize{$^\dagger${#1}}}} 

\newtheorem{theorem}{Theorem}[section]
\newtheorem{proposition}[theorem]{Proposition}
\newtheorem{lemma}[theorem]{Lemma}
\newtheorem{corollary}[theorem]{Corollary}
\newtheorem{addendum}[theorem]{Addendum}

\theoremstyle{definition}
\newtheorem{definition}[theorem]{Definition}

\theoremstyle{remark} 

\newtheorem{example}[theorem]{Example}

\newcommand\eu{\mathfrak} 
\newcommand\lie{\mathfrak} 
\renewcommand\a{\lie{a}} 
\newcommand\h{\lie{h}} 
\newcommand\g{\lie{g}}
\renewcommand\k{\lie{k}} 
\renewcommand\t{\lie{t}}
\newcommand\W{\eu{W}} 
\newcommand\bb{\mathbb} 
\newcommand\N{\bb{N}}
\newcommand\Z{\bb{Z}} 
\newcommand\Q{\bb{Q}}
\newcommand\R{\bb{R}} 
\newcommand\C{\bb{C}}
\newcommand\ca{\mathcal} 
\newcommand\F{\ca{F}} 
\newcommand\K{\ca{K}}
\renewcommand\P{\ca{P}}
\renewcommand\S{\ca{S}} 
\newcommand\eps{\varepsilon}

\newcommand\hull{\operatorname{hull}}
\newcommand\cone{\operatorname{cone}}
\renewcommand\star{\operatorname{star}}
\newcommand\supp{\operatorname{supp}}
\newcommand\inter{\operatorname{int}}
\newcommand\depth{\operatorname{depth}}
\newcommand\vol{\operatorname{vol}}
\newcommand\codim{\operatorname{codim}}
\newcommand\ind{\operatorname{index}}
\newcommand\Td{\operatorname{Td}}
\newcommand\Ch{\operatorname{Ch}} 
\newcommand\D{\operatorname{D}}
\newcommand\RR{\operatorname{RR}}
\newcommand\Bl{\operatorname{Bl}}
\newcommand\Rep{\operatorname{Rep}}
\newcommand\Ind{\operatorname{Ind}}
\newcommand\pr{\operatorname{pr}}
\newcommand\Diff{\operatorname{Diff}}
\newcommand\Iso{\operatorname{Iso}}
\newcommand\Aut{\operatorname{Aut}}
\newcommand\Hom{\operatorname{Hom}}
\newcommand\Car{\operatorname{Car}} 
\newcommand\hor{\mathrm{hor}} 

\newcommand\Spin{\operatorname{Spin}}

\newcommand\SU{\operatorname{SU}}
\newcommand\U{\operatorname{U}}
\renewcommand\Sp{\operatorname{Sp}}

\newcommand\qu{/\kern-.7ex/} 
\newcommand\bigqu{\big/\kern-.85ex\big/} 
\newcommand\longhookrightarrow{\lhook\joinrel\longrightarrow}

\newcommand\pre{\preccurlyeq}
\newcommand\suc{\succcurlyeq}

\newcommand\dirac{/\kern-1.2ex\partial} 
\renewcommand\d{\mathrm{d}} 
\newcommand\tplus{\lie t^*_+} 
\newcommand\inv{^{-1}} 
\newcommand\prin{_{\mathrm{prin}}} 
\newcommand\gen{_{\mathrm{gen}}} 
\newcommand\shift{^{\mathrm{shift}}}
\newcommand\id{\mathrm{id}} 
\newcommand\vertical{\mathrm{vert}}
\newcommand\pt{\mathrm{pt}}
\newcommand\basic{\mathrm{basic}}

\title{Singular Reduction and Quantization} 

\author[Eckhard Meinrenken]{Eckhard Meinrenken$^{*}$} 
\thanks{$^*$Partially supported by a Feodor-Lynen Fellowship of the
Humboldt Foundation}
\address{Massachusetts Institute of Technology, Department of
Mathematics, Cambridge, Massachusetts 02139-4307} 
\email{mein@math.mit.edu}
\author[Reyer Sjamaar]{Reyer Sjamaar$^{\dagger}$}
\thanks{$^\dagger$Partially supported by an Alfred P. Sloan Research
Fellowship, by NSF Grant DMS-9703947, and by the Institut Fourier,
Grenoble}
\address{Cornell University, Department of Mathematics, Ithaca, New
York 14853-7901}
\email{sjamaar@math.cornell.edu}
\date{July 1997}

\begin{document} 


\begin{abstract}
Consider a compact prequantizable symplectic manifold $M$ on which a
compact Lie group $G$ acts in a Hamiltonian fashion.  The
``quantization commutes with reduction'' theorem asserts that the
$G$-invariant part of the equivariant index of $M$ is equal to the
Riemann-Roch number of the symplectic quotient of $M$, provided the
quotient is nonsingular.  We extend this result to singular symplectic
quotients, using partial desingularizations of the symplectic quotient
to define its Riemann-Roch number.  By similar methods we also compute
multiplicities for the equivariant index of the dual of a prequantum
bundle, and furthermore show that the arithmetic genus of a
Hamiltonian $G$-manifold is invariant under symplectic reduction.
\end{abstract}


\maketitle

\tableofcontents

\section{Introduction}

Consider a compact symplectic manifold $(M,\omega)$ on which a compact
Lie group $G$ acts in a Hamiltonian fashion with equivariant moment
map $\Phi\colon M\to\g^*$.  A fundamental result due to Meyer and
Marsden-Weinstein says that if $0$ is a regular value of $\Phi$, then
the \emph{symplectic quotient} (also known as the \emph{reduced
space})
$$M_0=\Phi\inv(0)/G$$
is a symplectic orbifold, a symplectic space with finite-quotient
singularities.  However, if $0$ is a singular value of the moment map,
the symplectic quotient usually has more complicated singularities.
Singular symplectic quotients have been the subject of intensive study
over the past fifteen years.  For instance, it was proved by Arms et
al.~\cite{ar:sy} and Otto \cite{ot:re} that $M_0$ admits a finite
decomposition into smooth symplectic manifolds, labelled by orbit
types of $M$.  Sjamaar and Lerman \cite{sj:st} proved that this
decomposition is a stratification in the sense of Thom-Mather and gave
local normal forms for the singularities.

The object of this paper is twofold.  The first goal is to define
Riemann-Roch numbers of singular symplectic quotients $M_0$ with
coefficients in certain complex line bundles.  These bundles include
the trivial line bundle, the Riemann-Roch number of which we call the
arithmetic genus of $M_0$, and the prequantum line ``bundle'', the
Riemann-Roch number of which is the dimension of the quantization of
$M_0$.  (The prequantum line ``bundle'' is not a genuine fibre bundle
but an orbibundle.)

The second goal is to understand how these Riemann-Roch numbers are
related to the corresponding characteristic numbers of $M$.  Our
incentive is to extend the ``quantization commutes with reduction''
theorems of \cite{me:sym,je:lo2} to the singular case.  These theorems
arose from a conjecture of Guillemin and Sternberg \cite{gu:ge}.

A major obstacle to attaining our first goal is the fact that on a
singular space there is no obvious way to define a Riemann-Roch number
as the index of an elliptic operator.  To make matters worse,
symplectic quotients are seldom complex algebraic or even analytic
varieties, so that holomorphic Euler characteristics do not make sense
and the Riemann-Roch formulas of Baum et al.\ \cite{ba:ri2} and Levy
\cite{le:rr} do not apply.  Our attempt to surmount this obstacle
consists in (partially) resolving the singularities of $M_0$ and
defining the Riemann-Roch numbers of $M_0$ to be equal to those of its
desingularization.  This raises the question whether the result
depends on the way in which we resolve the singularities.  In contrast
to the situation in algebraic geometry this is not an easy question,
and the answer we find is incomplete.  One way of desingularizing a
symplectic quotient was discovered by Kirwan \cite{ki:pa}.  Another
way is simply to shift the value of the moment map to a nearby generic
value.  Our result says that these two desingularization methods lead
to the same Riemann-Roch numbers.  (Neither method yields a
desingularization of $M_0$ in the strict sense of the word, but only a
partial desingularization, which may have finite-quotient
singularities.)

We are far more successful in winning our second objective.  Let $L$
be a $G$-equivariant complex line bundle on $M$.  Let $\RR(M,L)$ be
the equivariant index of $M$ with coefficients in $L$, that is the
pushforward of $L$ to a point, viewed as an element of the equivariant
$K$-theory of a point.  Under favourable circumstances, e.~g.\ if $L$
is the trivial bundle or the prequantum bundle, $L$ induces a line
``bundle'' $L_0$ on the quotient $M_0$.  This enables us to define the
Riemann-Roch number of $M_0$ with coefficients in $L_0$ by means of
either of the two desingularization processes referred to above.  Our
results include:

1.  if $L$ is trivial, then $\RR(M,L)=\RR(M_0,L_0)$, i.~e.\ the
arithmetic genus of $M$ is equal to the arithmetic genus of $M_0$;

2. if $L$ is the prequantum line bundle, then
$\RR(M,L)^G=\RR(M_0,L_0)$, i.~e.\ quantization commutes with
reduction.  (The superscript $G$ denotes $G$-invariants.)  The latter
result leads to a geometric formula for the multiplicities of all
irreducible representations occurring in the quantization $\RR(M,L)$.
We obtain similar results for the ``negative'' quantization
$\RR(M,L\inv)$.

It also turns out that the multiplicities depend only on the weights
of the action of $T$ on the fibres of $L$ at the fixed-point set
$M^T$, where $T$ denotes a maximal torus of $G$.  This observation
enables us to generalize our results to a larger class of bundles.
  
The method of proof is an extension to the singular case of techniques
developed in \cite{me:sym} for the proof of the Guillemin-Sternberg
conjecture, the key tool being a gluing formula that relates the
equivariant index of $M$ to equivariant indices of simpler spaces
obtained from $M$ by symplectic cutting in the sense of Lerman
\cite{le:sy2}.  The gluing formula is an application of the
Atiyah-Segal-Singer equivariant index formula.

Because the operations of symplectic cutting and partial
desingularization give rise to orbifolds rather than manifolds and
because many ``bundles'' we shall consider are orbibundles rather than
bundles, we are obliged to place our discussion within the wider
framework of Hamiltonian $G$-orbifolds.  While this presents few
conceptual difficulties, the technicalities are sometimes rather
involved.  In the interests of clarity and brevity we shall at some
points treat in complete detail only the manifold case and indicate
succinctly how to extend the argument to orbifolds.  The relevant
versions of the index formula in this category are due to Kawasaki
\cite{ka:ri}, Vergne \cite{ve:eq} and Duistermaat \cite{du:he}.

The organization of this paper is as follows.  Section
\ref{section:results} contains detailed statements and a discussion of
our main results.  In Section \ref{section:singularquotients} we
review the local structure of singular symplectic quotients.  In
Section \ref{section:desingularization} we describe the two known
methods for desingularizing symplectic quotients.  We prove in detail
that Kirwan's partial desingularization is well-defined up to
deformation equivalence and discuss briefly how it is related to
shift-desingularizations.  We then present the proofs of our main
results, first in the abelian case (Section \ref{section:abelian}),
then in the nonabelian case (Section \ref{section:nonabelian}).  At
several points we illustrate our results by applying them to Delzant
spaces, a class of toric varieties with symplectic structures.  These
not only serve as an interesting example, but also play an important
part in symplectic cutting.  Appendix \ref{section:normalform}
contains a number of technical results concerning blowups and
constant-rank embeddings.  In Appendix \ref{section:product} we prove
a product formula for the Todd class of an almost complex fibre
bundle, which generalizes a classical result of Borel.  A table
listing our notational conventions is provided in Appendix
\ref{section:notation}.

\section{Statement of results}\label{section:results} 

In Section \ref{subsection:preliminaries} we introduce notation and
some basic notions concerning Hamiltonian actions.  This is standard
material except for Definitions \ref{definition:quasiregular} and
\ref{definition:aelt}.  Sections
\ref{subsection:singular}--\ref{subsection:dual} are a compendium of
the chief results of this paper.  It is important to note that these
results hold without any regularity assumptions on the values of the
moment map.

\subsection{Preliminaries}\label{subsection:preliminaries}

Throughout this paper $G$ denotes a compact connected Lie group.  We
choose once and for all a maximal torus $T$ of $G$ and a (closed) Weyl
chamber $\tplus$ in $\lie t^*$ and denote by $\lie W$ the Weyl group
$N_G(T)/T$.  Throughout $(M,\omega,\Phi)$ designates a connected
symplectic orbifold on which $G$ acts in a Hamiltonian fashion with a
proper and $G$-equivariant moment map $\Phi\colon M\to\lie g^*$.
(Many of our results hold only for compact $M$, but it is important
for technical reasons to allow $M$ to be noncompact.)  Our sign
convention for the moment map is as follows:
$$
\d\langle\Phi,\xi\rangle=\iota(\xi_M)\omega,
$$
where $\xi_M$ denotes the fundamental vector field induced by
$\xi\in\lie g$.  The pair $(\omega,\Phi)$ is an \emph{equivariant
symplectic form\/} on $M$.  An \emph{isomorphism\/} between two
Hamiltonian $G$-orbifolds is a $G$-equivariant symplectomorphism that
intertwines the moment maps on the two spaces.  Some basic material on
symplectic orbifolds can be found in \cite{le:co,le:ha,me:sym}.  Our
conventions concerning orbifolds are as in \cite{me:sym}.  (In
particular, the structure group of an orbifold at a point is not
required to act effectively, and the structure group of a suborbifold
at a point is the same as the structure group of the ambient orbifold
at that point.)  The set $\Phi(M)\cap\tplus$ will be denoted by
$\Delta$.  By a theorem of Kirwan \cite{ki:con} (cf.\ also
\cite{le:co,sj:co}) it is a convex rational polyhedron, referred to as
the \emph{moment polyhedron}.  For $\mu\in\g^*$ let $G_\mu$ be the
stabilizer group of $\mu$ under the coadjoint action and let $G\mu$ be
the coadjoint orbit through $\mu$.  The compact space
$$M_\mu=\Phi\inv(\mu)/G_\mu\cong\Phi\inv(G\mu)/G$$
is the \emph{symplectic quotient\/} of $M$ at level $\mu$.  The
symplectic quotient at level $0$ plays a special role and
(particularly in situations where there is more than one group acting)
will also be denoted by
$$M\qu G=M_0.$$
The symplectic quotients of $M$ have a natural stratification by
symplectic orbifolds determined by the infinitesimal orbit types of
$M$.  (See Section \ref{subsection:localmodels}.)

\begin{definition}\label{definition:quasiregular}
A point $\mu\in\lie g^*$ is a \emph{quasi-regular\/} value of $\Phi$
if the $G$-orbits in $\Phi\inv(G\mu)$ all have the same dimension.
\end{definition} 

Equivalently, $\mu$ is quasi-regular if the rank of $\Phi$ is constant
on $\Phi\inv(G\mu)$, or the dimension of the stabilizer $G_m$ is the
same for all $m$ in $\Phi\inv(G\mu)$, that is to say if
$\Phi\inv(G\mu)$ is contained in a single infinitesimal orbit type
stratum.  Consequently, if $\mu$ is a quasi-regular value, then the
orbifold stratification of $M_\mu$ consists of one piece only, and
$M_\mu$ is therefore a symplectic orbifold.  Here are some examples of
quasi-regular values: weakly regular values (i.~e.\ values $\mu$ for
which $\Phi$ intersects $\{\mu\}$ cleanly; see Proposition
\ref{proposition:weakquasi}), points in $\Delta$ of maximal norm (see
Lemma \ref{lemma:norm}) and, if $G$ is abelian, vertices of the moment
polytope.  See Section \ref{subsection:delzant1} for more examples.

Let $L$ be a $G$-equivariant complex line orbibundle (also known as an
orbifold line bundle) on $M$.  For $\mu\in\lie g^*$ define $L_\mu$ to
be the quotient of the restriction of $L$ to $\Phi\inv(G\mu)$,
$$L_\mu=\bigl(L|_{\Phi\inv(G\mu)}\bigr)\big/G.$$
For $\mu=0$ we shall also use the notation
$$
L\qu G=L_0.
$$
Suppose for a moment that $L$ is a true line bundle.  When is $L_\mu$
a topologically locally trivial complex line bundle on the topological
space $M_\mu$?  It is not hard to see that this is the case if and
only if for all $m\in \Phi\inv(G\mu)$ the stabilizer $G_m$ acts
trivially on the fibre $L_m$, in other words $L|_{\Phi\inv(G\mu)}$ is
$G$-equivariantly locally trivial.  We shall make a slightly weaker
assumption.

\begin{definition}\label{definition:aelt}
The line orbibundle $L$ is \emph{almost equivariantly locally
trivial\/} at $m$ if the action of the identity component of $G_m$ on
$L_m$ is trivial.  It is \emph{almost equivariantly locally trivial at
level\/} $\mu$ if it is almost equivariantly trivial at all $m\in
\Phi\inv(G\mu)$.
\end{definition} 

For instance, if $\mu$ is a regular value of $\Phi$, then $G$ acts
locally freely on $\Phi\inv(G\mu)$, so $L$ is almost equivariantly
locally trivial at $\mu$.  If $L$ is almost equivariantly locally
trivial at $\mu$, then the fibres of the induced map $L_\mu\to M_\mu$
are finite quotients of $\C$.  If in addition $\mu$ is a quasi-regular
value, then $L_\mu$ is a line orbibundle over the orbifold $M_\mu$.

Now assume that $M$ is compact.  Choose a $G$-invariant almost complex
structure $J$ on $M$ which is compatible with $\omega$ in the sense
that the symmetric bilinear form $\omega({\cdot},J{\cdot})$ is a
Riemannian metric.  Let $\bar{\partial}_L$ be the Dolbeault operator
with coefficients in $L$.  Also choose a $G$-invariant Hermitian fibre
metric on $L$.  The \emph{Dolbeault-Dirac operator\/} on $M$ with
coefficients in $L$ is defined by
$\dirac_L=\sqrt2\,(\bar{\partial}_L+\bar{\partial}_L^*)$, considered
as an operator from $\Omega^{0,\text{even}}(M,L)$ to
$\Omega^{0,\text{odd}}(M,L)$.  The \emph{equivariant Riemann-Roch
number} of $M$ with coefficients in $L$ is the equivariant index of
$\dirac_L$,
$$\RR(M,L)=\ind_G(\dirac_L),$$
viewed as an element of $\Rep G$, the character ring of $G$.  An
alternative definition goes as follows.  Every $G$-orbifold with an
invariant almost complex structure carries a canonical invariant
$\Spin_c$-structure.  The $\Spin_c$-Dirac operator of $M$ with
coefficients in $L$ has the same principal symbol as $\dirac_L$ (see
e.~g.\ \cite{du:he}), and therefore has the same equivariant index.

The character $\RR(M,L)$ does not depend on the choice of $J$ (because
any two compatible almost complex structures are homotopic), nor on
the choice of the fibre metric on $L$.  Indeed, $\RR(M,L)$ depends
only on the homotopy class of the almost complex structure and the
equivariant Chern class of $L$.

Let $\Lambda$ be the integral lattice $\ker\bigl(\exp|_{\lie t}\bigr)$
of $\lie t$.  Then $\Lambda^*=\Hom_\Z(\Lambda,\Z)$ is the lattice of
real infinitesimal weights and $\Lambda^*_+=\Lambda^*\cap\tplus$ is
the set of real dominant weights.  The \emph{multiplicity function\/}
of $L$ is the function $N_L$ (also denoted by $N$) on $\Lambda^*_+$
with values in $\Z$ defined by the orthogonal decomposition
\begin{equation}\label{equation:rr}
\RR(M,L)=\sum_{\mu\in\Lambda^*_+}N_L(\mu)\,\chi_\mu,
\end{equation}
where $\chi_\mu$ denotes the character of the irreducible
representation with highest weight $\mu$.

Even if $M$ is not compact, its symplectic quotients are, so if $\mu$
is a quasi-regular value of $\Phi$ and $L$ is almost equivariantly
locally trivial at $\mu$, then the Riemann-Roch number
$\RR(M_\mu,L_\mu)$ of the orbifold $M_\mu$ with coefficients in the
orbibundle $L_\mu$ is well-defined.  We shall now discuss how to
define $\RR(M_\mu,L_\mu)$ even when $M_\mu$ is not an orbifold.

\subsection{The singular case}\label{subsection:singular}

Consider a value $\mu$ of the moment map that is not quasi-regular.
Let $L$ be a $G$-equivariant line orbibundle on $M$ and suppose that
$L$ is almost equivariantly locally trivial at $\mu$.  Every point in
$M_\mu$ has an open neighbourhood $O$ which can be written as a
quotient of a space $\tilde{O}$ by a finite group $\Gamma$ such that
$L_\mu$ is the quotient by $\Gamma$ of a $\Gamma$-equivariant line
bundle on $\tilde{O}$.  (See Section \ref{subsection:linenearby}.)  We
shall call $L_\mu$ a line orbibundle over $M_\mu$, even though the
base space $M_\mu$ need not be an orbifold.  Now let $\tilde{M}_\mu$
be Kirwan's canonical partial desingularization of $M_\mu$.  To
construct $\tilde{M}_\mu$ one first performs a sequence of equivariant
symplectic blowups of a $G$-invariant neighbourhood $U$ of
$\Phi\inv(G\mu)$ to obtain a Hamiltonian $G$-orbifold
$(\tilde{U},\tilde{\omega},\tilde{\Phi})$ with the property that $\mu$
is a quasi-regular value of $\tilde{\Phi}$.  The space $\tilde{M}_\mu$
is then the reduction of $\tilde{U}$ at $\mu$, which is a symplectic
orbifold.  The pullback bundle $\tilde{L}$ on $\tilde{U}$ is almost
equivariantly locally trivial at $\mu$, so that $\tilde{L}_\mu$ is a
line orbibundle over $\tilde{M}_\mu$ and the Riemann-Roch number
$\RR(\tilde{M}_\mu,\tilde{L}_\mu)$ makes sense.

\begin{definition}\label{definition:singular}
$\RR(M_\mu,L_\mu)=\RR\bigl(\tilde{M}_\mu,\tilde{L}_\mu\bigr)$.
\end{definition}

In the algebraic case, where $M$ is a complex projective orbifold and
the quotients $M_\mu$ are complex projective varieties, this equality
is not the definition of $\RR(M_\mu,L_\mu)$, but a consequence of the
fact that the $M_\mu$ have rational singularities.  (Cf.\
\cite{sj:ho}.)

One problem with this definition is that the symplectic structures on
$\tilde{U}$ and $\tilde{M}_\mu$ depend on a long list of choices.
However, in Section \ref{subsection:deform} we prove the following
result.

\begin{theorem}\label{theorem:deform}
The germ at $\tilde{\Phi}\inv(G\mu)$ of the triple $\bigl(\tilde
U,\tilde\omega,\tilde\Phi\bigr)$ is unique up to deformation
equivalence.  This implies that the symplectic structure on
$\tilde{M}_\mu$ is unique up to deformation equivalence and hence that
the Riemann-Roch numbers of $\tilde{M}_\mu$ are well-defined.
\end{theorem}

A more difficult and as yet unresolved problem is that the partial
resolution $\tilde{M}_\mu$ depends on the way $M_\mu$ is written as a
quotient.  It is conceivable that $M_\mu$ could be presented in a
different way as a quotient $M'_{\mu'}$ of a Hamiltonian $G'$-orbifold
$M'$ and it is \emph{a priori\/} unclear if $\tilde{M}_\mu$ and
$\tilde{M}'_{\mu'}$ have the same Riemann-Roch numbers.  The theorem
below offers limited evidence that $\RR(M_\mu,L_\mu)$ is independent
of the partial resolution.  By Lemma \ref{lemma:aelt} almost
equivariant local triviality is an open condition, so for all $\nu$
near $\mu$ the quotient $L_\nu$ is a line orbibundle over $M_\nu$ and
$\RR(M_\nu,L_\nu)$ is well-defined (as
$\RR(\tilde{M}_\nu,\tilde{L}_\nu)$ if $M_\nu$ is singular).

\begin{theorem}\label{theorem:independent}
If $L$ is almost equivariantly locally trivial at level $\mu$\upn,
then
$$\RR(M_\mu,L_\mu)=\RR(M_\nu,L_\nu)$$
for all $\nu\in\Phi(M)$ sufficiently close to $\mu$.
\end{theorem}

The proof is in Section \ref{subsection:induction}.  If $M_\mu$ and
$M_\nu$ have the same dimension (e.~g.\ if $\mu$ and $\nu$ are both in
the interior of $\Delta$), then $\tilde{M}_\mu$ and $\tilde{M}_\nu$
are two different partial desingularizations of $M_\mu$, and we shall
refer to $\tilde{M}_\nu$ as a \emph{shift desingularization\/} of
$M_\mu$.  Theorem \ref{theorem:independent} asserts that the shift
desingularizations of $M_\mu$ give the same Riemann-Roch numbers as
the canonical partial desingularization.

In the next sections we shall make a detailed comparison between the
virtual character $\RR(M,L)$ and the numbers $\RR(M_\mu,L_\mu)$ for
three different types of bundle.

\subsection{Rigid bundles}\label{subsection:arithmetic}
  
A $G$-equivariant line orbibundle $L$ on $M$ is called \emph{rigid\/}
(or \emph{$G$-rigid\/}) if the action of $T$ on $L|_{M^T}$ is trivial.
This condition is obviously independent of the choice of the maximal
torus $T$.  A rigid bundle is almost equivariantly locally trivial
everywhere by Lemma \ref{lemma:line}, so that $\RR(M_\mu,L_\mu)$ is
well-defined for all $\mu$.

\begin{example}
The equivariantly trivial line bundle $\C$ is rigid.  Its equivariant
index is the \emph{arithmetic\/} or \emph{Todd genus\/} of $M$.  The
induced bundle $\C_\mu$ is of course the trivial line bundle on
$M_\mu$.
\end{example}

Notice that the definition of a rigid orbibundle makes sense for an
arbitrary almost complex $G$-orbifold.  The term ``rigid orbibundle''
is inspired by an observation of Lusztig (see \cite{at:sp}) stating
that the arithmetic genus of a almost complex $G$-manifold is rigid
(i.~e.\ a constant character).  Along the lines of
\cite{at:sp,bo:on,du:sy,ko:ap} we shall prove the following stronger
result.

\begin{theorem}\label{theorem:rigid} 
\begin{enumerate}
\item\label{part:almostcomplex}
Let $\ca M$ be a compact almost complex $G$-orbifold and let $\ca L$
be a rigid orbibundle on $\ca M$. Then the character $\RR(\ca M,\ca
L)$ is constant.  Choose a generic $\xi\in\lie t$.  Then
$$ 
\RR(\ca M,\ca L)= \sum_{\substack{F\\ \xi\in\check\ca C_F}}\RR(F,\ca
L|_F),
$$
where the summation is over all connected components $F$ of $\ca M^T$
such that $\langle\alpha,\xi\rangle<0$ for all orbiweights $\alpha$ of
the $T$-action on the normal bundle of $F$.
\item\label{part:hamiltonian}
Let $M$ be a compact Hamiltonian $G$-orbifold and let $L$ be a rigid
orbibundle on $M$.  Then
\begin{equation}\label{equation:rigidquotient}
\RR(M,L)=\RR(M_\mu,L_\mu)
\end{equation}
for all $\mu\in\Phi(M)$.  In particular\upn, $\RR(M_\mu,L_\mu)$ does
not depend on the value of $\mu$.  In the presence of an action of a
compact connected Lie group $H$ on $M$ and $L$ that commutes with the
action of $G$\upn, \eqref{equation:rigidquotient} holds as an equality
of virtual characters of $H$.
\end{enumerate}
\end{theorem}

The proof of \ref{part:almostcomplex} is in Section
\ref{subsection:almostcomplex} and the proof of \ref{part:hamiltonian}
is in Section \ref{subsection:induction}.  Setting $L=\C$ we obtain
from \ref{part:hamiltonian} that the arithmetic genus is invariant
under symplectic reduction and that all symplectic quotients of $M$
have the same arithmetic genus.  The latter fact is perhaps not very
surprising (and was known in special cases; see e.~g.\ \cite{au:to}),
since according to Guillemin and Sternberg \cite{gu:bi} all
nonsingular symplectic quotients are ``birationally equivalent'' in
the symplectic category.  For regular values $\mu$ the invariance of
the arithmetic genus under reduction was proved independently by Tian
and Zhang \cite{ti:sy}.

Applying this result to a coadjoint orbit $M=G\mu$, we find that the
arithmetic genus of $M$ is equal to 1 since the symplectic quotient
$M_\mu$ is a point.  (This observation follows also from Theorem 2.23
in \cite{sj:ho}.)  This implies that any two-dimensional symplectic
quotient of a coadjoint orbit is a sphere (possibly with orbifold
singularities).

\subsection{Moment bundles}\label{subsection:prequantum}

A $G$-equivariant line orbibundle $L$ over $M$ is called a
\emph{moment bundle} (or \emph{$G$-moment bundle\/}) if for all
components $F$ of the fixed-point set $M^T$ the orbiweight of the
$T$-action on $L|_F$ is equal to $\iota^*\Phi(F)$.  Here $\iota$
denotes the inclusion map $\lie t\to\lie g$ and $\iota^*\Phi(F)$ the
(constant) value of $\iota^*\Phi$ on $F$.  It is easy to see that this
condition is independent of the choice of the maximal torus $T$.  It
is also obvious that if $M$ admits a moment bundle $L$, then its
moment polytope is rational.  In fact, $\iota^*\Phi(F)\in
d\inv\Lambda^*$ if the generic fibre of $\ca L|_F$ is the folded line
$\C/(\Z/d\Z)$.  If $M$ is compact, a moment bundle is almost
equivariantly locally trivial at $0$ by Lemma \ref{lemma:line}, so
that $\RR(M_0,L_0)$ is well-defined.

\begin{example}\label{example:prequantum}
A \emph{prequantum line bundle} is a $G$-equivariant line orbibundle
whose equivariant Chern class is equal to the equivariant cohomology
class of the equivariant symplectic form $(\omega,\Phi)$.  The
equivariant index of a prequantum line bundle is called the
\emph{quantization\/} of $M$.  A prequantum line bundle always exists
if the cohomology class of $\omega$ is integral and $G$ is simply
connected and is then a true $G$-equivariant line bundle.  (A
nessesary and sufficient condition is the integrality of the
$G$-equivariant \emph{orbifold} cohomology class of $\omega$.)  On a
prequantum line bundle $L$ there exist a $G$-invariant Hermitian
metric and connection such that the curvature of the connection is
equal to the symplectic form and for every $\xi\in\lie g$ the
fundamental vector field $\xi_L$ is given by Kostant's formula
\begin{equation}\label{equation:kostant}
\xi_L= \operatorname{lift}\xi_M+
\langle\Phi,\xi\rangle\,\frac\partial{\partial\phi}.
\end{equation}
Here $\partial/\partial\phi$ is the generating vector field for the
scalar $S^1$-action on $L$.  This implies immediately that $L$ is a
moment bundle.  If $0$ is a quasi-regular value, then the reduced
bundle $L_0$ is a prequantum line bundle on the symplectic orbifold
$M_0$.  (See e.~g.\ \cite{gu:ge}.)  If $0$ is not quasi-regular we
shall still call the orbibundle $L_0$ a prequantum line bundle on the
stratified symplectic space $M_0$.  Note however that for $\mu\ne0$
the quotient bundles $L_\mu$ are \emph{not} prequantizing.
\end{example}

Our main result is as follows.  The proof is in Section
\ref{subsection:nonabelianmultiplicities}.

\begin{theorem}\label{theorem:moment} 
Let $M$ be a compact Hamiltonian $G$-orbifold and $L$ be a $G$-moment
bundle on $M$.  Then the multiplicity of the trivial representation in
$\RR(M,L)$ is equal to the Riemann-Roch number of the symplectic
quotient $M_0$\upn:
\begin{equation}\label{equation:redquant}
\RR(M,L)^G=\RR(M_0,L_0).
\end{equation}
In the presence of an action of a compact connected Lie group $H$ on
$M$ and $L$ that commutes with the action of $G$\upn,
\eqref{equation:redquant} holds as an equality of virtual characters
of $H$.  By Theorem \ref{theorem:independent} we also have
$\RR(M,L)^G=\RR(M_\mu,L_\mu)$ for small $\mu\in\Phi(M)$.
\end{theorem}

For prequantum line bundles this result goes by the name of
``quantization commutes with reduction'' and was conjectured (for
regular values of the moment map) by Guillemin and Sternberg
\cite{gu:ge}.  Results on the quantization conjecture were obtained in
\cite{gu:ge}, \cite{sj:ho} and \cite{br:qu2} in the context of
K\"ahler quantization.  In the above formulation the quantization
conjecture was first proved (for regular values) by Guillemin
\cite{gu:re}, Meinrenken \cite{me:on} and Vergne \cite{ve:qu} in the
abelian case and by Meinrenken \cite{me:sym} in the nonabelian case.
A similar result was obtained by Jeffrey and Kirwan \cite{je:lo2}.
For a presymplectic version see Canas et al.\ \cite{ca:ci}.  A proof
of the quantization conjecture using analytical methods was given (in
the regular case) by Tian and Zhang \cite{ti:sy}.  See \cite{me:lo}
for an application to loop group actions.  See \cite{sj:sy} for a
survey and further references.

Combined with the shifting trick Theorem \ref{theorem:moment} leads to
a complete decomposition of the virtual character $\RR(M,L)$ into
irreducible characters as follows.  Consider a moment bundle $L$ on
$M$.  Identify $\t^*$ with the subspace of $T$-fixed vectors in
$\g^*$.  Recall that every weight $\mu$ exponentiates to a character
$G_\mu\to S^1$, which gives rise to a line bundle
$$
E_\mu=G\times^{G_\mu}\C
$$
on the symplectic manifold $G\mu$.  The unique compatible invariant
almost complex structure on $G\mu$ is integrable, and by the
Borel-Weil-Bott theorem the equivariant Riemann-Roch number satisfies
\begin{equation}\label{equation:bwb}
\RR(G\mu,E_\mu)=\Ind_T^G\zeta_\mu.
\end{equation}
Here $\zeta_\mu$ is the character of $T$ defined by
$$
\zeta_\mu(\exp\xi)=\exp2\pi i\mu(\xi)
$$
and $\Ind_T^G\colon\Rep T\to\Rep G$ denotes the induction functor,
which is defined as follows.  Let $w\odot\mu=w(\mu+\rho)-\rho$
denote the affine action of the Weyl group $\eu W$, where $\rho$ is
half the sum of the positive roots.  If $\mu\in\Lambda^*$ then
$\Ind_T^G\zeta_\mu$ is nonzero if and only if there exists a Weyl
group element $w$ with $w\odot\mu\in\tplus$ and in this case
\begin{equation}\label{equation:induction}
\Ind_T^G\zeta_\mu=(-1)^{\mathrm{length}(w)}\chi_{w\odot\mu}.
\end{equation}
Let $*$ denote the involution of $\lie t$ defined by
$*\mu=\mu^*=-w_0\mu$, where $w_0$ is the longest Weyl group element.
Then for dominant $\mu$ we have by \eqref{equation:bwb} and the
K\"unneth formula
$$N(\mu)=\RR(M\times G\mu^*,L\boxtimes E_{\mu^*})^G,$$
where $N=N_L$ is the multiplicity function of $L$.  It is easy to see
that $L\boxtimes E_{\mu^*}$ is a moment bundle on the product $M\times
G\mu^*$, so we can use Theorem \ref{theorem:moment} to evaluate the
right-hand side.  According to the shifting trick, $(M\times
G\mu^*)\qu G\cong M_\mu$.

\begin{definition}\label{definition:shiftquantum}
The \emph{shifted quotient bundle} on $M_\mu$ is the orbibundle
$L_\mu\shift=(L\boxtimes E_{\mu^*})\qu G$.
\end{definition}
 
Observe that $L\shift_\mu$ is not equal to $L_\mu$ unless $\mu=0$.  In
fact, the shifted quotient bundle cannot even be defined unless $\mu$
is integral!  If $L$ is a prequantum line bundle, then $L\shift_\mu$
is a prequantum line bundle on $M_\mu$.  Theorem \ref{theorem:moment}
implies $N(\mu)=\RR(M_\mu,L\shift_\mu)$, so we have proved the
following statement.

\begin{corollary}\label{corollary:decomposition}
Let $L$ be a moment bundle on the compact Hamiltonian $G$-orbifold
$M$.  Then the decomposition of $\RR(M,L)$ into irreducible characters
is given by
$$
\RR(M,L)=
\sum_{\mu\in\Lambda^*\cap\Delta}\RR(M_\mu,L_\mu\shift)\,\chi_\mu.
$$
In particular\upn, the support of the multiplicity function is
contained in the moment polytope $\Delta$.  
\qed
\end{corollary}

While this formula may be difficult to evaluate in practice (unless
the quotients are zero- or two-dimensional), in combination with the
index theorem for orbifolds it yields interesting qualitative
information about the multiplicity diagram of a moment bundle $L$.
For instance, a weaker form of Corollary \ref{corollary:decomposition}
was used in \cite{me:on} to prove a quantum version of the
Duistermaat-Heckman Theorem.  We shall now use Corollary
\ref{corollary:decomposition} to improve on this result.

Let $M\prin$ be the principal infinitesimal orbit type stratum, i.~e.\
the set of all points at which the stabilizer has minimal dimension.
Then $\Phi$ has maximal rank on $M\prin$ and the quotients $M_\mu$ for
$\mu\in\Phi(M\prin)$ all have the same dimension, say $2k$.  Let
$\inter\Delta$ be the relative interior of the polytope $\Delta$.  The
set of \emph{generic values\/} of $\Phi$ is the set $\Delta\gen$
consisting of all $\mu\in\inter\Delta$ satisfying
$\Phi\inv(\mu)\subset M\prin$.  Let
\begin{equation}\label{equation:chamber}
\Delta\gen=\bigcup_i \Delta_i
\end{equation}
be its decomposition into connected components.  The closure of each
component $\Delta_i$ is a convex polytope (see e.~g.\ \cite{le:co}),
and $\Phi\colon\Phi\inv(\Delta_i)\to\Delta_i$ is a locally trivial
fibre bundle.  The \emph{Duistermaat-Heckman measure\/} is the measure
$\lambda_{\text{DH}}$ on $\lie t^*$ defined by
$$
\lambda_{\text{DH}}(\mu)=\vol(M_\mu)\,\lambda(\mu),
$$
where $\vol(M_\mu)$ denotes the $2k$-dimensional symplectic volume of
$M_\mu$ and $\lambda$ the normalized Lebesgue measure on the affine
subspace spanned by $\Delta$.  According to the Duistermaat-Heckman
Theorem the density function $\mu\mapsto\vol(M_\mu)$ is continuous on
$\Delta$ and is given by a polynomial on $\Delta_i$ for every $i$.

Recall that a function $f\colon \Xi\to \Z$ defined on a lattice
$\Xi\cong\Z^r$ is called \emph{quasi-polynomial\/} if there exists a
sublattice $\Xi'$ of finite index such that for all $\gamma\in\Xi$ the
translates $f_\gamma=f(\gamma+\cdot)\colon \Xi'\to \Z$ are polynomial
functions.  The degree of the $f_\gamma$ is called the \emph{degree\/}
of $f$.  If $\Xi'$ is chosen as large as possible, the number of
elements in $\Xi/\Xi'$ is the \emph{period\/} of $f$.  Consider for
instance the lattice $\Xi=\Z\times \Lambda^*$ and the function
$$f(m,\mu)=N^{(m)}(\mu),$$
where $N^{(m)}(\mu)$ is defined as the multiplicity of $\mu$ in the
character $\RR(M,L^m)$.  Replacing $\omega$ by $m\omega$ and $\Phi$ by
$m\Phi$ we obtain from Corollary \ref{corollary:decomposition}
expressions for $N^{(m)}(\mu)$ for all $m>0$ and from Theorem
\ref{theorem:rigid} for $m=0$.  (This does not work for $m<0$ because
then the almost complex structure on $M$ is not tame with respect to
$m\omega$.)  From Kawasaki's Riemann-Roch formula applied to the
orbifolds $M_\mu$ (or $\tilde{M}_\mu$ for $\mu$ that are not
quasi-regular) we then read off the following result.

\begin{corollary}[quantum DH]\label{corollary:quasipolynomial}
For every moment bundle $L$ on $M$ the function $(m,\mu)\mapsto
N^{(m)}(\mu)$ is quasi-polynomial on each of the closed cones
$$
\ca C_i=\{\,(t,t\mu):t\ge 0\mbox{ {\rm and} }\mu\in\bar{\Delta}_i\,\}.
$$
Each of these quasi-polynomials has degree $\le k$\upn, where $2k$ is
the dimension of the generic symplectic quotient\upn, and the degree
is equal to $k$ for all $i$ if $L$ is a prequantum bundle.  For
$m\ge0$ the function $m\mapsto N^{(m)}(0)$ is a quasi-polynomial\upn,
whose period is a divisor of the smallest positive integer $l$ such
that the quotient bundle $L^l\qu G$ is a genuine line bundle\upn,
i.~e.\ has fibre equal to $\C$ everywhere.  
\qed
\end{corollary}

The fact that the multiplicities $N^{(m)}(\mu)$ exhibit
quasi-polynomial behaviour even on the boundary of the cones $\ca C_i$
may be viewed as a quantum version of the continuity property of the
Duistermaat-Heckman measure at the walls of $\Delta$.

\subsection{Dual moment bundles}\label{subsection:dual}

Let $L$ be a moment bundle on $M$.  The dual orbibundle $L\inv$ is
called a \emph{dual moment bundle}.  Theorem \ref{theorem:moment}
fails for dual moment bundles.

\begin{example}[Vergne; cf.\ \cite{je:lo2}]\label{example:vergne}
Let $G=\SU(2)$.  Then $\tplus\cong i\R_+$, $\Lambda^*\cong2\pi i\Z$
and $\Lambda^*_+\cong2\pi i\N$.  Under these identifications the
positive root $\alpha=2\rho\in\Lambda^*_+$ corresponds to $4\pi i$.
Let $M$ be the projective line $\C P^1$ with $\omega$ equal to twice
the standard K\"ahler form.  The prequantum bundle on $M$ is
$L=\ca{O}(2)$, so that $H^0(M,L\inv)=\{0\}$ and $\dim H^1(M,L\inv)=1$.
It follows that $G$ acts trivially on $H^1(M,L\inv)$ and that
$\RR(M,L\inv)=\RR(M,L\inv)^G=-1$.  On the other hand,
$\Delta=\{2\rho\}$ and so $M_0$ is empty.  Thus
$\RR(M,L\inv)^G\ne\RR\bigl(M_0,L\inv_0\bigr)$.
\end{example}

It is nevertheless possible to generalize Corollaries
\ref{corollary:decomposition} and \ref{corollary:quasipolynomial} to
dual moment bundles.  As is to be expected, the correct multiplicity
formula involves some signs and shifts by half the sum of the positive
roots.

\begin{theorem}\label{theorem:dual}
Let $L$ be a moment bundle on the compact Hamiltonian $G$-orbifold
$M$.  Then
\begin{equation}\label{equation:dualmultiplicity}
\RR(M,L\inv)=(-1)^{\dim\Delta}\sum_{\mu\in\Lambda^*\cap\inter\Delta}
\RR\bigl(M_\mu,(L\shift_\mu)\inv\bigr)\Ind_T^G\zeta_{-\mu}.
\end{equation}
It follows that the support of the multiplicity function satisfies
\begin{equation}\label{equation:dualsupport}
\supp N_{L\inv}\subset
*\bigl(\inter\Delta-2(\rho-\rho_\sigma)\bigr)\cap\Lambda^*_+,
\end{equation}
where $\sigma$ is the principal wall of $M$.
\end{theorem}

Here $\Ind_T^G$ is the induction functor defined in
\eqref{equation:induction}; the principal wall of $M$ is the smallest
open wall $\sigma$ of the Weyl chamber such that
$\Delta\subset\bar\sigma$; and $\rho_\sigma$ denotes half the sum of
the positive roots of the centralizer $G_\sigma$ (so
$\rho-\rho_\sigma$ is equal to the orthogonal projection of $\rho$
onto $\sigma$).  See Section \ref{subsection:nonabelianmultiplicities}
for the proof.  Note that in contrast to Theorem \ref{theorem:moment}
the summation is only over the relative interior of the moment
polytope.  Similar formulas hold of course for all tensor powers
$L^{-m}$.  This result may be viewed as a generalization of Ehrhart's
reciprocity theorem for the number of lattice points in a convex
polytope, as we shall see in Section \ref{subsection:delzant3}.

It is not hard to see from \eqref{equation:dualsupport} that for every
$\nu\in\supp N_{L\inv}$ there exist \emph{unique\/} $w\in\eu W$ and
$\mu\in\Lambda^*_+\cap\inter\Delta$ such that $\nu=w\odot(-\mu)$.
In fact, $w=w_0w_\sigma$ and $\mu=-w_\sigma w_0\odot\nu$, where
$w_\sigma$ is the longest Weyl group element of $G_\sigma$.  (See
Lemma \ref{lemma:reflect}.)  Consequently, only the symplectic
quotient at $-w_\sigma w_0\odot\nu$ contributes to the multiplicity
at $\nu$:
$$
N_{L\inv}(\nu)=\RR\bigl(M_{-w_\sigma
w_0\odot\nu},(L\shift_{-w_\sigma w_0\odot\nu})\inv\bigr).
$$
Here are some special cases of \eqref{equation:dualsupport}: if
$\Delta$ contains strictly dominant points, then $\sigma=\inter\tplus$
and $\rho_\sigma=0$, so $\supp
N_{L\inv}\subset(*\inter\Delta-2\rho)\cap\Lambda^*_+$.  If
$\Delta=\{0\}$, then $\supp N_{L\inv}\subset\{0\}$.  If $G=\SU(3)$ and
$\sigma$ is the wall spanned by the fundamental weight $\lambda_1$,
then $\rho_\sigma=\frac1{2}\alpha_2$, so
$\rho-\rho_\sigma=\frac3{2}\lambda_1$ and $\supp
N_{L\inv}\subset(*\inter\Delta-3\lambda_2)\cap\Lambda^*_+$.  Finally
if $G$ is a torus, then $\supp N_{L\inv}\subset-\inter\Delta$.

\begin{example}
For $G=\SU(2)$ we have $\Ind_T^G\zeta_\mu=\chi_\mu$ for $\mu\ge0$,
$\Ind\zeta_{-\rho}=0$ and $\Ind_T^G\zeta_{-\mu}=-\chi_{\mu-2\rho}$ for
$\mu\ge2\rho$.  If $\dim\Delta=1$, then
$$
\RR\bigl(M,L\inv\bigr)= \sum_{\substack{\mu\ge2\rho \\
\mu\in\inter\Delta}} 
\RR\bigl(M_{\mu},(L\shift_{\mu})\inv\bigr)\,\chi_{\mu-2\rho}.
$$
For $\Delta=\{\mu\}$ we have $\RR(M,L\inv)=0$ unless $\mu=0$ or
$\mu\ge2\rho$.  If $\mu=0$, then
$\RR(M,L\inv)=-\RR\bigl(M_0,(L\shift_0)\inv\bigr)$ is a constant
character.  If $\mu\ge2\rho$, then 
$$
\RR(M,L\inv)=
-\RR\bigl(M_{\mu},(L\shift_{\mu})\inv\bigr)\,\chi_{\mu-2\rho}.
$$
In Example \ref{example:vergne} $\mu=2\rho$, $\Delta=\{2\rho\}$,
$M_{\mu}$ is a point, and $\chi_{\mu-2\rho}=\chi_0$ is the trivial
one-dimensional character, so $\RR(M,L\inv)=-1$.
\end{example}

\section{Singular symplectic quotients}\label{section:singularquotients}

In Section \ref{subsection:localmodels} we review the local normal
form theorem of \cite{sj:st} for quotients of Hamiltonian actions on
manifolds and generalize it to actions on orbifolds.  We also
investigate how orbibundles descend to orbibundles on the quotients
and, in easy cases, how nearby quotients are related to one another.
Our treatment differs from \cite{sj:st} in that we work with the
stratification of $M$ by infinitesimal orbit types (which leads to a
stratification of the symplectic quotient $M_\mu$ into orbifolds)
rather than by orbit types (which results in a stratification of
$M_\mu$ into manifolds).  This is more natural from our point of view,
because at generic levels of the moment map the symplectic quotient is
usually an orbifold, even if the original space is smooth.  Moreover,
it is impossible to remove orbifold singularities by the
desingularization process discussed in the next section.  In Section
\ref{subsection:delzant1} we apply some of our results to Delzant
spaces.

\subsection{Local normal form near a stratum}
\label{subsection:localmodels}

\subsubsection{Stratifying the level set}

For every Lie subalgebra $\h$ of $\g$ let $(\h)$ denote the conjugacy
class of $\h$.  The \emph{stratum of infinitesimal orbit type\/}
$(\h)$ is the subset of the Hamiltonian $G$-orbifold $M$ defined by
$$M_{(\h)}=\{\,m\in M:\g_m\text{ is $G$-conjugate to }\h\,\}.$$
The connected components of $M_{(\h)}$ are suborbifolds, and the set
of conjugacy classes $(\h)$ for which $M_{(\h)}$ is nonempty is
locally finite.  There is a unique conjugacy class $(\h)$ with the
property that $\h$ is subconjugate to every other stabilizer
subalgebra.  The corresponding stratum is denoted by $M\prin$ and is
called the \emph{principal infinitesimal orbit type stratum}.  It is
open, dense and connected, and it is precisely the set of points where
the moment map has maximal rank.

Let $Z$ denote the $G$-invariant subset
$$
Z=\Phi\inv(0)
$$
of $M$.  By partitioning $Z$ into sets where the dimension of the
stabilizer is constant and then partitioning further into connected
components we obtain a collection $\{Z_\alpha:\alpha\in\eu A\}$ of
locally closed subsets.  We define a partial order $\pre$ on the
indexing set $\eu A$ by putting $\alpha\pre\beta$ if and only if
$Z_\alpha\subseteq \bar{Z}_\beta$.  The decomposition
$$Z=\bigcup_{\alpha\in\eu A} Z_{\alpha}$$
is called the \emph{infinitesimal orbit type stratification} of $Z$.
Every $Z_\alpha$ arises as a connected component of some intersection
$Z\cap M_{(\h)}$ for some subalgebra $\h$ of $\g$.  We denote by
$(\g_\alpha)$ the conjugacy class of stabilizer subalgebras
corresponding to the stratum $Z_\alpha$, and by
$(G_\alpha)=(\exp\g_\alpha)$ the corresponding class of connected
subgroups.  If $\alpha\pre\beta$, then $\g_\beta$ is conjugate to a
subalgebra of $\g_\alpha$.

As we shall see, each $Z_\alpha$ is a suborbifold of $M$ and the
null-foliation of the restriction of $\omega$ to $Z_\alpha$ is given
by the $G$-orbits.  Thus the symplectic quotient $X=M_0=Z/G$ inherits
a decomposition
$$X=\bigcup_{\alpha\in\eu A}X_\alpha,$$
whose pieces are symplectic orbifolds.  It is shown in \cite{sj:st}
that $X$ has a unique open, dense and connected piece (even when $Z$
does not intersect $M\prin$).  If $0$ is a quasi-regular value (see
Definition \ref{definition:quasiregular}), then $Z=Z_\alpha$ for some
$\alpha$, so $X$ is a symplectic orbifold.

We now construct orbibundle charts for the quotient mapping $\pi\colon
Z_\alpha\to X_\alpha$.  Let $z\in Z_\alpha$ and put $x=\pi(z)$.  The
orbit $Gz$ is isotropic because $\Phi(z)=0$.  The \emph{symplectic
slice\/} at $z$ is the fibre at $z$ of the symplectic normal bundle to
$Gz$.  This is the symplectic vector orbispace $V/\Gamma$, where
$\Gamma$ is the orbifold structure group of $M$ at $z$ and
$$
V=\tilde{T}_z(Gz)^\omega/\tilde{T}_z(Gz).
$$
Here the superscript $\omega$ stands for symplectic orthogonal
complement, and $\tilde{T}_z(Gz)$ is the uniformized tangent space of
the orbit $Gz$, that is the tangent space of $\phi\inv(Gz)$ at
$\tilde{z}$ in an orbifold chart $\phi\colon\tilde U\to U$ around $z$
with $\phi(\tilde{z})=z$.  Let us choose a $G$-invariant almost
complex structure on $M$.  This induces a Hermitian structure on $V$.
The stabilizer $H=G_z$ does not necessarily act on $V$, but the
extension $\hat{H}$ of $H$ by $\Gamma$ determined by the following
commutative diagram with exact rows does:
\begin{equation}\label{equation:liftgroup}
\vcenter{\xymatrix{
{\Gamma} \ar[r]\ar@{=}[d] & {\smash{\hat
H}\vphantom{H}}\ar@{->>}[r]\ar[d] & H\ar[d]\\ 
{\Gamma}\ar[r]^-\tau & {\smash{\hat{\U}}(V)\mathstrut}\ar@{->>}[r] &
{\U}(V/\Gamma).
}}
\end{equation}
The bottom row is the definition of the ``unitary group''
$\U(V/\Gamma)$.  Here $\tau$ denotes the unitary representation of
$\Gamma$ on $V$ and $\hat{\U}(V)$ is the group of all $\phi\in\U(V)$
such that there exists a group isomorphism $f\colon\Gamma\to\Gamma$
(depending on $\phi$) satisfying
$\phi\bigl(\tau(\gamma)v\bigr)=\tau\bigl(f(\gamma)\bigr)\phi(v)$ for
all $\gamma\in\Gamma$ and $v\in V$.  Note that we do not assume
$\Gamma$ to act effectively on $V$, that is to say $\tau$ need not be
injective.  Clearly, $\hat{\U}(V)$ is a closed subgroup of the
normalizer of $\tau(\Gamma)$ inside $\U(V)$.  The action of $\hat H$
on $V$ is Hamiltonian with moment map $\Phi_V\colon V\to\hat\lie h^*$
given by $\langle\Phi_V(v),\xi\rangle= \frac1{2}\omega(\xi\cdot v,v)$.
Composing $\Phi_V$ with the natural isomorphism $\hat\lie h^*\to\lie
h^*$ we obtain a $\Gamma$-invariant map, which descends to a moment
map $\Phi_{V/\Gamma}$ for the $H$-action on $V/\Gamma$.  Let
$F(H,V/\Gamma)=G\times^H\bigl(\lie h^0\times V/\Gamma\bigr)$ be the
symplectic orbifold defined in \eqref{equation:bundle}.  (See Appendix
\ref{section:normalform}.)  The following result is a consequence of
Theorem \ref{theorem:minimalcoupling} and the isotropic embedding
theorem.

\begin{theorem}[symplectic slices, \cite{le:ha}]\label{theorem:slice} 
Assume that $\Phi(z)=0$.  A $G$-invariant neighourhood of $z$ in $M$
is isomorphic as a Hamiltonian $G$-orbifold to a neighbourhood of the
zero section in $F(H,V/\Gamma)$.  The moment map on $F(H,V/\Gamma)$ is
given by $[g,\beta,\Gamma v]\mapsto
g\bigl(\beta+\Phi_{V/\Gamma}(\Gamma v)\bigr)$.
\qed
\end{theorem}

Let the identity components of $H$ and $\hat H$ be $H^0$, resp.\ $\hat
H^0$, and let their component groups be $\pi_0(H)$, resp.\ $\pi_0(\hat
H)$.  Put $\Upsilon=\Gamma\cap\hat H^0$.  Diagram
\ref{diagram:stabilizer}, which is commutative and has exact rows and
columns, summarizes the relationships among these groups.  Each of
them depends on the point $z\in Z$.  Define the vector space $W$ by
the orthogonal splitting
\begin{equation}\label{equation:orthogonal}
V=V^{\hat H^0}\oplus W.
\end{equation}
Both summands carry a natural unitary representation of $\hat H$. 

\begin{figure}
$$
\xymatrix{
{\Upsilon}\mathstrut\ar@{((->}[r]\ar[d] &
{\Gamma}\mathstrut\ar@{->>}[r]\ar[d] &
{\Gamma}/\Upsilon\ar[d] \\
{\smash{\hat{H}^0}}\vphantom{\pi_0(\hat{H})^0}
\ar@{((->}[r]\ar@{->>}[d] & {\smash{\hat
H}}\vphantom{\pi_0(\hat{H})^0}\ar@{->>}[r]\ar@{->>}[d] &
{\smash{\pi_0(\hat{H})}}\vphantom{\pi_0(\hat{H})^0}\ar@{->>}[d] \\
{\smash{H^0}}\vphantom{\pi_0(H)^0}\ar@{((->}[r] & {\smash
H}\vphantom{\pi_0(H)^0}\ar@{->>}[r] &
{\smash{\pi_0(H)}}\vphantom{\pi_0(H)^0}.
}
$$
\caption{Stabilizers and orbifold structure groups}
\label{diagram:stabilizer}
\end{figure}

\begin{lemma}\label{lemma:w}
The origin is a weakly regular value of $\Phi$ if and only if $W=0$ at
all points in the fibre $Z=\Phi\inv(0)$.  It is a quasi-regular value
if and only if $\Phi\inv_W(0)=0$ at all points in $Z$.
\end{lemma}

\begin{proof}
Consider a point $z$ in $Z_\alpha\subset Z$.  Computing in the model
given by Theorem \ref{theorem:slice} we find
\begin{gather}
Z= G\times^H\Phi\inv_{V/\Gamma}(0)\subset G\times^H(V/\Gamma),
\label{equation:z}\\ 
Z_\alpha= G\times^H\bigl(V^{\hat H^0}/\Gamma\bigr),
\label{equation:zalpha}\\
\ker\d\Phi_z= \tilde{T}_z(Gz)^\omega/\Gamma=
\bigl(\tilde{T}_z(Gz)\oplus V\bigr)\big/\Gamma. \label{equation:dphi}
\end{gather}
Note that $V^{\hat H^0}/\Gamma$ is contained in the zero level set of
$\Phi_{V/\Gamma}$.  Weak regularity means that $Z$ is a suborbifold of
$M$ and that for all $z$ the tangent orbispace $T_zZ$ is equal to the
vector orbispace $\ker\d\Phi_z$.  By \eqref{equation:z} and
\eqref{equation:dphi}, this is equivalent to $Z$ being an open
suborbifold of $G\times^H(V/\Gamma)$, which is equivalent to
$\Phi\inv_{V/\Gamma}(0)=V/\Gamma$.  That is to say $\Phi_V=0$, which
means that $\hat H^0$ acts trivially on $V$, i.~e.\ $W=0$.

Quasi-regularity means that $Z=Z_\alpha$.  By \eqref{equation:z} and
\eqref{equation:zalpha}, this is equivalent to $\Phi\inv_V(0)=V^{\hat
H^0}$, in other words $\Phi\inv_W(0)=0$.
\end{proof}

It follows from \eqref{equation:zalpha} that the quotient
$X_\alpha=Z_\alpha/G$ is isomorphic near $x$ to the symplectic vector
orbispace $V^{\hat H^0}\big/\hat H=V^{\hat H^0}\big/\pi_0(\hat H)$.  We
conclude that $X_\alpha$ is a symplectic orbifold whose structure
group at $x$ is $\pi_0(\hat H)$ and that an orbibundle chart for the
map $Z_\alpha\to X_\alpha$ is given by
\begin{equation}\label{equation:quotient}
\vcenter{\xymatrix{
G/H^0\times V^{\hat H^0} \ar[r]\ar[d]^{/G} & G\times^H\bigl(V^{\hat
H^0}/\Gamma\bigr) \ar@{ >->}[r]\ar[d]^{/G} & Z_\alpha\ar[d]^{/G} \\ 
V^{\hat H^0}\ar[r] & V^{\hat H^0}\big/\pi_0(\hat H) \ar@{ >->}[r] &
X_\alpha.
}}
\end{equation}
Here the horizontal arrows on the left are quotient maps under the
action of $\pi_0(\hat H)$ and the horizontal arrows on the right
represent germs of equivariant embeddings at $z$ (in the top row) and
at $x$ (in the bottom row).  The fibre at $z$ is the orbit $G/H$ and
the general fibre is $G/H^0$.

\subsubsection{Neighbourhood of a stratum}

To write a normal form for a neighbourhood of a stratum $Z_\alpha$ we
examine the symplectic normal bundle $N_\alpha$ of $Z_\alpha$ in $M$.
This is the $G$-equivariant Hermitian vector orbibundle over
$Z_\alpha$ whose fibre at $z\in Z_\alpha$ is the vector orbispace
$W/\Gamma$, where $W$ is the Hermitian vector space defined by
\eqref{equation:orthogonal}.  From the symplectic slice theorem we
obtain a $G$-equivariant orbibundle chart
\begin{equation}\label{equation:normal}
\vcenter{\xymatrix{
{\bigl(G\times^{\hat H^0}W\bigr)\times V^{\hat H^0}} \ar[r]\ar[d] &
G\times^H(V/\Gamma) \ar@{ >->}[r]\ar[d] & N_\alpha\ar[d] \\
G/H^0\times V^{\hat H^0}\ar[r] & G\times^H\bigl(V^{\hat
H^0}/\Gamma\bigr) \ar@{ >->}[r] & Z_\alpha,
}}
\end{equation}
where again the leftmost horizontal maps are quotient maps under
$\pi_0(\hat H)$ and the rightmost horizontal maps are germs of
equivariant embeddings.  The fibre of $N_\alpha$ over $z$ is
$W/\Gamma$ and the fibre of the vertical map on the left is
$W/\Upsilon$ at every point.  (Note that \eqref{equation:normal} is
not a \emph{vector\/} orbibundle chart, as $W/\Upsilon$ is not a
vector space.)
Despite the fact that $\Upsilon$ depends on the point $z\in Z_\alpha$,
the vector orbispace $W/\Upsilon$ does not.  To see this, observe that
the spaces $G\times^{\hat H^0}W\cong G\times^{H^0}(W/\Upsilon)$ are
all the same because the infinitesimal representation of $\hat\lie
h\cong\lie h$ on $W$ does not depend on $z\in Z_\alpha$.

Now fix a point $z_\alpha$ in the orbifold $Z_\alpha$, say a smooth
point.  We shall decorate with a subscript $\alpha$ each of the above
groups and vector spaces evaluated at the basepoint $z_\alpha$.  Thus
$H_\alpha=G_{z_\alpha}$, $\Gamma_\alpha$ is the orbifold structure
group of $M$ at $z_\alpha$, $V_\alpha/\Gamma_\alpha$ is the symplectic
slice at $z_\alpha$, $\Upsilon_\alpha=\Gamma_\alpha\cap\hat
H_\alpha^0$, etc.  Then the conjugacy class of connected subgroups
associated to the stratum $Z_\alpha$ is $(G_\alpha)=(H_\alpha^0)$.
Moreover, the general fibres of the maps $N_\alpha\to Z_\alpha$ and
$Z_\alpha\to X_\alpha$ are $W_\alpha/\Upsilon_\alpha$, resp.\
$G/G_\alpha$.  The composition of these maps is an orbifold fibration
over $X_\alpha$.  We can easily compute the fibres and find orbibundle
charts by stacking diagram \eqref{equation:normal} on top of diagram
\eqref{equation:quotient}.  The fibre at $z$ is the associated
orbibundle $G\times^H(W/\Gamma)$ and the general fibre is
$G\times^{H^0}(W/\Upsilon)\cong
G\times^{G_\alpha}(W_\alpha/\Upsilon_\alpha)$.

We can now build a standard model $M_\alpha$ for $M$ near $Z_\alpha$
as follows.  First we construct a principal orbibundle $P_\alpha$ over
$Z_\alpha$ and then define $M_\alpha$ as an associated orbibundle.
For $z\in Z_\alpha$ let $P_{\alpha,z}$ be the set of all smooth maps
from $G\times^{G_\alpha}(W_\alpha/\Upsilon_\alpha)$ to
$G\times^H(W/\Gamma)$ that factor through $G$-equivariant Hermitian
vector orbibundle isomorphisms from
$G\times^{G_\alpha}(W_\alpha/\Upsilon_\alpha)$ to
$G\times^{H^0}(W/\Upsilon)$.  In other words,
$$
P_{\alpha,z}= \Iso\bigl(G\times^{G_\alpha}(W_\alpha/\Upsilon_\alpha),
G\times^{H^0}(W/\Upsilon)\bigr)^G\big/\pi_0(\hat H),
$$
where $\Iso$ stands for Hermitian vector orbibundle isomorphisms (that
is diffeomorphisms from one space to the other that map fibres
complex-linearly and isometrically to fibres).  By Lemma
\ref{lemma:automorphism}, $P_{\alpha,z}$ is a homogeneous space under
the group $K(G_\alpha,W_\alpha)=N_{G\times K}(G_\alpha)/G_\alpha$
defined in \eqref{equation:group}.  (Here we take $K$ to be the
unitary group $\U(W_\alpha/\Upsilon_\alpha)$, which acts on
$W_\alpha/\Upsilon_\alpha$ in a Hamiltonian fashion.)  We claim that
$P_\alpha=\coprod_{z\in Z_\alpha}P_{\alpha,z}$ is a principal
orbibundle over $X_\alpha$ with structure group
$K_\alpha=K(G_\alpha,W_\alpha)$.  Indeed, a $K_\alpha$-orbibundle
chart around $x=\pi(z)$ is given by
\begin{equation}\label{equation:principal}
\vcenter{\xymatrix{
{\hat{\ca X}}\ar[r]\ar[d]^{/K_\alpha} & {\ca X}\ar@{
>->}[r]\ar[d]^{/K_\alpha} & P_\alpha\ar[d]^{/K_\alpha} \\
V^{\hat H^0}\ar[r] & V^{\hat H^0}\big/\pi_0(\hat H) \ar@{ >->}[r] &
X_\alpha,
}}
\end{equation}
where
\begin{align*}
\ca X &= \Iso\bigl(G\times^{G_\alpha}(W_\alpha/\Upsilon_\alpha),
G\times^{H^0}(W/\Upsilon)\bigr)^G \times^{\pi_0(\hat H)}V^{\hat H^0},
\\
\hat{\ca X} &=
\Iso\bigl(G\times^{G_\alpha}(W_\alpha/\Upsilon_\alpha),
G\times^{H^0}(W/\Upsilon)\bigr)^G\times V^{\hat H^0}.
\end{align*}
By construction $Z_\alpha$ is the associated orbibundle
$P_\alpha\times^{K_\alpha}(G/G_\alpha)$.  Choose a principal
connection on $P_\alpha$ and let $F_\alpha$ be the Hamiltonian
$G\times K_\alpha$-orbifold $F(G_\alpha,W_\alpha/\Upsilon_\alpha)=
G\times^{G_\alpha}(\g_\alpha^0\times W_\alpha/\Upsilon_\alpha)$
defined in \eqref{equation:bundle}.  Our standard model is the
associated orbibundle
\begin{equation}\label{equation:model}
M_\alpha=P_\alpha\times^{K_\alpha}F_\alpha\cong
P_\alpha\times^{K_\alpha}\bigl(G\times^{G_\alpha}(\g_\alpha^0\times
W_\alpha/\Upsilon_\alpha)\bigr).
\end{equation}
By Theorem \ref{theorem:minimalcoupling} the minimal coupling form is
a closed two-form on $M_\alpha$ and is nondegenerate in a
neighbourhood of $Z_\alpha\subset M_\alpha$. By construction the
symplectic normal bundle of $Z_\alpha$ in $M_\alpha$ is isomorphic to
$N_\alpha$.  The constant rank embedding theorem implies the following
result.

\begin{theorem}[cf.\ \cite{sj:st}]\label{theorem:model}
There exist a $G$-invariant open neighbourhood $U_\alpha$ of
$Z_\alpha$ in $M_\alpha$ and an isomorphism of Hamiltonian
$G$-orbifolds $U_\alpha\to M$ onto a neighbourhood of $Z_\alpha$ in
$M$.  
\qed
\end{theorem} 

\subsubsection{Transverse structure of the singularities}

Taking symplectic quotients and identifying $F_\alpha\qu G$ with
$(W_\alpha/\Upsilon_\alpha)\qu G_\alpha$ as in Example
\ref{example:model} we obtain a local model for the quotient $X$ near
its stratum $X_\alpha$:
\begin{equation}\label{equation:quotientmodel}
M_\alpha\qu G\cong P_\alpha\times^{K_\alpha} (F_\alpha\qu G)\cong
P_\alpha\times^{K_\alpha}(W_\alpha/\Upsilon_\alpha)\qu G_\alpha,
\end{equation}
which fibres over the stratum $X_\alpha=P_\alpha/K_\alpha$ with
general fibre $(W_\alpha/\Upsilon_\alpha)\qu G_\alpha$.  It is
instructive to do this calculation directly by exhibiting $Z$ (locally
near $Z_\alpha$) as a bundle over $Z_\alpha$.  Consider the
suborbibundle $S_\alpha$ of $M_\alpha$ defined by
\begin{equation}\label{equation:minimalsymplectic}
S_\alpha=P_\alpha\times^{K_\alpha}T^*(G/G_\alpha)\cong
P_\alpha\times^{K_\alpha}\bigl(G\times^{G_\alpha}\g_\alpha^0\bigr).
\end{equation}
This space has the following properties: $S_\alpha\cap U_\alpha$ is
symplectic, $Z_\alpha$ is coisotropic in $S_\alpha$, and in fact $0$
is a regular value of $\Phi|_{S_\alpha}$ and $Z_\alpha$ is its zero
fibre.  Furthermore, the standard model $M_\alpha$ is a symplectic
vector orbibundle over $S_\alpha$ with general fibre
$W_\alpha/\Upsilon_\alpha$, and the projection $M_\alpha\to S_\alpha$
maps $Z$ onto $Z_\alpha$.  In other words, the following diagram
commutes:
\begin{equation}\label{equation:levelbundle}
\vcenter{\xymatrix{
Z\ar@{ >->}[r]\ar@{ >->>}[d] & M_\alpha\ar@{->>}[d] \\
**[l]Z\cap S_\alpha=Z_\alpha\ar@{((->}[r] & S_\alpha.
}}
\end{equation}
Upon dividing out the left-hand column by the action of $G$ we obtain
the fibration $M_\alpha\qu G\to X_\alpha$ of
\eqref{equation:quotientmodel}.

\begin{theorem}[cf.\ \cite{sj:st}]\label{theorem:quotientmodel}
A neighbourhood of $X_\alpha$ in $X$ is modelled by a neighbourhood of
$X_\alpha$ in the fibre bundle
$$
P_\alpha\times^{K_\alpha}(W_\alpha/\Upsilon_\alpha)\qu
G_\alpha\longrightarrow X_\alpha,
$$
whose general fibre is the symplectic cone
$(W_\alpha/\Upsilon_\alpha)\qu G_\alpha$.  
\qed
\end{theorem}

The scalar $S^1$-action on $W_\alpha$, which is generated by the
function $w\mapsto-\frac1{2}\lVert w\rVert^2$, induces a Hamiltonian
$S^1$-action on the symplectic cone $(W_\alpha/\Upsilon_\alpha)\qu
G_\alpha$.  The base of the cone is the level set at level $1$, which
is called the \emph{link\/} of the stratum:
$$
\bigl((\Phi_{W_\alpha}\inv(0)\cap
S^{2r_\alpha-1})/\Upsilon_\alpha\bigr)\big/G_\alpha,
$$
where $2r_\alpha=\dim W_\alpha$ and $S^{2r_\alpha-1}$ is the unit
sphere in $W_\alpha$.  The symplectic quotient at level $1$ is the
\emph{symplectic link\/} $(\bb P(W_\alpha)/\Upsilon_\alpha)\qu
G_\alpha$.  The quotient map from the link to the symplectic link is a
principal $S^1$-orbibundle.  Thus a singular symplectic quotient is
locally an iterated cone, where the base of each cone is an
$S^1$-orbibundle over a symplectic quotient of lower depth.

\subsection{Symplectic cross-sections}\label{subsection:crosssection}

Results similar to Theorems \ref{theorem:model} and
\ref{theorem:quotientmodel} hold at arbitrary levels of the moment
map.  One way to see this is to invoke the shifting trick.  A better
way is to appeal to Guillemin and Sternberg's symplectic cross-section
theorem.  The version we shall discuss is borrowed from \cite{gu:sy}
and \cite{le:co}.

Recall that every coadjoint orbit $G\mu$ of $G$ intersects the
positive Weyl chamber $\t^*_+$ in exactly one point, which we denote
by $q(\mu)$.  Then $q$ is a continuous quotient mapping for the
coadjoint action,
\begin{equation}\label{equation:q}
q\colon\lie g^*\longrightarrow\tplus\sim\lie g^*/G.
\end{equation}
Let $A$ be the identity component of the centre of $G$ and $[G,G]$ the
commutator subgroup or semisimple part.  Then the intersection of $A$
and $[G,G]$ is finite and $G=A[G,G]$.  Let $\g=\a\oplus[\g,\g]$ be the
corresponding decomposition of the Lie algebra.  All points in an open
wall $\sigma$ of the positive Weyl chamber have the same centralizer
$G_\sigma$.  The group $G_\sigma$ is connected and the stratum in
$\lie g^*$ of orbit type $G_\sigma$ is the saturation
$G\sigma=q\inv(\sigma)$ of $\sigma$.  We define a partial order $\pre$
on the open walls by putting $\tau\pre\sigma$ if $\tau$ is contained
in the closure of $\sigma$.  If $\tau\pre\sigma$ then $G_\sigma$ is a
subgroup of $G_\tau$, and if $\sigma$ is the top-dimensional open wall
$\inter\tplus$, then $G_\sigma$ is the maximal torus $T$.  Therefore
$T\subset G_\sigma$ for all $\sigma$.  We can write
$G_\sigma=A_\sigma[G_\sigma,G_\sigma]$, where $A_\sigma$ is the
identity component of the centre of $G_\sigma$.  Let
$\g_\sigma=\a_\sigma\oplus[\lie g_\sigma,\lie g_\sigma]$ be the
corresponding $G_\sigma$-invariant splitting of the Lie algebra.  We
shall identify $\lie g_\sigma^*$ with the subspace of $\lie g^*$
centralized by $A_\sigma$ and $\lie a_\sigma^*$ with the subspace of
$\lie g^*$ centralized by $G_\sigma$.  Then $\lie g^*$ is an
$A_\sigma$-invariant direct sum $\lie g^*=\lie a_\sigma^*\oplus[\lie
g_\sigma,\lie g_\sigma]^*\oplus\lie g_\sigma^0$.  The summand
$\a_\sigma^*$ is equal to the linear span of the wall $\sigma$.
Define
\begin{equation}\label{equation:star}
\eu S_\sigma=G_\sigma\cdot\star\sigma
\end{equation}
where $\star\sigma$ is the open star $\bigcup_{\tau\suc\sigma}\tau$ of
$\sigma$.  Then $\eu S_\sigma$ is a $G_\sigma$-invariant open
neighbourhood of $\sigma$ in $\g_\sigma^*$ and it is in fact a slice
for the coadjoint action at all points of $\sigma$.  For example, if
$\sigma=\a^*$ then $\eu S_\sigma=\g^*$, and if $\sigma=\inter\tplus$
then $\eu S_\sigma=\sigma$.  The \emph{symplectic cross-section\/} of
$M$ over $\sigma$ is the subset
$$Y_\sigma=\Phi\inv(\eu S_\sigma).$$
Note that $\Phi(Y_\sigma)\subset\g_\sigma^*$.  Let $M_\sigma$ denote
the $G$-invariant open subset $GY_\sigma$ of $M$.

\begin{theorem}[symplectic cross-sections]\label{theorem:crosssection} 
For every open wall $\sigma$ of $\tplus$ the symplectic cross-section
$Y_\sigma$ is a connected $G_\sigma$-invariant symplectic suborbifold
of $M$.  The action map $G\times Y_\sigma\to M$ induces a
diffeomorphism $G\times^{G_\sigma}Y_\sigma\to M_\sigma$.  If
$Y_\sigma$ is nonempty\upn, then its saturation $M_\sigma$ is open and
dense in $M$.  The $G_\sigma$-action on $Y_\sigma$ is Hamiltonian with
moment map $\Phi|_{Y_\sigma}$.  For all $\mu\in\eu S_\sigma$ the
inclusion map $\Phi\inv(\mu)\hookrightarrow Y_\sigma$ induces a
symplectomorphism $M_\mu\cong(Y_\sigma)_\mu$.
\end{theorem}

The \emph{principal wall\/} for $M$ is the minimal open wall $\sigma$
of the Weyl chamber such that $\Delta\subset\bar{\sigma}$.  The
cross-section $Y_\sigma$ over it is the \emph{principal
cross-section}.  The action of $[G_\sigma,G_\sigma]$ on the
principal cross-section is trivial, so that it is in effect a
Hamiltonian $A_\sigma$-space.

Consider the symplectic manifold $T^*G\cong G\times\lie g^*$ on which
$G$ acts by left multiplication.  The moment map is projection on the
second factor and so the cross-section over $\sigma$ is the $G\times
G_\sigma$-manifold $G\times\eu S_\sigma$.  Theorem
\ref{theorem:crosssection} tells us that $G\times\eu S_\sigma$ is a
symplectic submanifold of $T^*G$.  (The symplectic structure can also
be produced by minimal coupling.)  It is the \emph{universal
$\sigma$-cross-section\/} in the sense that the cross-section over
$\sigma$ of the Hamiltonian $G$-orbifold $M$ is isomorphic in a
natural way to the symplectic quotient
$$
Y_\sigma\cong\bigl(M\times(G\times\eu S_\sigma)^-\bigr)\bigqu G,
$$
where the superscript ``$-$'' means that we replace the symplectic
form with its opposite.  Furthermore, the bundle
$G\times^{G_\sigma}Y_\sigma$ is diffeomorpic to $(G\times\eu
S_\sigma\times Y_\sigma)\qu G_\sigma$ and so acquires a symplectic
form.  Theorem \ref{theorem:crosssection} can now be supplemented as
follows.

\begin{addendum}\label{addendum:crosssection}
For all $\sigma$ the open embedding $G\times^{G_\sigma}Y_\sigma\to
GY_\sigma\subset M$ is symplectic.
\end{addendum}

Now let $\mu$ be an arbitrary point in $\tplus$.  How to stratify the
fibre $\Phi\inv(\mu)$?  First assume that $\mu$ is fixed under the
coadjoint action.  This case can be reduced to the case $\mu=0$ by
replacing $\Phi$ by the equivariant moment map $\Phi-\mu$.  Now
consider an arbitrary point $\mu$ and let $\sigma$ be the open face of
the positive Weyl chamber containing $\mu$.  Observe that
$\Phi\inv(\mu)\subset Y_\sigma$ and that $G_\sigma$ fixes $\mu$.  By
the results of the previous section the intersections of
$\Phi\inv(\mu)$ with infinitesimal $G_\sigma$-orbit type strata in
$Y_\sigma$ are orbifolds and their quotients by $G_\sigma$ are
symplectic orbifolds.  (We can also stratify the $G$-invariant subset
$\Phi\inv(G\mu)$ into orbifolds by flowing out the strata of
$\Phi\inv(\mu)$ by the $G$-action.)  Therefore, the symplectic
quotient $M_\mu=(Y_\sigma)_\mu$ has a natural decomposition into
symplectic orbifolds and structure theorems analogous to
\ref{theorem:model} and \ref{theorem:quotientmodel} hold.

\subsection{Orbibundles and nearby quotients}
\label{subsection:linenearby}

We give some applications of Theorem \ref{theorem:model}, namely to
quotients of line orbibundles and to the variation of symplectic
quotients.  We also discuss a property of moment and rigid bundles.

Let $L$ be a $G$-equivariant line orbibundle on $M$ and let $z$ be a
point in $Z=\Phi\inv(0)$.  Let $\Gamma$ be the structure group of $z$
and let $\phi\colon\tilde U\to U$ be an orbifold chart around $z$,
where $\tilde U$ is a $\Gamma$-invariant ball about the origin in
$\tilde T_zM$ and $\phi(0)=z$.  Let $\tilde L$ be a line bundle on
$\tilde U$ such that $L|_U=\tilde L/\Gamma$.  We may assume that $U$
is invariant under the action of the stabilizer $H=G_z$.  As in
\eqref{equation:liftgroup} the action of $H$ lifts to actions of
groups $\hat H$ on $\tilde U$ and $\hat H^L$ on $\tilde L$.  The
extensions $\hat H$ and $\hat H^L$ are not necessarily isomorphic
(unless $\Gamma$ acts effectively), but the natural homomorphism
$\beta\colon\hat H^L\to\hat H$ restricts to a covering homomorphism
$(\hat H^L)^0\to\hat H^0$, so that we may identify the Lie algebras
$\hat\lie h^L$, $\hat\lie h$ and $\lie h$.

\begin{lemma}\label{lemma:aelt}
The line orbibundle $L$ is almost equivariantly locally trivial at $z$
if and only if $\tilde L$ has an $\lie h$-invariant section that does
not vanish at $0$.  If either of these conditions holds\upn, then $L$
is almost equivariantly locally trivial in a neighbourhood of $z$.
\end{lemma}

\begin{proof}
If $L$ is almost equivariantly locally trivial at $z$ (see Definition
\ref{definition:aelt}), then $H^0$ acts trivially on $L_z$, so that
$(\hat H^L)^0$ acts trivially on $\tilde L_0$.  We can then produce an
$\lie h$-invariant section $\tilde s$ of $\tilde L$ that does not
vanish at $0$ by starting with an arbitrary section $s$ such that
$s(0)\ne0$ and averaging:
$$\tilde s(y)=\int_{(\hat H^L)^0}gs\bigl((\beta g)\inv y\bigr)\d g.$$
Conversely, if $\tilde L$ has an $\lie h$-invariant section that does
not vanish at $0$, then $(\hat H^L)^0$ acts trivially on the fibre
$\tilde L_0$.  It follows that $H^0$ acts trivially on $L_z$.

By the same token, if $\tilde L$ has an $\lie h$-invariant section
that does not vanish at $0$, then for all $v$ near $0$ the identity
component of $\beta\inv(\hat H^0_v)$ acts trivially on the fibre
$\tilde L_v$.  This implies that $G_y^0$ acts trivially on $L_y$ for
all $y$ near $z$.
\end{proof}

Now assume that $L$ is almost equivariantly locally trivial at all
points of $Z$.  We want to construct orbibundle charts for the
quotient $L_0=(L|_Z)/G$.  Let $V/\Gamma$ be the symplectic slice at
$z\in Z$ and let $\Phi_V$ be the moment map for the $\hat H$-action on
$V$.  By the symplectic slice theorem we can identify a neighbourhood
of $z$ in $Z$ with a neighbourhood of $0$ in $\Phi_{V/\Gamma}\inv(0)$,
and a neighbourhood of $x=\pi(z)$ in $M_0=Z/G$ with a neighbourhood of
the apex of the symplectic cone $(V/\Gamma)\qu H$.  Let $E$ be the
restriction of $L$ to $\Phi_{V/\Gamma}\inv(0)$ (with respect to the
chosen embedding); then near $x$ the quotient $L_0$ can be identified
with the quotient $E/H$.  Choose a nowhere vanishing $\lie
h$-invariant section $s$ of $\tilde L$ over $\tilde U$.  Then $s$
induces a nowhere vanishing section of the induced map $\tilde E/(\hat
H^L)^0\to\Phi_V\inv(0)/\hat H^0$, where $\tilde E$ is the restriction
of $\tilde L$ to $\Phi_V\inv(0)$.  In other words, $\tilde E/(\hat
H^L)^0$ is a trivial complex line bundle over the space
$\Phi_V\inv(0)/\hat H^0$.  Dividing out by the action of $\pi_0(\hat
H)$ on $\Phi_V\inv(0)/\hat H^0$ and by $\pi_0(\hat H^L)$ on $\tilde
E/(\hat H^L)^0$ we obtain a commutative diagram
\begin{equation}\label{equation:quotientorbi}
\vcenter{\xymatrix{
{\tilde E/(\hat H^L)^0} \ar[r]^-{/\pi_0(\hat H^L)}\ar[d] & E/H\ar@{
>->}[r]\ar[d] & L_0\ar[d] \\
{\Phi_V\inv(0)/\hat H^0}\ar[r]^-{/\pi_0(\hat H)} &
{\Phi_{V/\Gamma}\inv(0)/H}\ar@{ >->}[r] & M_0,
}}
\end{equation}
where the horizontal maps on the right denote germs of embeddings.  If
$0$ is a quasi-regular value, then $\Phi_V\inv(0)/\hat H^0$ is a
symplectic manifold and \eqref{equation:quotientorbi} shows that
$L_0$ is a line orbibundle on the symplectic orbifold $M_0$.  If
$0$ is not quasi-regular, we consider an atlas on $L_0$ consisting
of diagrams as in \eqref{equation:quotientorbi} to be the
\emph{definition\/} of an orbibundle on the singular quotient $M_0$.
Because of the symplectic cross-section theorem these observations
generalize immediately to arbitrary values of the moment map.

\begin{proposition}
Let $\mu$ be an arbitrary value of $\Phi$.  If $L$ is a line
orbibundle on $M$ that is almost equivariantly locally trivial at
level $\mu$\upn, then the quotient $L_\nu$ is a line orbibundle on
$M_\nu$ for all $\nu$ in a neighbourhood of $\mu$.  
\qed
\end{proposition}

Henceforth assume that $\mu$ is a quasi-regular value, so that the
level set $\Phi\inv(\mu)$ consists of one single stratum and
$M_\mu=\Phi\inv(\mu)/G$ is an orbifold.  Assume also for the moment
that $\mu=0$.  Then $Z=Z_\alpha$ for some $\alpha\in\eu A$, so for all
$\nu$ sufficiently close to $0$ the fibre $\Phi\inv(\nu)$ is contained
in the standard neighbourhood $U_\alpha$ of Theorem
\ref{theorem:model}.  We conclude that the quotient $M_\nu$ is
symplectomorphic to the associated orbibundle
$P_\alpha\times^{K_\alpha}(F_\alpha)_\nu$, whose general fibre is the
(possibly singular) quotient $(F_\alpha)_\nu$.  Furthermore, if $L$ is
a $G$-equivariant line orbibundle on $M$, then the restriction of $L$
to $U_\alpha$ is equivariantly isomorphic to the pullback of $L|_Z$ to
$U_\alpha$, because $U_\alpha$ retracts equivariantly onto $Z$.
Consequently, if $L$ is almost equivariantly locally trivial at level
$0$, the orbibundle $L_\nu$ is isomorphic to the pullback of the
orbibundle $L_0$ along the map $M_\nu\to M_0$.  These statements are
true at all quasi-regular values of the moment map.

\begin{proposition}\label{proposition:varquot} 
Let $\mu$ be a quasi-regular value of $\Phi$.  For all $\nu\in \g^*$
sufficiently close to $\mu$ there is a symplectic fibre orbibundle
$M_\nu\to M_\mu$ whose general fibre is the symplectic quotient
$(F_\alpha)_\nu$\upn, where $\alpha$ indicates the infinitesimal orbit
type of the single stratum intersecting $\Phi\inv(\mu)$.  If $\nu$ is
also a quasi-regular value\upn, then the fibre $(F_\alpha)_\nu$ is an
orbifold.  If $L$ is a line orbibundle on $M$ which is almost
equivariantly locally trivial at level $\mu$\upn, then the line
orbibundle $L_\nu$ is isomorphic to the pullback of the line
orbibundle $L_\mu$ via the map $M_\nu\to M_\mu$.  
\qed
\end{proposition}

We conclude this section with two facts that were referred to in
Section \ref{subsection:preliminaries}.

\begin{proposition}\label{proposition:weakquasi}
Weakly regular values of $\Phi$ are quasi-regular.
\end{proposition}

It is easy to see that the converse of this statement is false.  For
instance, if $M$ is the complex line with the standard symplectic form
and the standard circle action, which is generated by the Hamiltonian
function $z\mapsto-\frac1{2}\lvert z\rvert^2$, then $0$ is a
quasi-regular, but not weakly regular, value.

\begin{proof}
Obvious from Lemma \ref{lemma:w} for those values of $\Phi$ that are
$G$-fixed.  The general case now follows from the symplectic
cross-section theorem.
\end{proof}

\begin{lemma}\label{lemma:line}
Assume that the Hamiltonian $G$-orbifold $M$ is compact.
\begin{enumerate}
\item\label{part:rigid}
If the orbibundle $L$ is rigid\upn, then it is almost equivariantly
locally trivial everywhere.
\item\label{part:moment}
Assume that $L$ is a moment bundle.  For $m$ in $M$ let
$\sigma_m\in\lie g_m^*$ denote the character defining the
$(G_m)^0$-action on $L_m$.  Then $\sigma_m$ is equal to the projection
of $\Phi(m)\in\lie g^*$ onto $\lie g_m^*$.  Hence $L$ is almost
equivariantly locally trivial at level $0$.
\item\label{part:commutator}
If $L$ is a moment bundle\upn, then for all $m\in M$ the action of the
identity component of $G_m\cap\bigl[G,G_{\Phi(m)}\bigr]$ on $L_m$ is
trivial.
\end{enumerate}
\end{lemma}

\begin{proof}
Let $m$ be an arbitrary point in $M$.  Let $H$ be a maximal torus of
$(G_m)^0$ contained in the maximal torus $T$ of $G$.

To prove \ref{part:rigid} it suffices to show that $H$ acts trivially
on $L_m$.  Let $F$ be the connected component of $M^H$ that contains
$m$.  Then the orbiweight of the $H$-action on $L_{m'}$ is the same
for all $m'$ in $F$, so it suffices to show that $H$ acts trivially on
$L_{m'}$ for some $m'$ in $F$.  Now $F$ is a closed, and hence
compact, $T$-invariant symplectic suborbifold of $M$, so the
fixed-point set $F^T$ is nonempty.  Let $m'\in F^T$; then, because $L$
is rigid, $T$ acts trivially on $L_{m'}$, and \emph{a fortiori\/} so
does $H$.

For the proof of \ref{part:moment} let $\iota_m$ denote the inclusion
of $\lie g_m$ into $\lie g$.  Note that $\iota_m^*\Phi(m)$ is a
character of $\lie g_m$, so it suffices to prove that
$\langle\iota_m^*\Phi(m),\xi\rangle= \sigma_m(\xi)$ for all $\xi$ in
the Lie algebra of $H$.  Take $m'\in F^T$ as above.  Since
$\langle\Phi,\xi\rangle$ is constant on $F$,
$\langle\iota_m^*\Phi(m),\xi\rangle= \langle\Phi(m),\xi\rangle=
\langle\Phi(m'),\xi\rangle= \sigma_{m'}(\xi)= \sigma_m(\xi)$, where
the third equality follows from the definition of a moment bundle.

Finally, \ref{part:commutator} follows from \ref{part:moment} and the
fact that $\bigl\langle\Phi(m),[\xi,\eta]\bigr\rangle=0$ for all
$\xi\in\lie g$ and $\eta\in\lie g_{\Phi(m)}$.
\end{proof}

\subsection{Delzant spaces I}\label{subsection:delzant1}

Delzant spaces are the symplectic counterparts of projective toric
varieties.  They provide an example of the theory expounded above and
play a role in multiple symplectic cutting, which will be discussed in
Section \ref{subsection:hamiltonian}.  See Sections
\ref{subsection:delzant2} and \ref{subsection:delzant3} for further
results on Delzant spaces.  A \emph{Delzant space\/} is a
multiplicity-free connected Hamiltonian $H$-space, where $H$ is a
torus.  (Recall that a Hamiltonian $H$-space is
\emph{multiplicity-free\/} if each of its symplectic quotients is
either a point or empty.  By a space we mean for the present purposes
an orbifold or a stratified space that arises as a symplectic quotient
of an orbifold.)  Delzant \cite{de:ha} proved that multiplicity-free
Hamiltonian $H$-manifolds are completely characterized by their moment
polytopes and showed how to reconstruct such a manifold from its
polytope.  His results were extended to orbifolds by Lerman and Tolman
\cite{le:ha}.  We shall adapt their construction to polyhedra that are
not necessarily simple or compact to obtain a larger class of spaces,
although in this more general setting it is not known if the
(labelled) polyhedron determines the space.  Our version of the
construction is based on symplectic cutting.

\subsubsection{Symplectic cutting}\label{subsubsection:cut}

Suppose that $S^1$ acts on $M$ in a Hamiltonian fashion with moment
map $\psi\colon M\to \R$ and that the $S^1$-action commutes with the
$G$-action.  Consider the diagonal action of $S^1$ on the product
$M\times\C$, which has moment map
$\tilde\psi(m,z)=\psi(m)-\frac1{2}\vert z\vert^2$.  Here $\C$ is the
complex line equipped with the standard cirle action and the standard
symplectic structure.  In \cite{le:sy2} Lerman defines the
\emph{symplectic cut\/} as the symplectic quotient at level $0$,
$$M_{\ge0}=(M\times\C)\qu S^1.$$
Let $M_{>0}$ be the set of points $m\in M$ with $\psi(m)>0$ and define
$\beta\colon M_{>0}\to M\times\C$ by
\begin{equation}\label{equation:beta}
\beta(m)=\bigl(m,\sqrt{2\,\psi(m)}\,\bigr).
\end{equation} 
The zero level set of $\tilde\psi$ is the disjoint union of
$\psi\inv(0)\times\{0\}$ and the image of the map $M_{>0}\times S^1\to
M\times\C$ sending the pair $(m,e^{i\theta})$ to
$e^{i\theta}\cdot\beta(m)$.  This implies that as a set $M_{\ge0}$ is
the disjoint union of a copy of $M_{>0}$ and a copy of the reduced
space $M_0=\psi\inv(0)/S^1$.  It implies also that if $0$ is a regular
value of $\psi$, then it is a regular value of $\tilde\psi$.

\begin{proposition}[\cite{le:sy2}]\label{proposition:cut}
The canonical embeddings
\begin{equation}\label{equation:isomorphism}
\iota_0\colon M_0\longhookrightarrow M_{\ge0},\qquad\iota_{>0}\colon
M_{>0}\longhookrightarrow M_{\ge0}
\end{equation}
are symplectic embeddings.  If $0$ is a regular value of $\psi$\upn,
then $M_{\ge0}$ is a symplectic orbifold.  The lifted $G$-action
$g(m,z)=(gm,z)$ on $M\times\C$ induces a Hamiltonian $G$-action on
$M_{\ge0}$\upn, whose moment map is the map $\Phi_{\ge0}$ induced by
the $S^1$-invariant map $(m,e^{i\theta})\mapsto\Phi(m)$.  In
particular\upn, the original $S^1$-action on $M$ induces a Hamiltonian
$S^1$-action on $M_{\ge0}$.  The image of its moment map $\psi_{\ge0}$
is $\psi(M)\cap\R_{\ge0}$.  
\qed
\end{proposition}

Performing these constructions with the product $M\times\C^-$ instead
of $M\times\C$ (where $\C^-$ is the complex line endowed with the
standard circle action but the opposite of the standard symplectic
form), and the moment map $\psi(m)+\frac1{2}\vert z\vert^2$ instead of
$\tilde{\psi}$, we obtain the opposite cut space $M_{\le0}$, which is
a Hamiltonian $S^1$-space whose moment map $\psi_{\le0}$ has image
$\psi(M)\cap\R_{\le0}$.

Now let $L$ be a $G\times S^1$-equivariant line orbibundle on $M$.
Assume that $L$ is almost equivariantly locally trivial at level $0$
with respect to the $S^1$-action.  The \emph{cut bundle} is the
$G$-equivariant orbibundle on the stratified symplectic space
$M_{\ge0}$ defined by
$$L_{\ge0}=(\pr_M^*L)\qu S^1.$$
where $\pr_M\colon M\times\C\to M$ is the projection onto the first
factor.  There are canonical isomorphisms (see \cite{me:sym})
\begin{equation}\label{equation:isoline}
\iota_0^*L_{\ge0}\cong L_0,\qquad \iota_{>0}^*L_{\ge0}\cong
L|_{M_{>0}}.
\end{equation}

In a similar vein we can define a cut bundle $L_{\le0}$ on the
opposite cut space $M_{\le0}$.

From \eqref{equation:isomorphism} and \eqref{equation:isoline} we see
that cutting is a local operation in the following sense.  Suppose
that the circle action is defined only on a $G$-invariant open
neighbourhood $U$ of $Z$, where $Z$ is a closed suborbifold of $M$ of
codimension one such that $M-Z$ consists of two connected components,
and that the action on $U$ is generated by a $G$-invariant Hamiltonian
function $\psi$ with $\psi\inv(0)=Z$.  The data $(U,Z,\psi)$ are
called the \emph{cutting data}.  The cut $U_{\ge0}$ is well-defined
and, letting $M_{>0}$ be the component of $M-Z$ containing $U_{>0}$,
we can construct a global symplectic cut
\begin{equation}\label{equation:localcut}
M_{\ge0}=M_{>0}\amalg_{U_{>0}}U_{\ge0}
\end{equation}
by pasting $M_{>0}$ to $U_{\ge0}$ along the open symplectic embeddings
$U_{>0}\hookrightarrow M_{>0}$ and $\iota_{\ge0}\colon
U_{>0}\hookrightarrow U_{\ge0}$.  The resulting orbifold depends only
on $M_{>0}$, $Z$, and the germs at $Z$ of the set $U$ and the function
$\psi$.  The opposite cut $M_{\le0}$ and the cut bundles $L_{\ge0}$
and $L_{\le0}$ are defined likewise.

\subsubsection{Cutting the cotangent bundle}

Let $H$ be a $k$-dimensional torus and let $\Lambda=\ker\exp$ be its
integral lattice.  A \emph{label\/} is an ordered pair consisting of a
nonzero lattice vector and an arbitrary real number.  Let $\S$ be a
set of labels,
$$
\S=\bigl\{(v_1,r_1),(v_2,r_2),\dots,(v_n,r_n)\bigr\}.
$$
The polyhedron $\P$ \emph{associated to\/} $\S$ is the subset of
$\h^*$ consisting of all points $\mu$ satisfying the inequalities
$$
\langle\mu,v_i\rangle\ge r_i\qquad\text{for $i=1$, $2,\dots$, $n$}.
$$
The pair $(\S,\P)$ is called a \emph{labelled polyhedron}.  Clearly
$\P$ is a convex polyhedron.  It is not necessarily compact or
nonempty, but it has finitely many faces.  For every collection of
labelling vectors $v_1,\dots$, $v_n$ there is a nonempty open set of
parameters $(r_1,\dots,r_n)$ such that the associated polyhedron is
nonempty and $k$-dimensional.  If the $r_i$ are rational, then $\P$ is
rational.

Now let $A\colon\R^n\to\h$ be the linear map sending the $i$-th
standard basis vector $e_i$ in $\R^n$ to the $i$-th labelling vector
$v_i$.  Since $A$ maps the standard lattice $\Z^n$ to $\Lambda$, it
exponentiates to a homomorphism $\bar{A}\colon T^n\to H$, where $T^n$
is the torus $\R^n/\Z^n$.
%
The left $H$-action on the cotangent bundle $T^*H\cong H\times\lie
h^*$ is Hamiltonian with moment map $\pr_{\lie h^*}$.  Via the
homomorphism $\bar{A}$, $T^n$ also acts on $T^*H$ in a Hamiltonian
fashion.  Then $T^*H$ inherits a $T^n$-moment map from $H$, which we
translate by the vector $r=(r_1,\dots,r_n)\in(\R^n)^*$ to obtain the
moment map $\psi_r\colon T^*H\to(\R^n)^*$ given by
\begin{equation}\label{equation:psi}
\psi_r=A^*\circ\pr_{\lie h^*}-r.
\end{equation}
Here $A^*\colon \h^*\to (\R^n)^*$ is the transpose map of $A$.  The
moment map for the $i$-th circle in $T^n$ sends $(h,\eta)$ to
$(A^*\eta)_i-r_i$, where $h\in H$, $\eta\in\lie h^*$, and
$(A^*\eta)_i$ is the $i$-th coordinate of the vector
$A^*\eta\in(\R^n)^*$.  Recall that the moment map for the standard
$T^n$-action on $\C^n$ is given by $\phi(z)=-\frac1{2}\bigl(\vert
z_1\vert^2,\dots,\vert z_n\vert^2\bigr)$.

\begin{definition}\label{definition:delzant}
The \emph{Delzant space $D_\S$ labelled by $\S$} is the symplectic
quotient at level $0$ of the product $T^*H\times\C^n$ with respect to
the diagonal $T^n$-action and the moment map
\begin{equation}\label{equation:tildepsi}
\tilde{\psi}_r(h,\eta,z)= \psi_r(h,\eta)+\phi(z).
\end{equation}
\end{definition} 

Alternatively, $D_\S$ can be thought of as the stratified symplectic
space obtained by performing successive symplectic cuts on $T^*H$ with
respect to each of the $n$ circles in $T^n$.

Since the $H$-action on $T^*H$ commutes with the $T^n$-action, $D_\S$
is a Hamiltonian $H$-space with moment map $\Psi_\S\colon D_\S\to\h^*$
induced by $\pr_{\lie h^*}$.  Because $T^*H$ is multiplicity-free
(with respect to the $H$-action), $D_\S$ is multiplicity-free as well.
Because the $H$-moment map on $T^*H$ is surjective, the image of
$\Psi_\S$ is precisely the polyhedron $\P$.

Features of the labelled polyhedron are reflected in features of the
associated Delzant space in an interesting way.  Consider for instance
an arbitrary open face $\F$ of $\P$.  There are two sets of labels
naturally associated with $\F$, namely the set $\S_\F$ of all labels
in $\S$ corresponding to hyperplanes containing $\F$, and the set
$\S|_\F$ of all labels giving the equations and inequalities for
$\bar{\F}$.  In other words,
\begin{equation}\label{equation:littlelabel}
\begin{split}
\S_\F&=\bigl\{\,(v_i,r_i)\in\S:\text{the hyperplane
$\langle{\cdot},v_i\rangle=r_i$ contains $\F$}\,\bigr\},\\
\S|_\F&=\S\cup\bigl\{\,(-v_i,-r_i):(v_i,r_i)\in\S_\F\,\bigr\}.
\end{split}
\end{equation}
(We write an equation as a pair of inequalities to present $\bar{\F}$
as an intersection of half-spaces.)  Let $\h_\F\subset\h$ be the
subspace of $\lie h$ annihilating the tangent space to $\F$, that is
the linear span of all $v_i$ such that $(v_i,r_i)\in\S_\F$, and let
$H_\F=\exp\h_\F$ be the corresponding subtorus.  It is not hard to see
that the preimage of $\F$ under $\Psi_\S$ is a connected component of
the stratum of orbit type $H_\F$.  Furthermore, there is a natural
identification of the preimage of $\bar\F$ with the Delzant space
associated to the set of labels $\S|_\F$.  The result of these
observations is as follows.

\begin{proposition}\label{proposition:orbistrat}
The stabilizer of every point in $D_\S$ is connected and the
decomposition $\P=\bigcup_{\F\pre\P}\F$ of $\P$ into open faces $\F$
gives rise to the decomposition of $D_\S$ into $H$-orbit type
strata\upn:
\begin{equation}\label{equation:decompose1}
D_\S=\bigcup_{\F\pre\P}\Psi_\S\inv(\F).
\end{equation}
The stratum $\Psi_\S\inv(\F)$ is an $H/H_\F$-bundle over $\F$.  Its
closure is the Delzant space $D_{\S|_\F}$.  Its dimension is $2\dim
\F$ and the dimension of $D_\S$ is $2\dim\P$.  The $H$-fixed points in
$D_\S$ are the preimages of the vertices of $\P$.  If $\P$ is
$k$-dimensional\upn, the $H$-action is free on
$\Psi_\S\inv(\inter\P)$.  
\qed
\end{proposition}

The following statements are obvious from the construction.

\begin{proposition}\label{proposition:dilate}
For $t\neq0$ let $t\S$ be the set of labels
$\bigl\{(v_1,tr_1),\dots,(v_n,tr_n)\bigr\}$.  The polyhedron
associated to $t\S$ is the dilated polyhedron $t\P$\upn, and $D_{t\S}$
is symplectomorphic to $D_\S$ equipped with $t$ times its original
symplectic form.

If $r=0$\upn, then $\P$ is a cone with apex at the origin\upn, $D_\S$
is a symplectic cone\upn, and the moment map $\Psi_\S$ is homogeneous
of degree one.
\end{proposition}

When is $D_\S$ nonsingular?  To determine the orbifold stratification
of $D_\S$, regarded as a symplectic quotient of the manifold
$H\times\lie h^*\times\C^n$ by $T^n$, we must calculate the
stabilizers of the points in $\tilde{\psi}_r\inv(0)$ with respect to
the $T^n$-action, where $\tilde{\psi}_r$ is as in
\eqref{equation:tildepsi}.  Let $K$ be the kernel of $\bar{A}$, so
that we have an exact sequence
$$
\xymatrix{K\ar@{((->}[r] & T^n \ar[r]^{\bar{A}} & H.}
$$
Note that $\bar{A}$ is not necessarily surjective and that $K$ is not
necessarily connected.  The stabilizer of $(h,\eta,z)\in H\times\lie
h^*\times\C^n$ is $K\cap T^n_z$, where $T^n_z$ is the stabilizer of
$z\in\C^n$ under the $T^n$-action.  Consider an open face $\F$ of $\P$
and a point $(h,\eta,z)\in\tilde{\psi}_r\inv(0)$ such that $\eta\in
\F$, in other words the image under $\Psi_\S$ of the orbit
$T^n\cdot(h,\eta,z)\in D_\S$ is contained in $\F$.  Then it is easy to
check that $T^n_z$ is the torus $T^n_\F$ generated by the span of all
vectors $e_i\in\R^n$ such that the label $(v_i,r_i)$ is in $\S_\F$,
where $\S_\F$ is as in \eqref{equation:littlelabel}.  The homomorphism
$\bar{A}$ restricts to a surjective map $\bar{A}_\F\colon T^n_\F\to
H_\F$, and the stabilizer of $(h,\eta,z)$ is exactly the kernel
$K_\F$.  These groups form an exact sequence
$$
\xymatrix{K_\F\ar@{((->}[r] & T^n_\F\ar@{->>}[r]^{\bar{A}_\F} & H_\F,}
$$
where again $K_\F=K\cap T^n_\F$ is not necessarily connected.  The
dimension of $T^n_\F$ is the number of labels $l$ in the set $\S_\F$.
The dimension of $K_\F$ is the \emph{excess\/} $e_\S(\F)=l-\codim \F$
of the open face $\F$.  The \emph{excess function\/} of $(\S,\P)$ is
the step function
$$
e_\S=\sum_{\F\pre\P}e_\S(\F)\,1_\F,
$$
where $1_\F$ denotes the indicator function of $\F$.  Evidently,
$e_\S$ is upper semicontinuous on $\P$ and for all $p$ the set
$e_\S\inv(p)$ of constant excess $p$ is a union of open faces.  If
$\P$ is $k$-dimensional, then $e_\S$ vanishes on the interior of $\P$.
Now let $\{\P_\alpha:\alpha\in\eu A\}$ be the collection of all
connected components of all sets $e_\S\inv(p)$, where $p$ ranges over
$\N$.  This decomposition of $\P$ is called the \emph{excess
decomposition\/} of $\P$.  We have proved the first part of the
following result.

\begin{proposition}\label{proposition:inforbistrat}
\begin{enumerate}
\item
The infinitesimal $T^n$-orbit type stratification of the space $D_\S$
is induced by the excess decomposition
$\P=\bigcup_\alpha\P_\alpha$\upn:
\begin{equation}\label{equation:decompose2}
D_\S=\bigcup_\alpha D_{\S,\alpha},
\end{equation}
where $D_{\S,\alpha}=\Psi_\S\inv(\P_\alpha)$.  The connected subgroup
of $T^n$ corresponding to the stratum $D_{\S,\alpha}$ is the identity
component of $K_\F$\upn, where $\F$ is any of the faces in
$\P_\alpha$.  The structure group of a point $p$ in the orbifold
$D_{\S,\alpha}$ is the component group of $K_\F$\upn, where $\F\subset
\P_\alpha$ is the face containing $\Psi_\S(p)$.
\item\label{part:images}
The subsets $\P_\alpha$ and $\bar{\P}_\alpha$ of $\P$\upn, which are
the images under $\Psi_\S$ of respectively $D_{\S,\alpha}$ and
$\bar{D}_{\S,\alpha}$\upn, are convex.  In fact\upn, $\bar{\P}_\alpha$
is the closure of a single face of $\P$.
\item\label{part:orbi}
Assume that the excess function $e_\S$ is constant on $\P$.  Then for
every open face $\F$ of $\P$ the Delzant space $D_{\S|_\F}$ associated
to the labelled polyhedron $\bigl(\S|_\F,\bar\F\bigr)$ is an orbifold.
In particular\upn, $D_\S$ itself is an orbifold.
\item\label{part:mani}
Assume that $e_\S$ is constant on $\P$ and that for every open face
$\F$ the vectors in $\S_\F$ generate the weight lattice
$\Lambda\cap\h_\F$ of $H_\F$.  Then the Delzant spaces $D_{\S|_\F}$
are manifolds.  In particular\upn, $D_\S$ itself is a manifold.
\item 
For all $\alpha$ the following statements are equivalent\/\upn:
\begin{enumerate}
\item\label{part:polyclosed} $\P_\alpha$ is closed\/\upn;
\item\label{part:stratumclosed} $D_{\S,\alpha}$ is closed\/\upn;
\item\label{part:localmax} if $\F_1$ and $\F_2$ are open faces of $\P$
such that $\F_1\pre\P_\alpha$ and
$\F_2\cap\bar{\P}_\alpha\neq\emptyset$ but
$\F_2\not\pre\P_\alpha$\upn, then $e_\S(\F_2)<e_\S(\F_1)$\upn; 
\item\label{part:excessconst} the excess function $e_{\S|_\F}$ of the
set of labels $\S|_\F$ is constant\upn, where $\F$ is the open face of
$\P_\alpha$ such that $\bar{\F}=\bar{\P}_\alpha$.
\end{enumerate}
\end{enumerate}
\end{proposition}
  
\begin{proof}[Outline of proof]
Part \ref{part:images} follows easily from the upper semicontinuity of
$e_\S$ and the following observation: let $\F_1$ and $\F_2$ be any
pair of open faces of $\P$.  Let $\F_1\vee \F_2$ and $\F_1\wedge \F_2$
denote respectively the smallest open face $\F$ such that $\F\succeq
\F_1\cup \F_2$ and the interior of the intersection
$\bar{\F}_1\cap\bar{\F}_2$.  Then
$$e_\S(\F_1\vee \F_2)+e_\S(\F_1\wedge \F_2)=e_\S(\F_1)+e_\S(\F_2).$$

For part \ref{part:orbi}, observe that the excess function is constant
if and only if $0$ is a quasi-regular value of $\tilde{\psi}_r$, and
in this case the orbifold decomposition \eqref{equation:decompose2}
consists of only one piece.  To see that the spaces $D_{\S|_\F}$ are
also orbifolds, note that for any $\F_1$ and $\F_2$ such that
$\F_2\pre\F_1$ one has
\begin{equation}\label{equation:faceexcess}
e_{\S|_{\F_1}}(\F_2)=e_\S(\F_2)+l,
\end{equation}
where $l$ is the number of labels in $\S_{\F_1}$.  Consequently, if
$e_\S$ is constant on $\P$, then $e_{\S|_\F}$ is constant on
$\bar{\F}$ for all $\F$.

The hypotheses in \ref{part:mani} imply that the groups $K_\F$ are all
connected.  Hence the one orbifold occurring in the decomposition
\eqref{equation:decompose2} has trivial structure groups at all
points, so $D_\S$ is a manifold.  The proof that the $D_{\S|_\F}$ are
manifolds is similar.

The equivalence of \ref{part:polyclosed} and \ref{part:stratumclosed}
is trivial and the equivalence of \ref{part:polyclosed} and
\ref{part:localmax} follows from the upper semicontinuity of $e_\S$
and the convexity of $\P_\alpha$.  The equivalence of
\ref{part:polyclosed} and \ref{part:excessconst} is a consequence of
the equality \eqref{equation:faceexcess}.
\end{proof}

Note that the stratification \eqref{equation:decompose2} is coarser
than the $H$-orbit type stratification \eqref{equation:decompose1},
which is a stratification into manifolds.  Note further that $0$ is a
regular value of $\tilde{\psi}_r$ if and only if the excess function
is identically zero, that is to say for all open faces $\F$ the
labelling vectors occurring in $\S_\F$ are linearly independent (and
hence form a basis of the vector space $\h_\F$).  In this case the set
of labels and its associated polyhedron are called \emph{simple}.  If
$\S$ is simple\upn, then $\P$ is $k$-dimensional and the Delzant
spaces associated to all the closed faces of $\P$ are orbifolds.  The
torus $T^n$ acts freely on $\tilde{\psi}_r\inv(0)$ if and only if for
all open faces $\F$ the labelling vectors in $\S_\F$ are a basis of
the weight lattice $\Lambda\cap\h_\F$.  In this case $\S$ and $\P$ are
called \emph{simply laced}.

To what extent are the set $\S$ and the space $D_\S$ determined by the
polyhedron $\P$?  If $(v_i,r_i)$ is a label such that $\P$ is
contained in the open half-space $\langle\mu,v_i\rangle>r_i$, then it
can be omitted from the set of labels without affecting the structure
of $D_\S$.  Call $\S$ \emph{minimal\/} if it does not contain any such
labels.  It is clear that every set of labels can be pruned down to a
unique minimal one by dropping all redundant labels.  If $\S$ is
minimal and simple, then the labels are in one-to-one correspondence
with the faces of $\P$ and every label is determined by the
corresponding face up to multiplication by a positive integer.  A
complete understanding of this special situation is obtained in
\cite{le:ha}.  Multiplying a label by an integer does not affect the
homeomorphism type of $D_\S$, but does affect its orbifold structure.
If $\S$ is minimal and simply laced (the case originally considered by
Delzant), then $\S$ and $D_\S$ are in fact uniquely determined by
$\P$.  It is not hard to see that the Delzant spaces associated to
general sets of labels are also determined by the underlying polyhedra
up to equivariant homeomorphism, but their classification up to
isomorphism of stratified symplectic spaces appears to be unknown.

We conclude this discussion by analysing the structure of $D_\S$
transverse to the singularities.  Observe that for every $\F$ the set
$\S_\F$ is a set of labels for the subtorus $H_\F$.  Now take any
$\alpha\in\eu A$ and let $\F$ be the open face such that
$\P_\alpha\subset\bar{\F}$.  The polyhedron $\P_\F\subset\h_\F^*$
associated to $\S_\F$ has two important properties: it is a strictly
convex cone with apex the point $\pr_{\h_\F^*}\mu$, where $\mu$ is any
point in $\F$, and $\P_\alpha$ has a neighbourhood in $\P$ that is the
product of $\P_\alpha$ and a neighbourhood of the apex in $\P_\F$.  By
Proposition \ref{proposition:dilate} the Delzant space associated to
$\S_\F$ is a symplectic cone and after subtracting $\pr_{\h_\F^*}\mu$
its moment map is homogeneous of degree one.  The upshot is as
follows.

\begin{proposition}
A neighbourhood of $D_{\S,\alpha}$ in $D_\S$ is symplectomorphic to
the product of $\Psi_\S\inv(\P_\alpha)$ with the symplectic cone
$D_{\S_\F}$\upn, where $\F$ is the largest open face contained in
$\P_\alpha$.  
\qed
\end{proposition}

The symplectic link of $\P_\alpha$ is also a Delzant space.  We leave
it to the reader to determine its labelled polytope.

\section{Partial desingularizations}\label{section:desingularization}

If $M$ is a nonsingular projective variety over $\C$ on which the
complex reductive group $G^\C$ acts by projective linear
transformations, then the compact group $G$ acts on $M$ in a
Hamiltonian fashion, where the symplectic form is the imaginary part
of the Fubini-Study metric.  Under the hypothesis that the set of
stable points is nonempty (which amounts to the symplectic condition
that the zero locus of the moment map contains regular points), Kirwan
showed in \cite{ki:pa} how to construct explicitly a ``partial''
desingularization of the categorical quotient $X=M\qu G^\C$.  It is
birationally equivalent to $X$ and possesses at worst finite-quotient
singularities.  The construction consists in judiciously blowing up
subvarieties of $M$ until all semistable points become stable, and
subsequently dividing out by the action of $G^\C$.  Kirwan pointed out
that this method works also in the symplectic category if one replaces
complex blowups by symplectic blowups in the sense of Gromov.
However, symplectic blowups depend on a number of choices and the
partial desingularization obtained by this method is not unique up to
symplectomorphism.  A further contrast with the algebraic case is that
there is no natural way to define a blowdown map.

In Sections \ref{subsection:canonical}--\ref{subsection:deform} we
review Kirwan's desingularization process in the symplectic setting
and show that the result is determined uniquely up to symplectic
deformations.  Deformation equivalence is much weaker than
symplectomorphism, but suffices for our purpose of calculating
characteristic numbers.  A problem is that the blowup centres involved
are not compact and are indeed not even uniquely defined.  We handle
this difficulty by applying a version of the constant-rank embedding
theorem ``with parameters'' proved in Appendix
\ref{section:normalform}.  Further we show that Kirwan's method works
equally well for $G$-actions on orbifolds and remove the hypothesis
that the zero locus of the moment map should contain a regular point.
In Section \ref{subsection:shift} we compare her desingularization to
the shift desingularizations discussed in Section
\ref{subsection:singular}.  Finally in Section
\ref{subsection:delzant2} we apply both methods to Delzant spaces.

\subsection{The canonical partial desingularization}
\label{subsection:canonical}

\subsubsection{Symplectic blowing up}\label{subsubsection:blow}

Let $S$ be a locally closed $G$-invariant symplectic suborbifold of
$M$ and let $\K$ be a $G$-invariant closed subset of $M$ such that
$\K\cap S$ is compact.  It is possible to define a symplectic blowup
of $M$ along $S$, at least in a neighbourhood of $\K$. (If $S$ is
compact, we can take $\K=M$, but we are mainly interested in cases
where $S$ is not compact.)  This is most easily accomplished in terms
of symplectic cutting and involves certain auxiliary data
$(j,\theta,\iota,\eps)$.  Here $j$ is a $G$-invariant compatible
complex structure on the normal bundle $N$ of $S$ and $\theta$ a
$G$-invariant principal connection on the orbibundle $P$ of unitary
frames in $N$.  To define $\iota$ and $\eps$ we note that $N$ is
isomorphic to the associated vector orbibundle $P\times^{\U(k)}\C^k$,
where $2k$ is the codimension of $S$ in $M$ and $\C^k$ is Euclidean
space equipped with its standard symplectic form
$\frac{i}{2}\sum_{l=1}^k\d z_l\wedge\d\bar{z}_l$.  According to
Theorem \ref{theorem:minimalcoupling} the minimal coupling form on $N$
defined by means of $\theta$ is nondegenerate in a neighbourhood $N'$
of the zero section.  By the symplectic embedding theorem there exists
a $G$-equivariant symplectic embedding $\iota\colon N'\to M$ such that
the diagram
$$
\xymatrix{S\ar[d]\ar[dr]^{\subset}\\N'\ar[r]_\iota & M}
$$
commutes.  Here the vertical arrow denotes the zero section of $N$.
By \ref{part:fibre} of Theorem \ref{theorem:minimalcoupling} the
$S^1$-action on $N$ defined by scalar multiplication on the fibres is
Hamiltonian with moment function $\psi(v)=-\frac1{2}\Vert v\Vert^2$,
where $\Vert{\cdot}\Vert$ denotes the fibre metric on $N$.  For
$\delta>0$ we have respectively the open and closed disc orbibundles
and the sphere orbibundle,
\begin{align*}
N(\delta) &=\{\,v\in N:\psi(v)>-\delta\,\}, \\
\bar N(\delta) &=\{\,v\in N:\psi(v)\ge-\delta\,\}, \\
\bb{S}N(\delta) &=\{\,v\in N:\psi(v)=-\delta\,\}.
\end{align*}
Because $\K\cap S$ is compact, there exist $\delta>0$ and a
$G$-invariant open neighbourhood $U'$ of $\K\cap S$ such that $\bar
N(\delta)|_{U'\cap S}$ is contained in $N'$ and $\iota$ embeds $\bar
N(\delta)|_{U'\cap S}$ properly into $U'$.  In addition, because $\K$
is closed there exists a $G$-invariant open subset $U''$ of $M$ such
that the union $U=U'\cup U''$ contains $\K$ and $U''$ does not
intersect the image of $\bar N(\delta)|_{U'\cap S}$ under $\iota$.
Then $U\cap S=U'\cap S$ is a closed suborbifold of $U$,
$\iota\bigl(\bar N(\delta)|_{U\cap S}\bigr)$ is a closed subset of
$U$, and the complement in $U$ of $\iota\bigl(N(\delta)|_{U\cap
S}\bigr)$ is a suborbifold with boundary
$\iota\bigl(\bb{S}N(\delta)|_{U\cap S}\bigr)$.  Now let
$0<\eps<\delta$ and put $\psi_\eps(v)=\psi(v)+\eps$.  Clearly, $0$ is
a regular value of $\psi_\eps$.  Let us identify $\bar
N(\delta)|_{S\cap U'}$ with its image under $\iota$.  Then the
symplectic cut of $U$ with respect to $\psi_\eps$ is well-defined.

\begin{definition}\label{definition:blowup}
The \emph{blowup\/} $\Bl(U,S,\omega,j,\theta,\iota,\eps)$, or
$\Bl(U,S,\eps)$, of $U$ with \emph{centre\/} $S$ is the Hamiltonian
$G$-orbifold $U_{\le0}$ obtained by cutting $U$ with respect to the
function $\psi_\eps$.  The \emph{exceptional divisor\/} is the
symplectic orbifold $U_0$.
\end{definition}

Thus the orbifold $\Bl(U,S,\omega,j,\theta,\iota,\eps)$ is obtained by
excising from $U$ the open disc orbibundle about $U\cap S$ and
collapsing the orbits of the $S^1$-action on the bounding sphere
orbibundle.  The exceptional divisor embeds symplectically into the
blowup.  It is symplectomorphic to the total space of (the restriction
to $U\cap S$ of) the orbibundle $\bb{P}N=P\times^{\U(k)}\C P^{k-1}$,
the projectivization of $N$, whose general fibre is the
$k-1$-dimensional complex projective space with $\eps$ times the
Fubini-Study form.  The blowup is diffeomorphic to the orbifold
$(U-S)\coprod\ca I$ obtained by gluing together $U-S$ and the
incidence relation $\ca I\subset(U\cap N')\times\bb{P}N|_{U\cap S}$
along the obvious map.  This space is called the \emph{null-blowup\/}
and denoted by $\Bl(U,S,\omega,j,\iota,0)$ or $\Bl(U,S,0)$.  It is
plainly independent of the data $(\theta,\eps)$.  Note that unlike the
symplectic blowups the null-blowup has a canonical blowdown map
$\pi\colon\Bl(U,S,0)\to U$.

The symplectomorphism type (in fact, the symplectic volume) of a
symplectic blowup depends on the parameter $\eps$.  On the other hand,
it is not hard to see that given a fixed $\eps>0$ the choice of $j$
and $\theta$ does not affect the blowup.  How the symplectic structure
of the blowup depends on the embedding $\iota$ is a delicate problem
and is only partly understood even in the case where $S$ is a point.
(See \cite{mc:re,mc:fr}.)  Observe however that the blowup does not
actually depend on the embedding, but only on the cutting data, that
is the image of the sphere orbibundle and the distance function under
the embedding.  The following elementary result, which is proved in
\cite{mc:in} for $S$ a point, says roughly that the blowup does not
depend on $(j,\theta,\iota)$ as long as $\eps$ is sufficiently small.
In Section \ref{subsection:deform} we consider what happens if
$\omega$, $S$ and $\eps$ are allowed to vary.

\begin{proposition}\label{proposition:smalleps}
For all triples $(j_0,\theta_0,\iota_0)$ and $(j_1,\theta_1,\iota_1)$
there exist $\delta>0$ and $G$-invariant open neighbourhoods $U_0$ and
$U_1$ of $\K$ such that for all $\eps<\delta$ the blowups
$\Bl(U_0,S,\omega,j_0,\theta_0,\iota_0,\eps)$ and
$\Bl(U_1,S,\omega,j_1,\theta_1,\iota_1,\eps)$ are isomorphic
Hamiltonian $G$-orbifolds.
\end{proposition}

\begin{proof}
First we reduce the problem to the case where $j_0=j_1$ and
$\theta_0=\theta_1$.  For $i=0$, $1$, let $\omega_i$ and
$\Vert{\cdot}\Vert_i$ denote the minimal coupling form and fibre
metric on $N$ associated to $(j_i,\theta_i)$.  By the symplectic
embedding theorem we can find $\delta'<\delta$, $U\supset\K$ and an
embedding $\upsilon\colon \bigl(N(\delta')|_{U\cap
S},\omega_0\bigr)\to \bigl(N(\delta)|_{U\cap S},\omega_1\bigr)$ of
Hamiltonian $G\times S^1$-orbifolds fixing the zero section.  Since
$\upsilon$ intertwines the $S^1$-moment maps, it maps the sphere
orbibundles associated to the fibre metric $\Vert{\cdot}\Vert_0$ to
those for $\Vert{\cdot}\Vert_1$.  This means that for $\eps<\delta'$
the cutting data used to define the blowup
$\Bl(U',S,\omega,j_0,\theta_0,\iota'_1,\eps)$ are identical to those
for the blowup $\Bl(U',S,\omega,j_1,\theta_1,\iota_1,\eps)$, where
$\iota'_1$ is the composite embedding $\iota_1\circ\upsilon$.  It
follows that these two blowups are isomorphic.

It remains to find $U'_0$ and $U'_1$ such that
$\Bl(U'_0,S,\omega,j_0,\theta_0,\iota'_1,\eps)$ is isomorphic to
$\Bl(U'_1,S,\omega,j_0,\theta_0,\iota_0,\eps)$ for small $\eps$.  If
$S$ is a point, the map $\iota_0\inv\circ\iota'_1$ is symplectically
isotopic to the identity in a small enough neighbourhood of $S$.  This
is in general false, but by Lemma \ref{lemma:circle} we can find
$\delta''<\delta'$ and a $G$-equivariant symplectic isotopy $H\colon
N(\delta'')\times[0,1]\to N(\delta')$ fixing $S$, starting at
$\iota_0\inv\circ\iota'_1$ and ending at an $S^1$-equivariant map.  As
before, this implies that $H_1$ preserves the sphere orbibundles and
the distance function in $N(\delta'')$.  In other words, the map
$\iota_0 H_1(\iota'_1)\inv$ maps the cutting data used to define the
blowup with respect to the embedding $\iota'_1$ to those for the
embedding $\iota_0$.  To finish the proof it suffices show that this
map can be extended to a global equivariant symplectomorphism
$U_0''\to U_1''$ of suitable open $U_0''$ and $U_1''$ containing $\K$.
This is accomplished by applying Lemma \ref{lemma:hamiltonian}, which
says that the isotopy $\iota_0 H_t(\iota'_1)\inv$ is generated by a
time-dependent Hamiltonian $f_t$, which in the present situation we
can arrange to be $G$-invariant.  Let $\chi$ be a $G$-invariant
cut-off function on $\iota'_1N(\delta'')$ that is compactly supported
in the fibre directions and is equal to $1$ on $\iota'_1N(\delta''')$
for some $\delta'''<\delta''$.  Because of the compactness of $\K\cap
S$, on a sufficiently small open $U_0''$ containing $\K$ the
Hamiltonian vector field of $\chi f_t$ integrates to a globally
defined Hamiltonian isotopy $F\colon U_0''\to U'$ which restricts to
$\iota_0 H_t(\iota'_1)\inv$ on $\iota'_1N(\delta''')$.  We now take
$U_1''$ to be the image of $U_0''$ under $F_1$ and conclude that for
$\eps<\delta'''$ the blowup
$\Bl(U''_0,S,\omega,j_0,\theta_0,\iota'_1,\eps)$ is isomorphic to
$\Bl(U''_1,S,\omega,j_0,\theta_0,\iota_0,\eps)$.
\end{proof}

Consider a $G$-equivariant line orbibundle $L$ on $U$.  Then the
restriction of $L$ to $N(\delta)$ is isomorphic to $\pr_S^*(L|_S)$,
which is trivial in the fibre directions.  Fix an isomorphism
$L|_{N(\delta)}\cong\pr_S^*(L|_S)$ and lift the action of $S^1$ on
$N(\delta)$ to $\pr_S^*(L|_S)$ by letting it act trivially on the
fibres.

\begin{definition}\label{definition:lineblow}
For $0<\eps<\delta$ the \emph{blowup}
$\Bl(L,S,\omega,j,\theta,\iota,\eps)$, or $\Bl(L,S,\eps)$, of
$L$ along $S$ is the $G$-equivariant orbibundle $L_{\le0}$ on
$\Bl(U,S,\eps)$ obtained by cutting $L$ with respect to the function
$\psi_\eps$.  The \emph{null-blowup\/} $\Bl(L,S,\omega,j,\iota,0)=
\Bl(L,S,0)$ is the pullback of $L$ to $\Bl(U,S,0)$.
\end{definition}

It is easy to see that the blowup of $L$ does not depend on the
identification of $L|_{N(\delta)}$ with $\pr_S^*(L|_S)$.  We emphasize
that the blowup of a prequantum line bundle is \emph{not\/}
prequantizing, because the curvature form of $\Bl(L,S,\eps)$ is
degenerate along the exceptional divisor. Indeed, under a
diffeomorphism $\Bl(U,S,0)\to\Bl(U,S,\eps)$ the pullback of
$\Bl(L,S,\eps)$ is isomorphic to $\Bl(L,S,0)$, whose curvature is
equal to the presymplectic form $\pi^*\omega$.

Analogous constructions can be carried out in the presence of an
arbitrary $S^1$-action that is defined on an open neighbourhood of $S$
and leaves $S$ fixed.  The weights
$\alpha=(\alpha_1,\alpha_2,\dots,\alpha_k)$ of the action are assumed
to be positive.  In this situation $N$ can be written as
$P\times^K\C^k$, where $K$ is the centralizer in $\U(k)$ of the circle
subgroup defined by the weights $\alpha$, and $P$ is now the
orbibundle of $S^1$-equivariant unitary frames.  The symplectic cut of
$U$ with respect to the $S^1$-action is referred to as a
\emph{weighted blowup\/} and is denoted by
$\Bl(U,S,\omega,j,\theta,\iota,\eps,\alpha)$ or
$\Bl(U,S,\eps,\alpha)$.  The exceptional divisor of a weighted blowup
of $\C^k$ is a weighted projective space.

\subsubsection{Resolving singularities}\label{subsubsection:resolve}

Suppose that $0$ is not a quasi-regular value of the moment map
$\Phi\colon M\to\g^*$.  Then the stratification of the quotient
$X=M\qu G$ consists of more than one piece.  By performing a
succession of equivariant blowups of an open neighbourhood $U$ of
$Z=\Phi\inv(0)$ we shall define a Hamiltonian $G$-orbifold
$\bigl(\tilde U,\tilde\omega,\tilde\Phi\bigr)$ such that $0$ is a
quasi-regular value of $\tilde{\Phi}$.  The blowup centres are ``the''
minimal $G$-invariant symplectic suborbifolds containing the strata of
maximal depth in $Z$.  (If $G$ is abelian, there is a canonical choice
for the blowup centres and we can furthermore arrange for them to be
closed, so that the blowups are globally defined.  The reason is that
for all $\alpha$ the closure of $M_{(\lie g_\alpha)}$ is a symplectic
suborbifold of $M$ and that $Z_\alpha$ is coisotropic in $\bar
M_{(\lie g_\alpha)}$.  The component of $\bar M_{(\lie g_\alpha)}$
containing $Z_\alpha$ is therefore a closed minimal symplectic
suborbifold around $Z_\alpha$.)  Following Kirwan \cite{ki:pa} we call
the symplectic quotient $\tilde X$ of $\tilde U$ the \emph{canonical
partial desingularization\/} of $X$, although the process is by no
means as canonical as in the algebraic case.  We investigate in
Section \ref{subsection:deform} to what extent the result is
well-defined.

Recall that the \emph{depth} of a stratum $Z_\alpha$ is the largest
integer $i$ for which there exists a strictly ascending chain of
strata $Z_{\alpha_0}\prec Z_{\alpha_1}\prec\dots\prec Z_{\alpha_i}$
with $\alpha_0=\alpha$.  The depth of $Z$ is the maximum of the depths
of all its strata.

Let $Z_\alpha$ be a stratum of maximal depth.  Then $Z_\alpha$ is
closed and hence, because $\Phi$ is proper, $Z_\alpha$ and $X_\alpha$
are compact.  Choose an embedding of an open set $U_\alpha\subset
M_\alpha$ as in Theorem \ref{theorem:model}.  Let $S_\alpha$ be the
minimal symplectic suborbifold containing $Z_\alpha$ defined in
\eqref{equation:minimalsymplectic}.  Then $U_\alpha\cap S_\alpha$ is
locally closed in $M$ and by \eqref{equation:levelbundle} its
intersection with $Z$ is the closed stratum $Z_\alpha$.  We conclude
that there exists a $G$-invariant open neighbourhood $U$ of the
compact set $\K=Z$ such that for $0<\eps<\delta$ the blowup
$\Bl(U,S_\alpha,\eps)$ is well-defined.  Recall that the symplectic
normal bundle of $S_\alpha$ in $U_\alpha$ is the model space
$M_\alpha$.  A neighbourhood of the exceptional divisor in
$\Bl(U,S_\alpha,\eps)$ is therefore modelled by
$$
\Bl(M_\alpha,S_\alpha,\eps)=
P_\alpha\times^{K_\alpha}\bigl(G\times^{G_\alpha}(\g_\alpha^0\times
\Bl(W_\alpha/\Upsilon_\alpha,0,\eps))\bigr).
$$
The point of blowing up is that it reduces the depth of the
stratification.

\begin{lemma}
\begin{enumerate}
\item\label{part:ray}
Let $H$ be a connected Lie group acting on a vector orbispace
$W/\Upsilon$ in a unitary fashion.  Let $w\ne0$ be a point in
$\Phi_{W/\Upsilon}\inv(0)$ and let $[w]\in\bb{P}(W/\Upsilon)$ be the
ray through $w$.  Then $w$ and $[w]$ have the same infinitesimal
stabilizer\upn, $\h_w=\h_{[w]}$.
\item\label{part:depth}
$\depth\Bl(M_\alpha,S_\alpha,\eps)\qu G=\depth M_\alpha\qu G-1$.
\end{enumerate}
\end{lemma}

\begin{proof}
Since $[w]$ is fixed under $\h_{[w]}$, there is an infinitesimal
character $\sigma\in(\h_{[w]})^*$ such that
\begin{equation}\label{equation:projective}
\exp (\eta)\cdot w =e^{2\pi i\langle\sigma,\eta\rangle}w
\end{equation} 
for all $\eta\in \h_{[w]}$.  Hence
$\langle\Phi_{W/\Upsilon}(w),\eta\rangle= \frac1{2}\omega(\eta\cdot
w,w)= -\pi\langle\sigma,\eta\rangle\lVert w\rVert^2$ for all $\eta\in
\h_{[w]}$.  Now $w\ne0$ and $\Phi_{W/\Upsilon}(w)=0$, so $\sigma=0$.
By \eqref{equation:projective} we conclude that $\h_{[w]}\subset\h_w$.
The reverse inclusion is obvious.  This proves \ref{part:ray}.

For the proof of \ref{part:depth} note that the depths of $M_\alpha\qu
G$ and $\Bl(M_\alpha,S_\alpha,\eps)\qu G$ are equal to those of
$(W_\alpha/\Upsilon_\alpha)\qu G_\alpha$ and
$\Bl(W_\alpha/\Upsilon_\alpha,0,\eps)\qu G_\alpha$, respectively.
Recall that the set of $G_\alpha$-fixed points in
$W_\alpha/\Upsilon_\alpha$ consists of the origin only.  In other
words, the stratum of maximal depth in
$\Phi_{W_\alpha/\Upsilon_\alpha}\inv(0)$ consists of the origin only.
Assertion \ref{part:ray} is therefore tantamount to saying that the
zero level set of $\Bl(W_\alpha/\Upsilon_\alpha,0,\eps)$ contains
exactly the same infinitesimal orbit types as the zero level set of
$W_\alpha/\Upsilon_\alpha$, except for the maximal element
$\g_\alpha$.
\end{proof}

This lemma tells us that by successively blowing up neighbourhoods of
the zero level set along minimal symplectic suborbifolds containing
the strata of maximal depth in the zero level set, we arrive
eventually at a Hamiltonian $G$-orbifold $\bigl(\tilde
U,\tilde\omega,\tilde\Phi\bigr)$ whose zero level set has depth $0$,
that is to say $0$ is a quasi-regular value of $\tilde\Phi$.  The
canonical partial desingularization of $X=M\qu G$ is by definition the
orbifold
$$ \tilde X=\tilde U\qu G=\tilde Z/G.$$

For abelian $G$, at each stage in the process the blowup centre is
uniquely defined and closed, so if $M$ is compact, the successive
blowups are globally defined.

If $L$ is an almost equivariantly locally trivial line orbibundle on
$M$, then by blowing it up at each step along the way, we obtain an
almost equivariantly locally trivial line orbibundle $\tilde L$ on
$\tilde U$.  We define the canonical desingularization of $L_0=L\qu G$
to be the line orbibundle $\tilde L_0=\tilde L\qu G$.

Note finally that by choosing appropriate symplectic cross-sections we
can apply the process described above to partially resolve the
singularities of $M_\mu$ for every value $\mu$ of the moment map.

\subsection{Deformation equivalence}\label{subsection:deform}

A \emph{deformation\/} of the equivariant symplectic form
$(\omega,\Phi)$ on $M$ is a smooth path of equivariant symplectic
forms $(\omega_t,\Phi_t)$ defined for $0\le t\le1$ such that
$(\omega_0,\Phi_0)=(\omega,\Phi)$ and $\Phi_t\inv(0)=\Phi\inv(0)$ for
all $t$.  The endpoints $(\omega_0,\Phi_0)$ and $(\omega_1,\Phi_1)$ of
the path are called \emph{deformation equivalent}.  A
\emph{deformation equivalence\/} between $(M,\omega,\Phi)$ and a
second Hamiltonian $G$-orbifold $(M',\omega',\Phi')$ is a
diffeomorphism $F$ from $M$ to $M'$ such that $(\omega,\Phi)$ is
deformation equivalent to $(F^*\omega',F^*\Phi')$.

For trivial $G$-actions (where $\Phi=0$) this reduces to the usual
notion of a deformation or pseudo-isotopy; see \cite[Ch.~6]{mc:in}.
For the purposes of this paper its interest is first of all that it is
preserved under symplectic reduction and that a deformation
equivalence class determines a homotopy class of almost complex
structures.

\begin{lemma}\label{lemma:deform}
\begin{enumerate}
\item\label{part:deformregular}
Let $(\omega_t,\Phi_t)$ be a deformation of $(\omega,\Phi)$.  If $0$
is a quasi-regular value of $\Phi$\upn, then for all $t$ it is a
quasi-regular value of $\Phi_t$\upn, and the $\omega_t$ induce a
smooth path of symplectic forms on $\Phi_t\inv(0)/G=\Phi\inv(0)/G$.
\item\label{part:deformcomplex}
The space of all $G$-invariant almost complex structures on $M$ that
are compatible with some equivariant symplectic structure in the
deformation equivalence class of $(\omega,\Phi)$ is nonempty and
path-connected.
\end{enumerate}
\end{lemma}

\begin{proof}
Part \ref{part:deformregular} follows directly from the observation
that the zero level set is by definition fixed under a deformation and
that its stratification depends not on the symplectic form but on the
$G$-action alone.  Part \ref{part:deformcomplex} follows from the fact
that for every path $\omega_t$ of invariant symplectic forms on $M$
there exists a path of invariant complex structures $J_t$ on $TM$ such
that $J_t$ is compatible with $\omega_t$.
\end{proof}

The second important property is that the operation of symplectic
cutting (and hence the operation of blowing up) is well-behaved with
respect to deformations. 

\begin{definition}\label{definition:family}
A \emph{family of cutting data\/} consists of sextuples
$(\omega_t,\Phi_t,U_t,Z_t,\psi_t,F_t)$ defined for $0\le t\le1$. Here
$(\omega_t,\Phi_t)$ is a path of equivariant symplectic forms (we do
not require $\Phi_t\inv(0)=\Phi\inv(0)$), $(U_t,Z_t,\psi_t)$ are
$G$-invariant cutting data with respect to the symplectic form
$\omega_t$ (see Section \ref{subsubsection:cut}), and $F\colon
U_0\times[0,1]\to M$ is an isotopy of the open subset $U_0$.  These
data are subject to the following conditions: $F_0$ is the identity
map of $U_0$, $F_t$ is $G$-equivariant and maps $U_0$ onto $U_t$, and
the path of $G$-equivariant symplectic forms
$(F_t^*\omega_t,F_t^*\Phi_t)$ on $U_0$ is a deformation of
$(\omega,\Phi)$.  Furthermore, $F_t$ is to be equivariant with respect
to the given $S^1$-actions on $U_0$ and $U_t$, and the path of
$S^1$-equivariant symplectic forms $(F_t^*\omega_t,F_t^*\psi_t)$ on
$U_0$ is to be a deformation of $(\omega,\psi_0)$.  Finally $\psi_t$
is required to depend smoothly on $t$ and $0$ is required to be a
regular value of $\psi_t$ for all $t$.
\end{definition}

These conditions entail
$$
F_t\bigl(\Phi_0\inv(0)\bigr)=\Phi_t\inv(0), \qquad F_t(Z_0)=Z_t,
\qquad (F_t)_*\Xi_0=\Xi_t,
$$
where $\Xi_t$ is the Hamiltonian vector field of $\psi_t$ with respect
to $\omega_t$.  We denote by $(M,Z_t)_{\ge0}$ the symplectic cut of
$M$ with respect to $(U_t,Z_t,\psi_t)$ and by $(M,Z_t)_{0}$ the
symplectic quotient.  For brevity let us denote the extension of the
equivariant symplectic form $(\omega_t,\Phi_t)$ to $(M,Z_t)_{\ge0}$
also by $(\omega_t,\Phi_t)$.  Mark that $F_t$ does not necessarily
pull back the function $\psi_t$ to $\psi_0$, but that it does map
$\psi_0\inv\bigl([0,\infty)\bigr)$ to
$\psi_t\inv\bigl([0,\infty)\bigr)$.

We wish to show that a family of cutting data gives rise to
deformation equivalent symplectic cuts.  To do this, we need first to
extend $F$ to a global $G$-equivariant diffeotopy $\check F$ of $M$
leaving the fibre $\Phi\inv(0)$ invariant.  Let $\eta_t$ be the
infinitesimal generator of $F$; it is a time-dependent vector field
supported on the track $\eu U=\bigcup_tU_t\times\{t\}$ of the isotopy.
Let $\eu U'$ and $\eu U''$ be $G\times S^1$-invariant tubular
neighbourhoods of the hypersurface $\bigcup_tZ_t\times\{t\}$ in $\eu
U$ such that $\eu U''\subset\eu U'$, and choose a $G\times
S^1$-invariant bump function $\chi\colon\eu U\to[0,1]$ that is
supported on $\eu U'$ and identically equal to $1$ on $\eu U''$.
Extend the vector field $\chi\eta_t$ by $0$ to a global smooth vector
field on $M\times[0,1]$.  Its flow, $\check F$, is defined for $0\le
t\le1$ and supported on $\eu U$.  It is clearly $G$-equivariant and
$S^1$-equivariant (where the $S^1$-action is defined).  In addition,
$\check F$ is equal on $\eu U''$ to the previously defined flow $F$,
and on $\eu U$ its trajectories are subsets of those of $F$, so
$\check F$ likewise preserves the set $\Phi\inv(0)$.

For simplicity we ignore henceforth the distinction between $F$ and
$\check F$.  Note that $F_t$ maps $(M,Z_0)_{>0}$ to $(M,Z_t)_{>0}$.

\begin{proposition}\label{proposition:deformcut}
For every $t$ the restriction of $F_t$ to $(M,Z_0)_{>0}$ extends
uniquely to a diffeomorphism
$\bar{F}_t\colon(M,Z_0)_{\ge0}\to(M,Z_t)_{\ge0}$.  The restriction of
$\bar{F}_t$ to $(M,Z_0)_{0}$ is equal to the map
$(M,Z_0)_{0}\to(M,Z_t)_{0}$ induced by $F_t$.  The path
$(\bar{F}_t^*\omega_t,\bar{F}_t^*\Phi_t)$ of equivariant symplectic
forms on $(M,Z_0)_{\ge0}$ is a deformation of $(\omega_0,\Phi_0)$.
\end{proposition}

\begin{proof}
Define $\tilde{F}_t$ from $(U_0-Z_0)\times\C$ to $(U_t-Z_t)\times\C$
by
$$
\tilde{F}_t(u,z)= \biggl(F_t(u),
\biggl\lvert\frac{\psi_t\bigl(F_t(u)\bigr)}{\psi_0(u)}
\biggr\rvert^{1/2} z\biggr).
$$
This map is clearly $S^1$-equivariant and maps $\tilde{\psi}_0\inv(0)$
to $\tilde\psi_t\inv(0)$.  (Recall
$\tilde \psi(m,z)=\psi(m)-\frac1{2}\vert z\vert^2$.)  We assert that it
extends to a diffeomorphism from $U_0\times\C$ to $U_t\times\C$.  To
show this, we observe that $0$ is a regular value of the composite
function $\psi_t\circ F_t\colon U_0\to\R$, and therefore, in suitable
local coordinates $(x_1,\dots,x_{2n})$ about a point in $Z$ in which
$Z$ is given by $x_1=0$, it can be written as
$\psi_t\bigl(F_t(x_1,\dots,x_{2n})\bigr)= a(t)x_1$ with $a$ smooth and
nowhere vanishing.  It follows that
$$
\tilde{F}_t(u,z)= \biggl(F_t(u), \biggl\lvert\frac{a(t)}{a(0)}
\biggr\rvert^{1/2} z\biggr)
$$
extends smoothly to the locus $\{x_1=0\}$.  Notice also that if
$\beta_t$ is the embedding given by \eqref{equation:beta} (with $\psi$
replaced by $\psi_t$), then $\tilde{F}_t\bigl(\beta_0(u)\bigr)=
\beta_t\bigl(F_t(u)\bigr)$, so the maps $F_t$ and
$\tilde{F}_t|_{\tilde{\psi}_0\inv(0)}$ can be glued together to give a
smooth map $\bar{F}_t$ from $(M,Z_0)_{\ge0}$ to $(M,Z_t)_{\ge0}$.  The
proof that the inverse of $F_t$ extends smoothly is similar.  If
$z=0$, then $\tilde{F}_t(u,0)=\bigl(F_t(u),0\bigr)$ for all $u\in
U_0$, which implies that the restriction of $\bar{F}_t$ to
$(M,Z_0)_{0}$ is equal to the map induced by $F_t$.

To show that the path $(\bar{F}_t^*\omega_t,\bar{F}_t^*\Phi_t)$ is a
deformation of $(\omega_0,\Phi_0)$, we need merely show that
$\bar{F}_t$ maps $\Phi_0\inv(0)$ onto $\Phi_t\inv(0)$.  This follows
from the fact that its restriction to $(M,Z_0)_{>0}$ is equal to
$F_t$, which sends $\Phi_0\inv(0)\cap(M,Z_0)_{>0}$ onto
$\Phi_t\inv(0)\cap(M,Z_t)_{>0}$.
\end{proof}

\begin{proof}[Proof of Theorem \ref{theorem:deform}]
The first statement of the theorem is to be interpreted as follows:
suppose we apply Kirwan's desingularization process (at the level
$\mu$) and that by making two sets of choices of the parameters
involved in the process we obtain two Hamiltonian $G$-orbifolds
$(\tilde U,\tilde\omega,\tilde\Phi)$ and $(\tilde
U',\tilde\omega',\tilde\Phi')$.  Then there exist invariant open
neighbourhoods $\check U$ of $\tilde\Phi\inv(G\mu)$ in $\tilde U$ and
$\check U'$ of $(\tilde\Phi')\inv(G\mu)$ in $\tilde U'$ and a
deformation equivalence $\check U\to\check U'$.

After choosing an appropriate cross-section we may assume that
$\mu=0$.  Observe that the first statement of the theorem together
with \ref{part:deformregular} of Lemma \ref{lemma:deform} implies the
uniqueness of the symplectic structure on $\tilde{M}_0$ up to
deformations.  Combined with \ref{part:deformcomplex} of Lemma
\ref{lemma:deform} this shows that the Chern classes of the tangent
bundle of $\tilde{M}_0$ are well-defined.  Hence the Riemann-Roch
numbers of $\tilde{M}_0$ are well-defined.

It remains to prove the first statement of the theorem.  The proof is
by induction on the depth of $Z$.  If $\depth Z=0$ there is nothing to
prove.  The inductive step is taken by establishing the following two
facts.  Let $Z_\alpha$ be a stratum of maximal depth in $Z$.  Then

\begin{enumerate}
\item\label{fact:1}
for all minimal $G$-invariant symplectic submanifolds $S_0$ and $S_1$
containing $Z_\alpha$ and for all blowup data
$(j_0,\theta_0,\iota_0,\eps_0)$ and $(j_1,\theta_1,\iota_1,\eps_1)$
relative to the blowup centres $S_0$, resp.\ $S_1$, there exist
invariant open neighbourhoods $U_0$ and $U_1$ of $\Phi\inv(0)$ such
that the blowup $\Bl(U_0,S_0,\omega,j_0,\theta_0,\iota_0,\eps_0)$ is
deformation equivalent to
$\Bl(U_1,S_1,\omega,j_1,\theta_1,\iota_1,\eps_1)$;
\item\label{fact:2}
for every deformation $(\omega_t,\Phi_t)$ of $(\omega,\Phi)$ there
exist $\delta>0$ and a family of data
$(U_t,S_t,E_t,j_t,\theta_t,\iota_t,\psi_t,F_t)$ such that the
following conditions are satisfied: $U_t$ is a $G$-invariant open
subset containing $\Phi_t\inv(0)=\Phi\inv(0)$, $S_t$ is a closed
$G$-invariant minimal $\omega_t$-symplectic submanifold of $U_t$
containing $Z_\alpha$; $E_t$ is the $\omega_t$-symplectic normal
bundle of $S_t$ in $U_t$; $j_t$ and $\theta_t$ are resp.\ an invariant
compatible almost complex structure and an invariant connection on
$E_t$; $\iota_t$ is a proper $G$-invariant $\omega_t$-symplectic
embedding of $\bar{E}_t(\delta)$ into $U_t$; $\psi_t$ is the function
$-\frac1{2}\Vert{\cdot}\Vert_t^2$, where $\Vert{\cdot}\Vert_t$ denotes
the fibre metric on $E_t$; and finally $F$ is a $G$-equivariant
isotopy $E_0(\delta)\times[0,1]\to M$ starting at the identity.  We
require that $F_t$ maps $E_0(\delta)$ onto $E_t(\delta)$ and is
equivariant for the $S^1$-actions on $E_0(\delta)$ and $E_t(\delta)$;
that $F_t$ maps $E_0(\delta)\cap\Phi\inv(0)$ into $\Phi_t\inv(0)$, and
$\psi_0\inv(-\eps)$ into $\psi_t\inv(-\eps)$ for $\eps<\delta$.
\end{enumerate}

Notice that the first fact suffices to prove the theorem if $\depth
Z=1$.  The second fact says that a deformation of the equivariant
symplectic form on $M$ gives rise to a family of cutting data in the
sense of Definition \ref{definition:family}.  Proposition
\ref{proposition:deformcut} then implies that the resulting symplectic
cuts (i.~e.\ blowups) are deformation equivalent.  Combined with fact
\ref{fact:1} this says that the relation of being deformation
equivalent is preserved when going through a single stage in Kirwan's
desingularization process.  This completes the inductive step when the
depth is greater than $1$.

We now proceed to prove facts \ref{fact:1} and \ref{fact:2}.
Proposition \ref{proposition:smalleps} and Theorem
\ref{theorem:centre} show that for all quadruples
$(S_0,j_0,\theta_0,\iota_0)$ and $(S_1,j_1,\theta_1,\iota_1)$ there
exist $\delta>0$ and $G$-invariant open $U_0$ and $U_1$ containing
$\Phi\inv(0)$ such that for all $\eps<\delta$ the blowup
$\Bl(U_0,S_0,\omega,j_0,\theta_0,\iota_0,\eps)$ is \emph{isomorphic\/}
to $\Bl(U_1,S_1,\omega,j_1,\theta_1,\iota_1,\eps)$ as a Hamiltonian
$G$-orbifold.  Furthermore, according to \eqref{equation:levelbundle}
the projection $M_\alpha\to S_\alpha$ in the model space $M_\alpha$
preserves the zero fibre of $\Phi$.  Therefore, by Lemma
\ref{lemma:blowdeform} the blowups
$\Bl(U_0,S_0,\omega,j_0,\theta_0,\iota_0,\eps)$ and
$\Bl(U_0,S_0,\omega,j_0,\theta_0,\iota_0,\eps_0)$ are deformation
equivalent as Hamiltonian $G$-orbifolds, and so are the blowups
$\Bl(U_1,S_1,\omega,j_1,\theta_1,\iota_1,\eps)$ and
$\Bl(U_1,S_1,\omega,j_1,\theta_1,\iota_1,\eps_1)$.  This proves fact
\ref{fact:1}.

The proof of fact \ref{fact:2} invokes the relative constant-rank
embedding theorem, Theorem \ref{theorem:constantrank}.  We use the
notation introduced before Theorem \ref{theorem:constantrank} and take
$\eu M=M\times[0,1]$ and $\eu Z=Z_\alpha\times[0,1]$.  The relative
symplectic form $\omega_{\eu M}$ is defined by $\omega_{\eu
M}|_{M\times\{t\}}= \omega_t$, so that the form $\tau=\omega_{\eu
M}|_{\eu Z}$ has constant rank $\frac1{2}\dim Z_\alpha/G$ on $\eu Z$.
Let $\eu N$ be the relative symplectic normal bundle of $\eu Z$ in
$\eu M$.  For the construction of the data
$(S_t,E_t,j_t,\theta_t,\iota_t,\psi_t,F_t)$ we may assume that we are
working in the standard model $\eu Y=\eu S\oplus\eu N$.  Let $Y_t=\eu
Y|_{\{t\}}$; then $Y_t$ is symplectomorphic to $(M,\omega_t)$ near
$Z_\alpha$.  Recall that $\eu S$ is defined as $(\ker\tau)^*$ and that
$\ker\tau$ is equal to the distribution tangent to the $G$-orbits on
$Z_\alpha\times[0,1]$, which does not depend on the value of the base
point in $[0,1]$.  Therefore $\eu S$ is equal to the product
$S_0\times[0,1]$, where $S_0=\eu S|_{Z_\alpha\times\{0\}}$.  Now
define $S_t$ to be $\eu S|_{\eu Z\times\{t\}}=S_0\times\{t\}$.  Then
$S_t$ is independent of $t$ as a manifold, although its symplectic
form may depend on $t$.  The symplectic normal bundle $\eu E$ of $\eu
S$ is equal to $\eu Y$, considered as an orbibundle over $\eu S$.  The
unit interval being contractible, there exists a $G$-equivariant
isomorphism of symplectic vector orbibundles $F\colon
E_0\times[0,1]\to\eu E$, where $E_0=\eu E|_{S_0}$.  Let $F_t\colon
E_0\to E_t$ be the map $E_0\to\eu E|_{S_t}$ defined by $F$; then $F_t$
is an isomorphism of symplectic vector orbibundles covering the
diffeomorphism $S_0\to S_t=S_0\times\{t\}$.  By construction $\eu
E|_{S_t}$ is simply the $\omega_t$-symplectic normal bundle of $S_t$
in $Y_t$.  Choose an invariant compatible almost complex structure
$j_0$ and an invariant connection $\theta_0$ on $E_0$ and put
$j_t=(F_t\inv)^*j_0$ and $\theta_t=(F_t\inv)^*\theta_0$.  We let
$\iota_t$ be the standard embedding from $E_t$ into $Y_t$ and
$\psi_t=-\frac1{2}\Vert{\cdot}\Vert_t^2$, where $\Vert{\cdot}\Vert_t$
is the fibre metric on $E_t$ with respect to the complex structure
$j_t$.  Again by construction, $F_t$ is a $G\times S^1$-equivariant
isomorphism of Hermitian vector orbibundles and therefore maps disc
bundles into disc bundles and sphere bundles into sphere bundles.
Recall that by \eqref{equation:levelbundle} the projection
$M_\alpha\to S_\alpha$ preserves $\Phi_t\inv(0)$ for all $t$, so
$\Phi_t\inv(0)=\Phi_{S_t}\inv(0)\cap\Phi_{\theta_t}\inv(0)$ for all
$t$ by Example \ref{example:blowdeform}.  Furthermore
$\Phi_{S_t}\inv(0)=Z_\alpha$ and $F_t$ maps $\theta_0$ to $\theta_t$
for all $t$, so we conclude that $F_t$ maps
$E_0(\delta)\cap\Phi_0\inv(0)$ into $\Phi_t\inv(0)$.

It remains to define the open neighbourhoods $U_t$ of $\Phi\inv(0)$.
We do this by starting with an invariant open $U_0$ containing
$\Phi\inv(0)$ such that the embedding $\bar E_0(\delta)\to U_0$ is
proper.  We extend the infinitesimal generator of $F_t$ to a globally
defined time-dependent vector field on $U_0$ by means of a suitable
cut-off function.  Because the level set $\Phi\inv(0)$ is compact,
after shrinking $U_0$ and $\delta$ if necessary the resulting vector
field is integrable for $0\le t\le1$.  As a result we obtain an
equivariant isotopy $\check F\colon U_0\times[0,1]\to M$ preserving
$\Phi\inv(0)$; and we put $U_t=F_t(U_0)$.
\end{proof}

\subsection{Shift desingularizations}\label{subsection:shift}

The process delineated in Section \ref{subsubsection:resolve} is
usually not the most economical method for resolving the singularities
of the quotient.  A simpler desingularization is often obtained by
shifting the value of the moment map to a nearby quasi-regular value.
Let $\Delta_i$ be one of the open chambers of the moment polyhedron as
in \eqref{equation:chamber} and assume $\mu$ is in its closure.  Since
$\Phi$ has maximal rank on $\Delta_i$, all $\nu\in\Delta_i$ are
quasi-regular values of $\Phi$ and all $M_\nu$ have the same
dimension.  If this dimension is the same as that of $M_\mu$, the
$M_\nu$ are called \emph{shift desingularizations\/} of $X$.  We
discuss briefly the relationship between $M_\nu$ and the canonical
desingularization $\tilde M_\mu$ (even if they do not have the same
dimension).  We say that two symplectic orbifolds $Q_0$ and $Q_1$ are
\emph{related\/} by a weighted symplectic blowup (resp.\ blowdown) if
$Q_1$ is deformation equivalent to a weighted blowup (resp.\ blowdown)
of $Q_0$ at a closed symplectic suborbifold.

\begin{theorem}\label{theorem:shiftcanonical}
Let $\Delta_i$ be an open chamber of $\Delta$ such that $\mu$ is in the
closure of $\Delta_i$.  There exists a symplectic fibre orbibundle $E$
over $\tilde M_\mu$ with the property that for all $\nu$ in $\Delta_i$
the symplectic orbifold $M_\nu$ is related to $E$ by a sequence of
weighted symplectic blowups and blowdowns.  The general fibre of $E$
is a generic reduced space of the space
$F(G_\alpha,W_\alpha/\Upsilon_\alpha)$\upn, where $\alpha$ is the
infinitesimal orbit type of the open stratum of $Z$.
\end{theorem}
 
\begin{proof}[Sketch of proof]
After applying the symplectic cross-section theorem we may assume that
$\mu=0$.  By the implicit function theorem the symplectic quotients
$M_\nu$ for $\nu\in\Delta_i$ are all deformation equivalent, so it
suffices to prove the theorem for $\nu$ close to $0$.  If $\nu$ is
sufficiently close to $0$ we can arrange, by choosing the parameters
$\eps$ in the desingularization process small enough, that $\nu$ is
contained in the image of $\tilde U$ under $\tilde\Phi$ and lies
outside the images of the exceptional divisors arising in the process.
Then $\nu$ is a generic value of $\tilde\Phi$ and the quotients
$M_\nu$ and $\tilde U_\nu$ are symplectomorphic.  We now apply
Proposition \ref{proposition:varquot} to the space $\tilde U$ and
conclude that for small generic values $\nu'$ of $\tilde\Phi$ the
quotient $E=\tilde U_{\nu'}$ is a fibre orbibundle over $\tilde M_0$
with general fibre $F(G_\alpha,W_\alpha/\Upsilon_\alpha)_{\nu'}$.
According to the symplectic cross-section theorem, the preimage
$Y=\tilde\Phi\inv\bigl(\inter(\tilde\Phi(\tilde U)\cap\tplus)\bigr)$
is a Hamiltonian $T$-orbifold and $Y_\nu\cong M_\nu$ and
$Y_{\nu'}\cong E$.  A result of Guillemin and Sternberg \cite{gu:bi}
now says that $Y_\nu$ is related to $Y_{\nu'}$ by a sequence of
weighted blowups and blowdowns.
\end{proof}

\subsection{Delzant spaces II}\label{subsection:delzant2}

Let $\S=\bigl\{(v_1,r_1),(v_2,r_2),\dots,(v_n,r_n)\bigr\}$ be a set of
labels for the torus $T$.  Assume that the polyhedron $\P$ associated
to $\S$ is nonempty.  We can apply both desingularization processes to
the Delzant space $D_\S$ defined in Section \ref{subsection:delzant1}.
Let us first do the canonical desingularization.  By Proposition
\ref{proposition:inforbistrat} a closed stratum $D_{\S,\alpha}$ in
$D_\S$ corresponds to a piece $\P_\alpha$ in the excess decomposition
of $\P$ that is the closure of a single open face $\F$.  Let $\S_\F$
be its associated set of labels as defined in
\eqref{equation:littlelabel}.  Blowing up at $D_{\S,\alpha}$ has the
effect of adding one label to the set $\S$:
$$
\tilde\S_{\eps_1}=
\bigl\{(v_1,r_1),\dots,(v_n,r_n),(v_{n+1},r_{n+1}+\eps_1)\bigr\},
$$
where
$$ 
v_{n+1}=\sum_{(v_i,r_i)\in\S_\F}v_i,\qquad
r_{n+1}=\sum_{(v_i,r_i)\in\S_\F}r_i
$$
and $\eps_1>0$ is sufficiently small.  The excess decomposition of the
polyhedron $\tilde\P_{\eps_1}$ associated to $\tilde\S_{\eps_1}$ has
one piece less than that of $\P$.  Iterating this process we obtain
eventually a labelled polyhedron
$\bigl(\tilde\S_\eps,\tilde\P_\eps\bigr)$ with constant excess
function, where $\eps$ denotes a vector with small positive entries.
The canonical desingularization of $D_\S$ is then the Delzant space
$\tilde D_\S=D_{\tilde\S_\eps}$.

Shift desingularization has the effect of replacing $\S$ by a set of
labels
\begin{equation}\label{equation:shiftlabel}
\S_\eta=\bigl\{(v_1,r_1+\eta_1),(v_2,r_2+\eta_2),\dots,
(v_n,r_n+\eta_n)\bigr\},
\end{equation}
where $\eta$ is a small vector chosen in such a way that the
associated polyhedron $\P_\eta$ is nonempty and has constant excess
function.  The shape of $\P_\eta$ depends on the choice of $\eta$, but
the directions of the faces of codimension one do not.

\begin{figure}
\setlength{\unitlength}{0.1mm}
$$
\begin{picture}(695,560)
\path(2,333)(182,333)(302,393)(152,543)(2,333) \path(182,333)(152,543)
\dashline{4.000}(2,333)(122,393)(302,393)
\dashline{4.000}(122,393)(152,543)
\path(392,333)(572,333)(692,393)(584,500)
\path(584,500)(550,482)(500,482)(533,500)(584,500)
\path(392,333)(500,482) \path(572,333)(550,482)
\dashline{4.000}(392,333)(512,393)(692,393)
\dashline{4.000}(513,393)(533,500)
\path(2,33)(182,33)(302,93) \path(182,33)(167,213)(302,93)
\path(2,33)(137,213)(167,213) \dashline{4.000}(2,33)(122,93)(302,93)
\dashline{4.000}(122,93)(137,213)
\path(392,33)(572,33)(692,93) \path(392,33)(533,213)(572,33)
\path(692,93)(552,219)(533,213) \dashline{4.000}(513,92)(552,219)
\dashline{4.000}(392,33)(513,92)(692,93)
\put(122,290){\makebox(0,0)[lb]{\smash{$\P$}}}
\put(512,290){\makebox(0,0)[lb]{\smash{$\tilde\P$}}}
\put(122,-10){\makebox(0,0)[lb]{\smash{$\P'$}}}
\put(512,-10){\makebox(0,0)[lb]{\smash{$\P''$}}}
\end{picture}
$$
\caption{Egyptian pyramid and desingularizations}
\label{diagram:pyramid}
\end{figure}
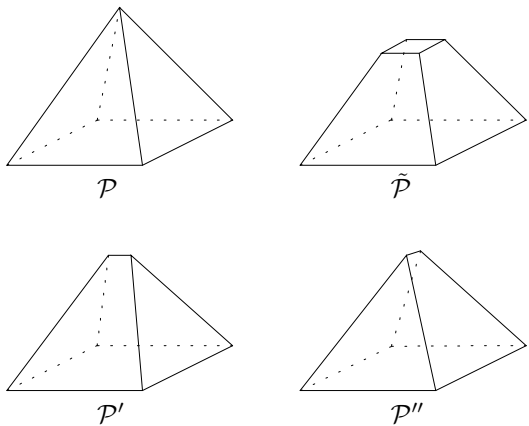

See Diagram \ref{diagram:pyramid} for the canonical desingularization
$\tilde\P$ and two different shift desingularizations $\P'$ and $\P''$
of the Egyptian pyramid.  The manifold $X$ corresponding to the
truncated pyramid $\tilde\P$ is a $\C P^1$-bundle over $\C P^1\times\C
P^1$.  The section at infinity $Y$ is a product $Y'\times Y''$, where
$Y'$ and $Y''$ are copies of $\C P^1$, and the manifolds corresponding
to $\P'$ and $\P''$ are obtained by blowing down $X$ at $Y'$, resp.\
$Y''$.  The space corresponding to $\P$ is obtained by blowing down
$X$ at the divisor $Y$.

\section{The abelian case}\label{section:abelian}

In this section we prove the results stated in Section
\ref{section:results} in the case of abelian group actions.  In
Section \ref{subsection:almostcomplex} the setting is that of an
almost complex orbifold $\ca M$ with an arbitrary equivariant line
orbibundle $\ca L$.  Here the Atiyah-Segal-Singer fixed-point formula
for the equivariant index yields imprecise, but useful, qualitative
information about the multiplicity diagram of $\ca L$.  The discussion
is inspired by that in \cite{du:sy}, which in turn bears great
similarity to the arguments of \cite{ko:ap,at:sp,bo:on}.  We deduce
the first part of Theorem \ref{theorem:rigid} (for general groups) as
an immediate corollary.  We then specialize to the setting of
Hamiltonian actions in Section \ref{subsection:hamiltonian}.  The
result of Section \ref{subsection:almostcomplex} furnishes precise
information on the Riemann-Roch numbers of symplectic quotients at
certain vertices of the moment polytope.  A cut-and-paste argument
involving symplectic cutting and a gluing formula for the equivariant
index then produces formul\ae\ for the Riemann-Roch numbers of
symplectic quotients at other points.  In preparation for Section
\ref{section:nonabelian} we then generalize our results to
``asymptotic'' moment bundles.  Finally in Section
\ref{subsection:delzant3} we apply them to Delzant spaces, where they
lead to a counting formula for lattice points in rational polytopes.

\subsection{Almost complex $T$-orbifolds}
\label{subsection:almostcomplex}

Let $\ca M$ be a \emph{compact\/} almost complex $T$-orbifold and $\ca
L\to\ca M$ an arbitrary $T$-equivariant line orbibundle.  As in
Section \ref{subsection:preliminaries}, let $\Lambda=\ker\exp$ denote
the integral lattice of $T$, $\RR(\ca M,\ca L)\in\Rep T$ the
equivariant index of $\dirac_{\ca L}$, and $N_{\ca
L}\colon\Lambda^*\to\Z$ the multiplicity function of $\ca L$.  For all
components $F$ of the fixed-point set $\ca M^T$, let
$\sigma_F\in\Lambda^*\otimes\Q$ be the orbiweight of the $T$-action on
$\ca L|_F$.  The \emph{weight polytope\/} $\Delta_{\ca L}$ of $\ca L$
is the rational convex polytope
$$
\Delta_{\ca L}=\hull\bigl\{\sigma_F:F\subset \ca M^T\bigr\}.
$$
Let $N_F$ be the normal bundle of $F$ in $\ca M$ and let
$\alpha_{jF}\in\Lambda^*\otimes\Q$ be the orbiweights for the action
of $T$ on $N_F$, where $j=1$, $2,\dots$, $\codim_\C F$.  Denote by
$\ca C_F$ be the rational cone in $\t^*$ spanned by the $-\alpha_{jF}$,
and by $\check{C}_F$ its dual cone
$$ 
\check\ca C_F=\{\,\xi\in\t:\langle\alpha_{jF},\xi\rangle\le0 \mbox{
for all $j$}\,\}.
$$
The data $\Delta_{\ca L}$ and $\ca C_F$ record certain information on
the multiplicity diagram of $\ca L$.  A vector $\xi\in\lie t$ is
called \emph{generic\/} if the one-parameter subgroup generated by
$\xi$ has the same fixed points as $T$.  This is equivalent to
$\langle\alpha_{jF},\xi\rangle\ne0$ for all $F$ and all $j$.  For
instance, $\xi$ is generic if its one-parameter subgroup is dense in
$T$.

\begin{theorem}\label{theorem:vertex}
\begin{enumerate}
\item\label{part:support}
The support of the multiplicity function is contained in the weight
polytope.
\item\label{part:vertexweight}
Let $\mu$ be a vertex of the weight polytope.  Choose a generic $\xi$
in $\lie t$ such that $\langle\mu-\sigma_F,\xi\rangle<0$ for all $F$
with $\sigma_F\ne\mu$. Then
\begin{equation}\label{equation:fixed}
N_{\ca L}(\mu)=\sum_{\substack{F\\ \sigma_F=\mu\\ \xi\in
\check{C}_F}}\RR(F,\ca L|_F).
\end{equation}
\item\label{part:tcrossh}
If $H$ is a compact connected Lie group acting on $\ca M$ and $\ca
L$\upn, and the action of $H$ commutes with that of $T$ and preserves
the almost complex structure on $\ca M$\upn, then
\eqref{equation:fixed} holds as an equality of virtual characters of
$H$.
\end{enumerate}
\end{theorem}

The support of $N_{\ca L}$ is by definition contained in the weight
lattice, so we have in fact $\supp N\subset\Lambda^*\cap\Delta_{\ca
  L}$.  This implies that if the vertex $\mu$ under
\ref{part:vertexweight} is not integral, then $N_{\ca L}(\mu)=0$, so
the right-hand side of \eqref{equation:fixed} vanishes.

\begin{proof} 
We assume first that $\ca M$ has no orbifold singularities.  According
to the Atiyah-Segal-Singer equivariant index formula the character is
equal to a sum
\begin{equation}\label{equation:index}
\RR(\ca M,\ca L)=\sum_F\chi_F
\end{equation}
over all fixed-point components $F$ of $M$.  The functions $\chi_F$
are for generic $\xi\in\lie t$ given by
\begin{equation}\label{equation:fix}
\chi_F(\exp\xi)=e^{2\pi
i\langle\sigma_F,\xi\rangle}\int_F\frac{\Ch(\ca
L|_F)\Td(F)}{\D_T(N_F,\xi)}.
\end{equation}
Here $\Ch(\ca L|_F)$ is the Chern character of $\ca L|_F$, $\Td(F)$ is
the Todd class of $F$, and $\D_T(N_F,\xi)$ is the $T$-equivariant
characteristic class
\begin{equation}\label{equation:split}
\D_T(N_F,\xi)=\prod_j\bigl(1-\exp(-2\pi
i\langle\alpha_{jF},\xi\rangle-c_{jF})\bigr)
\end{equation}
of $N_F$.  The $c_{jF}$ are the virtual Chern roots of $N_F$, that is
to say, $c_{jF}=c_1(N_{jF})$ if $N_F$ decomposes into a sum
$\bigoplus_jN_{jF}$ of $T$-equivariant line bundles where $T$ acts
with weight $\alpha_{jF}$ on the $j$th summand.  Because
\eqref{equation:split} is symmetric in the $c_{jF}$, $\D_T(N_F,\xi)$
is well-defined by the splitting principle.

The character $\RR(\ca M,\ca L)\colon T\to\C$ extends to a holomorphic
function on the complexified torus $T^\C$ and the $\chi_F$ extend to
rational functions on $T^\C$.  Substitute $-it\xi$ for $\xi$ in
\eqref{equation:fix} and \eqref{equation:split}, where $t$ is real,
and put $x=e^{2\pi t}$.  Notice that \eqref{equation:split} has a
finite limit as $x\to0$.  This limit is $1$ if
$\langle\alpha_{jF},\xi\rangle<0$ for all $j$, and $0$ otherwise.
Applying the Riemann-Roch Theorem to $F$ and $\ca L|_F$ we find that
as $x\to 0$
\begin{equation}\label{equation:limit}
\chi_F\bigl(\exp(-it\xi)\bigr)=
\begin{cases}
x^{\langle\sigma_F,\xi\rangle}\RR(F,\ca
L|_F)+o\bigl(x^{\langle\sigma_F,\xi\rangle}\bigr) &\text{ if
$\langle\alpha_{jF},\xi\rangle<0$ for all $j$}, \\
o\bigl(x^{\langle\sigma_F,\xi\rangle}\bigr)& \text{ otherwise}.
\end{cases}
\end{equation}
Compare this estimate to the following expression obtained from
\eqref{equation:rr}:
\begin{equation}\label{equation:character}
\RR(\ca M,\ca L)\bigl(\exp(-it\xi)\bigr)= \sum_{\nu\in\Lambda^*}N_{\ca
L}(\nu)x^{\langle\nu,\xi\rangle}.
\end{equation}
Assume that $\mu$ is not contained in the weight polytope $\Delta_{\ca
L}$.  Then we can select a generic $\xi$ such that
$\langle\mu-\sigma_F,\xi\rangle<0$ for all $F$, i.~e.\ 
$\langle\mu,\xi\rangle<\min_F\langle\sigma_F,\xi\rangle$.  According
to \eqref{equation:index} and \eqref{equation:limit},
$\RR(M,L)\bigl(\exp(-it\xi)\bigr)=
O\bigl(x^{\min_F\langle\sigma_F,\xi\rangle}\bigr)$ as $x\to 0$, so all
terms in \eqref{equation:character} with exponent strictly less than
$\min_F\langle\sigma_F,\xi\rangle$ vanish, and therefore $N_{\ca
L}(\mu)=0$.  This proves \ref{part:support}.

If $\mu$ is a vertex of $\Delta_{\ca L}$, we can choose a generic
$\xi$ such that $\langle\mu-\sigma_F,\xi\rangle<0$ for all $F$ with
$\sigma_F\ne\mu$, and in addition
$\langle\nu,\xi\rangle\ne\langle\mu,\xi\rangle$ for those (finitely
many) $\nu\in\Lambda^*$ for which $N_{\ca L}(\nu)\ne0$.  If
$\langle\nu,\xi\rangle<\langle\mu,\xi\rangle$, then
$\langle\nu-\sigma_F,\xi\rangle<0$, so $\nu$ is not in $\Delta_{\ca
L}$ and therefore $N_{\ca L}(\nu)=0$ by \ref{part:support}.  This
implies that $\langle\nu,\xi\rangle>\langle\mu,\xi\rangle$ whenever
$\nu\ne\mu$ and $N_{\ca L}(\nu)\ne0$, and therefore, by
\eqref{equation:character}, $\RR(M,L)\bigl(\exp(-it\xi)\bigr)= N_{\ca
L}(\mu)x^{\langle\mu,\xi\rangle}+o\bigl(x^{\langle\mu,\xi\rangle}\bigr)$
as $x\to0$.  By \eqref{equation:index} and \eqref{equation:limit}, the
coefficient of the term of order $\langle\mu,\xi\rangle$ in the
character is equal to the sum of the $\RR(F,\ca L|_F)$ over all $F$
with the property that $\mu=\sigma_F$ and
$\langle\alpha_{jF},\xi\rangle<0$ for all $j=1$, $2,\dots$, $\codim_\C
F$.  This proves \ref{part:vertexweight}.
 
In the presence of orbifold singularities and the action of an
additional group $H$ we invoke the results of
\cite{be:eq,ve:eq,du:he}, according to which $\RR(\ca M,\ca L)$,
considered as an element of $\Rep(T\times H)$, is given by a sum
\eqref{equation:index}, where
$$ 
\chi_F(\exp\xi,\exp\eta)= e^{2\pi
i\langle\sigma_F,\xi\rangle} \int_{\tilde
F} \frac1{d_{\tilde F}}\frac{\Ch_H^{\tilde F}\bigl(\tilde{\ca
L}|_{\tilde F},\eta\bigr)\Td_H\bigl(\tilde
F,\eta\bigr)}{\D_H^{\tilde F}(N_{\tilde
F},\eta)\D_{T\times H}^{\tilde F}\bigl(\tilde N_{\tilde
F},\xi,\eta\bigr)}
$$
for $\eta\in\lie h$ sufficiently small.  Here $\Ch_H$ etc.\ are the
$H$-equivariant counterparts of the characteristic classes considered
above, $\tilde F$ is the ``unwrapping'' of the orbifold $F$, and
$d_{\tilde F}$ is its multiplicity, which is a locally constant
function on $\tilde F$.  See \cite{me:sym} for a detailed
discussion.  Our previous argument goes through with trivial
modifications and the upshot is that \ref{part:support} and
\ref{part:vertexweight} hold for orbifolds and that
\ref{part:vertexweight} holds $H$-equivariantly.
\end{proof}

\begin{proof}[Proof of part \ref{part:almostcomplex} of Theorem
\ref{theorem:rigid}] Note first that the result for arbitrary $G$
follows from the case where $G$ is abelian.  The weight polytope of a
rigid orbibundle is by definition $\{0\}$, so the abelian case is an
immediate consequence of Theorem \ref{theorem:vertex}.
\end{proof}

\subsection{Hamiltonian $T$-orbifolds}\label{subsection:hamiltonian}

In this section $(M,\omega,\Phi)$ denotes a compact connected
Hamiltonian $T$-orbifold and $L$ a $T$-equivariant line orbibundle on
$M$.  In this situation we have two polytopes, namely the moment
polytope $\Delta=\Phi(M)$ and the weight polytope $\Delta_L$, and
Theorem \ref{theorem:vertex} yields information on the index
$\RR(M_\mu,L_\mu)$ for certain vertices $\mu$ of $\Delta$.  For those
bundles $L$ for which there is a simple relationship between the
polytopes $\Delta$ and $\Delta_L$, namely rigid, moment, and dual
moment bundles, $\RR(M_\mu,L_\mu)$ can then be calculated when $\mu$
is not a vertex of $\Delta$ by means of multiple symplectic cutting
and the gluing formula, which we review in Section
\ref{subsubsection:glue}.  In Section
\ref{subsubsection:multiplicities} we put these ingredients together
to obtain proofs of the remaining theorems of Section
\ref{section:results} for abelian groups.

\subsubsection{Multiple symplectic cutting}\label{subsubsection:glue}

Let $\S$ be a set of labels in $\lie t^*$ and let $\P$ be its
associated polyhedron.  Assume that the excess function of $(\S,\P)$
is constant.  Then the Delzant space $D_\S$ is an orbifold by
Proposition \ref{proposition:inforbistrat} and the Hamiltonian
$T$-orbifold $D_{-\S}$ is symplectomorphic to $D_\S$ with the opposite
symplectic form by Proposition \ref{proposition:dilate}.  The
\emph{symplectic cut of $M$ with respect to\/} $\P$ is the Hamiltonian
$T$-space
\begin{equation}\label{equation:multiplecut}
M_\P=(M\times D_{-\S})\qu T
\end{equation}
obtained by reduction at $0$ with respect to the diagonal $T$-action
and the moment map $\Phi\times-\Psi_\S$.  There are several
alternative ways to think of $M_\P$.  Firstly, by Definition
\ref{definition:delzant} and reduction in stages,
$$
M_\P= (M\times T^*T\times\C^n)\qu(T^n\times T)= (M\times\C^n)\qu T^n.
$$
In particular, the Delzant space $D_\S$ itself is equal to the
symplectic cut of $T^*T$ with respect to $\P$.  Secondly, $M_\P$ is
the space obtained by performing successive symplectic cuts on $M$
with respect to each of the labels in $\S$.  Lastly, as a topological
space it is equal to the inverse image $\Phi\inv(\P)$ in which, for
each open face $\F$ of $\P$, one divides out the preimage
$\Phi\inv(\F)$ by the $T_\F$-action.  Thus we have a decomposition
\begin{equation}\label{equation:cut}
M_\P=\bigcup_{\F\pre\P}\Phi\inv(\F)/T_\F.
\end{equation}
We designate the moment map for the $T$-action on $M_\P$ by $\Phi_\P$.
Its image is equal to $\Delta\cap\bar\P$.

\begin{definition}\label{definition:tadmissible}
The pair $(\S,\P)$ is \emph{admissible\/} or \emph{$T$-admissible\/}
with respect to $M$ if $\S$ has constant excess and the reduction in
\eqref{equation:multiplecut} is regular.
\end{definition}

The reduction being regular is equivalent to $\lie t_m\cap\lie
t_x=\{0\}$ for all $(m,x)$ in $M\times D_{-\S}$ such that
$\Phi(m)=\Psi_\S(x)$.  If $\F$ is the open face of $\P$ that contains
$\Psi_\S(x)$, then $\lie t_x=\lie t_\F$ by Proposition
\ref{proposition:orbistrat}, so admissibility amounts to the condition
\begin{equation}\label{equation:admissible}
\lie t_m\cap\lie t_\F=\{0\}\qquad\text{for all $m\in M$ and all
$\F\pre\P$ such that $\Phi(m)\in\F$.}
\end{equation}
This has two consequences: firstly, because $\lie t_\F$ depends only
on the face $\F$, admissibility depends only the polyhedron $\P$ and
not on the set of labels defining it; and secondly, if $\P$ is
admissible, then for every open face $\F$, $0$ is a regular value of
the moment map $\Phi\times-\Psi_{\S_\F}$ on the Hamiltonian
$T$-orbifold $M\times D_{\S_\F}$.  In other words, every closed face
of $\P$ is admissible as well.

\begin{example}\label{example:cut}
If $T=S^1$, then $\lie t=i\R$ and $\Lambda=2\pi i\Z$.  Let
$\S=\{(v,0)\}$ where $v=2\pi i\in\Lambda$, then $\P=[0,\infty)$ and
$D_{-\S}$ is $\C$ with the standard symplectic structure and circle
action.  So $M_\P$ is equal to the symplectic cut $M_{\ge0}$ and
$M_{-\P}$ is the opposite cut $M_{\le0}$.  If $\P'=\{0\}$ then
$M_{\P'}=M_0$.  Admissibility of $\P$ is equivalent to $0$ being a
regular value of $\Phi$.  If we multiply the labelling vector $v$ by
$k$, then $\P$ does not change, but the orbifold atlas on $M_\P$
changes and $M_{\P'}$ becomes the symplectic quotient $M_0$ counted
$k$ times.
\end{example}

It is useful to rephrase \eqref{equation:admissible} in combinatorial
terms.  Consider the infinitesimal orbit type
stratification
\begin{equation}\label{equation:stratification}
M=\bigcup_{\beta\in\eu B}M_\beta
\end{equation}
of $M$ and denote the subalgebra corresponding to $\beta$ by $\lie
t_\beta$.  Choose a basepoint $m_\beta\in M_\beta$ and put
$\mu_\beta=\Phi(m_\beta)$ for all $\beta$.  The closure of $M_\beta$
is a $T$-invariant symplectic suborbifold of $M$; the affine span of
the convex polytope $\Phi(\bar M_\beta)$ is $\mu_\beta+\lie t_\beta^0$
and its relative interior is $\Phi(M_\beta)$.  (See e.~g.\
\cite{gu:co1}.)  The sets $\Phi(M_\beta)$ are called the \emph{virtual
open faces\/} of $\Delta$.  Since $\lie t_\F$ is the annihilator of
the tangent space of $\F$, the following statement is clear.

\begin{lemma}\label{lemma:transverse}
A polyhedron $\P$ defined by a set of labels of constant excess is
admissible if and only if its open faces are transverse to the virtual
open faces of $\Delta$.  Consequently\upn, admissibility is a generic
condition.
\end{lemma}

If $\P$ is admissible, the pullback $\pr_M^*L$ of $L$ to $M\times
D_{-\S}$ is almost equivariantly locally trivial at level $0$, and so
the \emph{cut bundle}
$$L_\P=\pr_M^*L\qu T$$
is a well-defined line orbibundle on $M_\P$.  Likewise, $L$ induces
well-defined orbibundles $L_{\P'}$ on each of the cuts $M_{\P'}$.
From \eqref{equation:isoline} we obtain for every open face $\F$ of
$\P$ a canonical isomorphism
\begin{equation}\label{equation:cutbundles}
L_\P\big|_{\Phi_\P\inv(\F)/T_\F}
\cong\bigl(L|_{\Phi\inv(\F)}\bigr)\big/T_\F.
\end{equation}

\begin{definition}\label{definition:tsubdivision}
An \emph{admissible\/} or \emph{$T$-admissible polyhedral
subdivision\/} of $\lie t^*$ is a collection $\lie P$ satisfying the
following conditions: every element of $\lie P$ is a $T$-admissible
polyhedron in $\t^*$, their union is $\lie t^*$, for every element of
$\lie P$ all its closed faces are in $\lie P$, and the intersection of
any two elements of $\lie P$ is a closed face of each.
\end{definition}

\begin{theorem}[gluing formula, \cite{me:sym}]\label{theorem:glue} 
Let $\lie P$ be an admissible polyhedral subdivision of $\lie t^*$.
Then
\begin{equation}\label{equation:glue}
\RR(M,L)=\sum_{\P\in\lie P}(-1)^{\codim\P}\RR(M_\P,L_\P)
\end{equation}
as virtual characters of $T$.  If $T=S^1$ and the cutting data are
only locally defined\upn, then we have a numerical identity
\begin{equation}\label{equation:localglue}
\RR(M,L)=\RR(M_{\ge0},L_{\ge0})+\RR(M_{\le0},L_{\le0})-\RR(M_0,L_0).
\end{equation}
In the presence of a compact connected Lie group $H$ that acts on $M$
and $L$ in such a way that $H$ commutes with $T$ and the $H$-action on
$M$ is symplectic\upn, the equalities \eqref{equation:glue} and
\eqref{equation:localglue} hold as identities of virtual characters of
$H$.  
\qed
\end{theorem}

The orbifold structures of $M_\P$ and $L_\P$ depend not only on $\P$,
but also on the sets of labels defining them; cf.\ Example
\ref{example:cut}.  However, by Proposition 4.4 of \cite{me:sym} the
equivariant character $\RR(M_\P,L_\P)$ depends only on the underlying
polyhedron, so \eqref{equation:glue} and \eqref{equation:localglue}
make sense.

\subsubsection{Multiplicities}\label{subsubsection:multiplicities}

If the origin is a vertex of the weight polytope $\Delta_L$, then the
right-hand side of \eqref{equation:fixed} has at most one nonzero
summand, which corresponds to a certain vertex of the moment polytope
$\Delta$.

\begin{proposition}\label{proposition:vertex}
Suppose that $0$ is a vertex of the weight polytope $\Delta_L$.
Choose a generic $\xi\in\lie t$ with the property that
$\langle\sigma_F,\xi\rangle>0$ for all $F$ with $\sigma_F\not=0$.  Let
$\nu$ be the vertex of the moment polytope $\Delta$ where the function
sending $\rho$ to $\langle\rho,\xi\rangle$ attains its minimum.  Then
\begin{equation}\label{equation:zeromult}
\RR(M,L)^T=
\begin{cases}
\RR\bigl(\Phi\inv(\nu),L|_{\Phi\inv(\nu)}\bigr)=\RR(M_\nu,L_\nu) &
\text{if $\sigma_{\Phi\inv(\nu)}=0$}, \\ 0 & \text{otherwise}.
\end{cases}
\end{equation}
If $L$ is rigid\upn, then $\RR(M,L)=\RR(M_\mu,L_\mu)$ for all vertices
$\mu$ of $\Delta$.  If $L$ is a moment bundle\upn, then
$\RR(M,L)^T=\RR(M_0,L_0)$ and
$$
\RR\bigl(M,L\inv\bigr)^T=
\begin{cases}
\RR\bigl(M_0,L_0\inv\bigr) & \text{if $\Delta=\{0\}$}, \\ 0 &
\text{otherwise}.
\end{cases}
$$
\end{proposition}

\begin{proof}
Note first that for every fixed-point component $F$ the cone $\ca C_F$
is equal up to a translation to the cone with vertex $\Phi_F$ spanned
by $\Delta$:
\begin{equation}\label{equation:twocones}
\ca C_F=-\Phi(F)+\cone_{\Phi(F)}(\Delta).
\end{equation}
In particular, $\ca C_F=\t^*$ for $\Phi(F)$ in the interior of
$\Delta$.  Secondly, if $\mu$ is any vertex of $\Delta$, then
$\Phi\inv(\mu)$ is a connected component of $M^T$, and so
$M_\mu=\Phi\inv(\mu)$.  Further, if $\xi$ is generic, then
$\rho\mapsto\langle\rho,\xi\rangle$ attains its minimum at a unique
vertex of $\Delta$, so $\nu$ is well-defined.  If
$\sigma_{\Phi\inv(\nu)}=0$, then $T$ acts trivially on
$L|_{\Phi\inv(\nu)}$ and so $L_\mu=L|_{\Phi\inv(\nu)}$.

If $F$ is a fixed-point component for which $\xi\in\check\ca C_F$,
then by \eqref{equation:twocones} the moment polytope is contained in
the halfspace given by
$\langle\rho,\xi\rangle\ge\langle\Phi(F),\xi\rangle$.  In other words,
$\rho\mapsto\langle\rho,\xi\rangle$ attains its minimum at $\Phi(F)$,
so $\Phi(F)=\nu$.  The equality \eqref{equation:zeromult} now follows
immediately from Theorem \ref{theorem:vertex} by setting $\mu=0$.

If $L$ is rigid, then $\Delta_L=\{0\}$, so $\sigma_{\Phi\inv(\nu)}=0$.
Moreover, $\RR(M,L)$ is a constant character by
\ref{part:almostcomplex} of Theorem \ref{theorem:rigid}, so
\eqref{equation:zeromult} implies
$\RR(M,L)=\RR(M,L)^T=\RR(M_\nu,L_\nu)$.  By varying the choice of
$\xi$ we obtain this equality for all vertices $\nu$.  If $L$ is a
moment bundle, then $\Delta_L=\Delta$, so $\nu=0$.  Hence
$\RR(M,L)^T=\RR(M_0,L_0)$ by \eqref{equation:zeromult}.  The weight
polytope of the dual moment bundle $L\inv$ is
$\Delta_{L\inv}=-\Delta$.  So if $\Delta=\{0\}$, then $\nu=0$ and
$\RR\bigl(M,L\inv\bigr)^T= \RR\bigl(M_0,L_0\inv\bigr)$, but if
$\Delta\ne\{0\}$, then $\nu$ is a vertex distinct from $0$, in which
case $\sigma_{\Phi\inv(\nu)}=-\nu\ne0$, so
$\RR\bigl(M,L\inv\bigr)^T=0$ by \eqref{equation:zeromult}.
\end{proof}

An immediate consequence of this result and the gluing formula is the
invariance of the index under blowing up.  Let $S$ be a closed
symplectic suborbifold of $M$ and suppose that $S^1$ acts on an open
neighbourhood $U$ of $S$ with fixed-point set $S$ and with positive
weights $\alpha$.  Let $\Bl(M,S,\eps,\alpha)$ be a weighted blowup of
$M$ at $S$ as defined at the end of Section \ref{subsubsection:blow}
and let $\Bl(L,S,\eps,\alpha)$ be the weighted blowup bundle.

\begin{theorem}\label{theorem:invariant}
\begin{equation}\label{equation:invariant}
\RR\bigl(\Bl(M,S,\eps,\alpha),\Bl(L,S,\eps,\alpha)\bigr)=\RR(M,L).
\end{equation}
\end{theorem}

\begin{proof}
The blowup is by definition the symplectic cut $M_{\le0}$ with respect
to the function $\psi+\eps$, where $\psi$ is the function generating
the circle action on $U$.  The exceptional divisor $M_0$ is the
weighted projectivization $\bb{P}_\alpha N$, where $N$ is the
symplectic normal bundle of $S$, and $M_{\ge0}$ is the weighted
projectivization $\bb{P}_\alpha(N\oplus\C)$.  By Definition
\ref{definition:lineblow}, $L_{\ge0}$ is $S^1$-rigid with respect to
the residual $S^1$-action on $M_{\ge0}$.  The minimum of the moment
function on $M_{\ge0}$ is $0$ and the fibre over $0$ is $M_0$.  Hence
$\RR(M_{\ge0},L_{\ge0})=\RR(M_0,L_0)$ by Proposition
\ref{proposition:vertex} and therefore
$\RR(M,L)=\RR(M_{\le0},L_{\le0})$ by \eqref{equation:localglue}.
\end{proof}

Let $H$ be a compact connected Lie group that acts on $M$ and $L$ in
such a way that $H$ commutes with $T$ and the $H$-action on $M$ is
symplectic.  The following assertion is evident from
\ref{part:tcrossh} of Theorem \ref{theorem:vertex}.

\begin{addendum}\label{addendum:tcrossh}
The equalities \eqref{equation:zeromult} and
\eqref{equation:invariant} hold as identities of virtual characters of
$H$.  
\qed
\end{addendum}

To put \eqref{equation:glue} and Proposition \ref{proposition:vertex}
together we need to calculate, for any admissible polytope $\P$, the
fixed points of the $T$-action on $M_\P$ and the weights of the
$T$-action on the fibres of $L_\P$ at $M_\P^T$.  Let $m$ be in the
stratum $M_\beta$ defined by \eqref{equation:stratification}, so that
$\lie t_m=\lie t_\beta$, and suppose that $\Phi(m)$ lies in an open
face $\F$ of $\P$.  Then by \eqref{equation:cut} the infinitesimal
stabilizer of the image of $m$ in $M_\P$ is equal to $\lie
t_\beta+\lie t_\F$.  By \eqref{equation:admissible} this sum is
direct. Therefore, the connected components of $M_\P^T$ are the
orbifolds
\begin{equation}\label{equation:cutfix}
\bigl(M_\beta\cap\Phi\inv(\F)\bigr)\big/T_\F
\end{equation}
for all $\beta\in\eu B$ and all open faces $\F$ of $\P$ such that
$\lie t_\beta\oplus\lie t_\F=\lie t$.  (These sets are connected,
because $M_\beta\cap\Phi\inv(\F)$ is exactly the open stratum of a
fibre of the $T_\F$-moment map on $\bar{M}_\beta$.)  Let
$\sigma_\beta\in\t_\beta^*$ denote the orbiweight of the $\lie
t_\beta$-action on $L|_{M_\beta}$ and $\bar\sigma_{\beta\F}$ the
orbiweight of the $\lie t$-action on the restriction of $L_\P$ to the
fixed-point component \eqref{equation:cutfix}.  Then
\eqref{equation:cutbundles} implies
\begin{equation}\label{equation:cutweight}
\bar\sigma_{\beta\F}= (\sigma_\beta,0)\in\lie t_\beta^*\oplus\lie
t_\F^*=\lie t^*,
\end{equation}
which proves the first part of the following lemma.

\begin{lemma}\label{lemma:cutweight}
\begin{enumerate}
\item\label{part:cutweight}
The orbiweight $\bar\sigma_{\beta\F}$ is the projection of
$\sigma_\beta$ onto the tangent space of $\F$ along the affine space
spanned by the virtual face $\Phi(M_\beta)$.
\item\label{part:aeltrigid}
If $L$ is almost equivariantly locally trivial at level $\mu$\upn,
then there exists a neighbourhood $O$ of $\mu$ in $\lie t^*$ such that
for every admissible polyhedron $\P\subset O$ the cut bundle $L_\P$ is
rigid on $M_\P$.
\item\label{part:cutrigid}
If $L$ is rigid\upn, then so is $L_\P$ for any admissible polyhedron
$\P$.
\item\label{part:cutmoment}
Let $\P$ be an admissible polyhedron and assume $L$ is a moment
bundle.  Then $L_\P$ is a moment bundle on $M_\P$ if and only if the
affine subspace spanned by $\F$ contains the origin for every open
face $\F$ of $\P$ such that $\Delta\cap\F$ is nonempty.
\item\label{part:shiftcone}
If $L$ is a moment bundle and $\P$ is an admissible cone with apex at
the origin\upn, then $L_\P$ is a moment bundle.  If $\P'$ is the
shifted cone $\mu+\P$\upn, where $\mu\in\t^*$\upn, then the weight
polytope of $L_{\P'}$ is contained in $\P$ when $\mu$ is sufficiently
small.  Hence\upn, by Theorem \ref{theorem:vertex}\upn, the support of
the multiplicity function of $L_{\P'}$ is contained in $\P$.
\end{enumerate}
\end{lemma}

\begin{proof}
If $L$ is almost equivariantly locally trivial at level $\mu$, then by
Lemma \ref{lemma:aelt} there exists a neighbourhood $O$ of $\mu$ in
$\lie t^*$ such that $L$ is almost equivariantly locally trivial on
$\Phi\inv(O)$.  This implies that if $\P$ lies in $O$, then
$\sigma_\beta=0$ for all $\beta$ and $\F$ such that
$M_\beta\cap\Phi\inv(\F)$ is nonempty. Therefore $L_\P$ is rigid by
\ref{part:cutweight}.  This proves \ref{part:aeltrigid}.

If $L$ is rigid, then $\sigma_\beta=0$ for all $\beta$ by Lemma
\ref{lemma:line} and hence $L_\P$ is rigid for any admissible $\P$ by
\ref{part:cutweight}.  This proves \ref{part:cutrigid}.

For the proof of \ref{part:cutmoment}, let $\iota_\F$ and
$\iota_\beta$ denote the inclusions of $\lie t_\F$, resp.\ $\lie
t_\beta$, into $\lie t$.  Since $L$ is a moment bundle,
$\sigma_\beta=\iota_\beta^*\Phi(m)$ by Lemma \ref{lemma:line}.  The
cut bundle $L_{\ca P}$ is therefore a moment bundle if and only if
$\sigma_\beta= \Phi(m)=
\bigl(\iota_\F^*\Phi(m),\iota_\beta^*\Phi(m)\bigr)$ for all open faces
$\F$ of $\P$ and all $m\in\Phi\inv(\F)$ such that $\lie t_m=\lie
t_\beta$.  This condition amounts to $\iota_\F^*\Phi(m)=0$ for all
$m\in\Phi\inv(\F)$, which is equivalent to $\F\subset(\lie t_\F)^0$
whenever $\Phi(M)\cap\F$ is nonempty.

The first statement under \ref{part:shiftcone} is evident from
\ref{part:cutmoment}.  For the second statement consider an arbitrary
orbiweight $\bar\sigma_{\beta\F'}$ of $L_{\P'}$.  By
\ref{part:cutweight} it is equal to the projection of $\Phi(m)$ onto
the tangent space of $\F'$ for some $m\in M_\beta$.  Here the open
face $\F'$ of $\P'$ and $\beta\in\eu B$ are such that $\F'$ and
$\Phi(M_\beta)$ intersect transversely at $\Phi(m)$.  Now $\F'$ is of
the form $\mu+\F$, where $\F$ is an open face of $\P$.  If $\F=\{0\}$
then $\bar\sigma_{\beta\F'}=0$, which is in $\P$.  If $\F\ne\{0\}$,
let us choose an inner product on $\lie t^*$ such that the
decomposition $\lie t^*=\lie t_\beta^*\oplus\lie t_\F^*$ is
orthogonal.  Then the distance $d$ of $\Phi(m)\in\F'$ to the boundary
of $\F'$ is positive.  This implies that as long as
$\lvert\mu\rvert<d$ the projection of $\Phi(m)$ onto $\F$ is contained
in the interior of $\F$, so in particular $\bar\sigma_{\beta\F'}$ is
in $\P$.  Since the number of orbiweights $\bar\sigma_{\beta\F'}$ is
finite, we get only finitely many such conditions on $\mu$ and
conclude that $\Delta_{L_{\P'}}\subset\P$ if $\mu$ is sufficiently
small.
\end{proof}

\begin{proof}[Proof of Theorem \ref{theorem:rigid} \upn(abelian
case\upn)] 
The proof of \ref{part:almostcomplex} is in Section
\ref{subsection:almostcomplex}.  For the proof of
\ref{part:hamiltonian} we may assume that $T$ acts effectively on $M$
and, after shifting the moment map if necessary, that $\mu=0$.  Let
$\Delta'$ be the unique closed face of $\Delta$ that contains $0$ in
its (relative) interior.  Then $\Delta'$ is the image under $\Phi$ of
a component $M'$ of the fixed-point set of a certain subtorus $T'$ of
$T$.  Now $\pr_{(\t')^*}(\Delta')$ is a vertex of
$\pr_{(\t')^*}(\Delta)$, which is the moment polytope of $M$ for the
$T'$-action, and the symplectic quotient of $M$ at
$\pr_{(\t')^*}(\Delta')$ with respect to $T'$ is $M'$.  Proposition
\ref{proposition:vertex} implies that $\RR(M,L)=\RR(M',L|_{M'})$.  We
may therefore assume that $0$ is in the interior of $\Delta$.

Consider first the case that $0$ is a quasi-regular value of $\Phi$.
Then $0$ is in fact a regular value, because every fibre over an
interior point of $\Delta$ intersects the principal stratum of $M$ and
$T$ acts generically freely.  This implies that there exists an
admissible polyhedral subdivision $\lie P$ of $\lie t^*$ such that
$0\in\lie P$.  (For instance, we can take each element of $\lie P$ to
be a suitable simplicial cone with vertex $0$.)  Then $L_\P$ is rigid
for every $\P\in\lie P$ by Lemma \ref{lemma:cutweight} and $0$ is a
vertex of every $\P$.  Consequently $\RR(M_\P,L_\P)=\RR(M_0,L_0)$ by
Proposition \ref{proposition:vertex}.  Combining the Euler identity
\begin{equation}\label{equation:euler}
\sum_{\P\in\lie P}(-1)^{\codim\P}=1
\end{equation}
with \eqref{equation:glue} we infer that $\RR(M,L)=\RR(M_0,L_0)$.

Consider finally the case that $0$ is a singular value of $\Phi$.
Since $T$ is abelian and $M$ is compact, the blowups at each stage of
the desingularization process of Section \ref{subsubsection:resolve}
are globally defined and give rise to a compact Hamiltonian
$T$-orbifold $\bigl(\tilde M,\tilde\omega,\tilde\Phi\bigr)$ and a
globally defined rigid orbibundle $\tilde L$ on $\tilde M$.  Theorem
\ref{theorem:invariant} implies that $\RR(M,L)=\RR\bigl(\tilde
M,\tilde L\bigr)$.  Because $0$ is a quasi-regular value of
$\tilde\Phi$, the $T$-equivariant part of the latter is equal to
$\RR\bigl(\tilde M_0,\tilde L_0\bigr)$, which is equal to
$\RR(M_0,L_0)$ by Definition \ref{definition:singular}.
\end{proof} 

\begin{proof}[Proof of Theorem \ref{theorem:independent} \upn(abelian
case\upn)]
By Lemma \ref{lemma:cutweight} there exists a neighbourhood $O$ of
$\mu$ such that $L_\P$ is rigid on $M_\P$ for all admissible
$\P\subset O$.  Let $\nu$ be any point in $O$.  Choose any set of
labels $\S=\bigl\{(v_1,r_1),(v_2,r_2),\dots,(v_n,r_n)\bigr\}$ such
that the associated polyhedron $\P$ has dimension $\dim T$, is
contained in $O$, and contains $\mu$ and $\nu$ in its interior.  By
Lemma \ref{lemma:transverse} admissibility is a generic condition, so
a small perturbation of the parameters $(r_1,\dots,r_n)$ will change
$\P$ into an admissible polytope that still satisfies $\mu$,
$\nu\in\P$ and $\P\subset O$.  Then $L_\P$ is rigid on $M_\P$, so
Theorem \ref{theorem:rigid} implies $\RR(M_\mu,L_\mu)= \RR(M_\P,L_\P)=
\RR(M_\nu,L_\nu)$.
\end{proof}

\begin{proof}[Proof of Theorem \ref{theorem:moment} \upn(abelian
case\upn)]
As in the proof of Theorem \ref{theorem:rigid} we can reduce the
general case to the case where $T$ acts effectively and that $0$ is an
interior point of $\Delta$.  (Here we apply Proposition
\ref{proposition:vertex} and Addendum \ref{addendum:tcrossh}, where
$H=T/T'$.)

Again, we handle first the case that $0$ is a regular value of $\Phi$.
Let $\lie P$ be an admissible polyhedral subdivision of $\lie t^*$
consisting of cones centred at the origin.  Then $L_\P$ is a moment
bundle on $M_\P$ for every $\P\in\lie P$ by Lemma
\ref{lemma:cutweight}.  By Proposition \ref{proposition:vertex},
\begin{equation}\label{equation:conecut}
\RR(M_\P,L_\P)^T=\RR(M_0,L_0),
\end{equation}
because $0$ is a vertex of $\Delta\cap\P$.  Applying the gluing
formula and \eqref{equation:euler} we find $\RR(M,L)^T=\RR(M_0,L_0)$.
  
If $0$ is a singular value, we can choose $\lie P$ such that the
shifted cones $\P'=\mu+\P$ for $\P\in\lie P$ are admissible for $\mu$
sufficiently close to $0$.  The weight for the $T$-action on
$L_\mu=L_{\P'}|_{\Phi_{\P'}\inv(\mu)}$ is trivial and the set of
weights for the $T$-action on $L_{\P'}|_{M_{\P'}^T}$ is contained in
$\P$ by Lemma \ref{lemma:cutweight}.  Hence
\begin{equation}\label{equation:shiftconecut}
\RR(M_{\P'},L_{\P'})^T=\RR(M_\mu,L_\mu)
\end{equation}
by Proposition \ref{proposition:vertex}.  Putting together the gluing
formula, \eqref{equation:euler} and Theorem \ref{theorem:independent}
we obtain $\RR(M,L)^T=\RR(M_0,L_0)$.
\end{proof}

\begin{proof}[Proof of Theorem \ref{theorem:dual} \upn(abelian
case\upn)]
Assuming that the theorem is true, by taking $T$-invariants on both
sides we deduce that
\begin{equation}\label{equation:dual}
N_{L\inv}(0)=
\begin{cases}
(-1)^{\dim\Delta}\RR\bigl(M_0,L_0\inv\bigr) & \text{if
$0\in\inter\Delta$}, \\ 0 & \text{otherwise}.
\end{cases}
\end{equation}
We assert that the theorem is in fact equivalent to
\eqref{equation:dual}.  This follows from a variant of the shifting
trick, which in this abelian situation allows us to write
\begin{multline}\label{equation:dualshift}
N_{L\inv}(\mu)=
\RR\bigl(M\times\{\mu\},L\inv\boxtimes(E_{-\mu})\inv\bigr)^T=
\RR\bigl(M\times\{\mu\},(L\boxtimes E_{-\mu})\inv\bigr)^T \\ =
\begin{cases}
(-1)^{\dim\Delta}\RR\bigl(M_{-\mu},(L_{-\mu}\shift)\inv\bigr) &
\text{if $-\mu\in\inter\Delta$}, \\ 0 & \text{otherwise}.
\end{cases}
\end{multline}
The first equality follows from the K\"unneth formula and the fact
that the character $\zeta_{-\mu}$ is dual to $\zeta_\mu$; and the
third equality follows from \eqref{equation:dual}.  It is clear that
\eqref{equation:dualshift} implies Theorem \ref{theorem:dual}.

The proof of \eqref{equation:dual} proceeds in the same way as the
proof of Theorem \ref{theorem:moment}.  The only difference is that
\eqref{equation:conecut} is replaced by
$$
\RR\bigl(M_\P,L_\P\inv\bigr)^T=
\begin{cases}
\RR\bigl(M_0,L_0\inv\bigr) &
\text{if $\P=\{0\}$}, \\ 0 & \text{otherwise},
\end{cases}
$$
and \eqref{equation:shiftconecut} by
$$
\RR\bigl(M_{\P'},L_{\P'}\inv\bigr)^T= 
\begin{cases}
\RR\bigl(M_\mu,L_\mu\inv\bigr) & \text{if $\P=\{0\}$}, \\ 0 &
\text{otherwise},
\end{cases}
$$
both of which follow from Proposition \ref{proposition:vertex}.
\end{proof}

\subsubsection{Asymptotic moment bundles}
\label{subsubsection:asymptotic}

The above proofs of Theorems \ref{theorem:moment} and
\ref{theorem:dual} do not generalize directly to the nonabelian case.
We need to extend the discussion to a class of orbibundles that are
``almost'' moment bundles.  Let $\eu M$ be a $T$-orbifold fibring over
the interval $(0,1]$.  Assume that the bundle projection has compact
fibres and is $T$-invariant.  Let $\eu L$ be a line orbibundle over
$\eu M$ and $(\omega,\Phi)$ a relative equivariant symplectic form on
$\eu M$ in the sense of Appendix \ref{subsection:darboux}.  Let us
denote the fibre of $\eu M$ over $t$ by $M^t$ and the restrictions of
$\eu L$, $\omega$ and $\Phi$ to $M^t$ by $L^t$, $\omega^t$ and
$\Phi^t$, respectively.  Then all fibres $M^t$ are Hamiltonian
$T$-orbifolds; they are equivariantly diffeomorphic (but not
necessarily symplectomorphic) to one another; and the orbibundles
$L^t$ are equivariantly isomorphic to one another.  It follows that
the $T$-character $\RR\bigl(M^t,L^t\bigr)$ is independent of $t$.

We call $\eu L$ an \emph{asymptotic moment bundle\/} if for all
components $\eu F$ of the fixed-point set $\eu M^T$ the limit of
$\Phi(\eu F\cap M^t)$ as $t\to0$ exists and is equal to the orbiweight
of the $T$-action on $\eu L|_{\eu F}$.  Here $\Phi(\eu F\cap M^t)$ is
the (constant) value of $\Phi$ on $\eu F\cap M^t$.  Let $\Delta^t$
denote the moment polytope of $M^t$; then the limit polytope
$\Delta=\lim_{t\to0}\Delta^t$ is well-defined.  As in Lemma
\ref{lemma:line} one shows that $\sigma_m=
\lim_{t\to0}\iota_m^*\Phi^t(m)$, where $\sigma\in\lie t_m^*$ is the
character defining the action of $(T_m)^0$ on $L_m$ and
$\iota_m\colon\lie t_m\to\lie t$ is the inclusion.  

Let us now for each $\mu\in\Delta$ select a path $\gamma_\mu(t)$ in
$\lie t^*$ defined for $0<t\le1$ such that $\gamma_\mu(t)\in\Delta^t$
and $\lim_{t\to0}\gamma_\mu(t)=\mu$.  Then for all $m\in M^t$ such
that $\Phi(m)=\gamma_\mu(t)$ we have
$\sigma_m=\lim_{t\to0}\gamma_\mu(t)=\mu$.  It follows that for all
$\mu$ and all $t$ the orbibundle $L^t\boxtimes E_{-\mu}$ is almost
equivariantly locally trivial at level $\gamma_\mu(t)$, so that
$$
(L^t)_{\gamma_\mu(t)}\shift= (L^t\boxtimes
E_{-\mu})\big|_{(\Phi^t)\inv(\gamma_\mu(t))}\Big/T
$$
is a well-defined orbibundle on $M^t_{\gamma_\mu(t)}$.  The proof of
the following result is completely analogous to that of Theorems
\ref{theorem:moment} and \ref{theorem:dual} in the abelian case.

\begin{theorem}\label{theorem:asymptotic}
Let $\eu L$ be an asymptotic moment bundle on $\eu M$.  Then for 
$0<t\le1$
\begin{align*}
\RR(M^t,L^t) &=
\sum_{\mu\in\Lambda^*\cap\Delta} 
\RR\bigl(M^t_{\gamma_\mu(t)},(L^t)_{\gamma_\mu(t)}\shift\bigr)
\,\zeta_\mu
\qquad\text{and}
\\
\RR\bigl(M^t,(L^t)\inv\bigr)
&=(-1)^{\dim\Delta}\sum_{\mu\in\Lambda^*\cap\inter\Delta}
\RR\Bigl(M^t_{\gamma_\mu(t)},
\bigl((L^t)_{\gamma_\mu(t)}\shift\bigr)\inv\Bigr)
\,\zeta_{-\mu}.
\qed
\end{align*}
\end{theorem}

\subsection{Delzant spaces III}\label{subsection:delzant3}

Let $\S=\bigl\{(v_1,r_1),(v_2,r_2),\dots,(v_n,r_n)\bigr\}$ be a set of
labels for the torus $T$.  Assume that the associated polyhedron $\P$
is nonempty and compact and that its dimension is equal to $k=\dim T$.
Suppose that the $r_i$ are integers, so that $\P$ is a rational
polyhedron.  We assert that the Delzant space $D_\S$ is prequantizable
in the sense of Example \ref{example:prequantum}.

Indeed, the vector $r=(r_1,\dots,r_n)\in(\Z^n)^*$ defines a real
infinitesimal weight of the torus $T^n=\R^n/\Z^n$.  This implies that
the cotangent bundle $T^*T$ is $T^n\times T$-equivariantly
prequantizable.  The prequantum line bundle is the trivial line bundle
$L_{T^*T}=T^*T\times\C$, where $T^n$ acts with weight $-r$ on the
fibre and $T$ acts trivially on the fibre.  By
\eqref{equation:kostant} the $T^n$-moment map corresponding to this
equivariant line bundle is the map $\psi_r$ given by
\eqref{equation:psi}.  A $T^n\times T$-equivariant prequantum line
bundle on $\C^n$ is the trivial line bundle $L_{\C^n}=\C^n\times\C$,
where $T^n$ and $T$ both act trivially on the fibre.  The bundle
$L=L_{T^*T}\boxtimes L_{\C^n}$ is then a $T^n\times T$-equivariant
prequantum line bundle on $T^*T\times\C^n$ and the associated moment
map for the $T^n$-action is given by \eqref{equation:tildepsi}.  The
upshot is that the quotient $L_\S=L\qu T^n$ is a $T$-equivariant
prequantum line orbibundle on the Delzant space $D_\S$ and that the
associated moment map is $\Psi_\S$.

If $\P$ is a lattice polytope (i.~e.\ all its vertices are in
$\Lambda^*$), then it is not difficult to see that $L_{T^*T}$ is
$T^n$-equivariantly locally trivial, so $L_0$ is in fact a genuine
line bundle.

By Proposition \ref{proposition:dilate}, for $m\in\Z$ the Delzant
space $D_{m\S}$ is symplectomorphic to $D_\S$ with $m$ times its
symplectic form.  It is not hard to check that under this
symplectomorphism $L_{m\S}$ pulls back to the $m$th tensor power
$L_\S^m$.  The next result follows immediately from Theorems
\ref{theorem:moment} and \ref{theorem:dual} and the fact that $D_\S$
is multiplicity-free.

\begin{proposition}\label{proposition:lattice}
For all $m\in\N$
\begin{align*}
\RR(D_\S,L_\S^m)(z) &=\sum_{\mu\in\Lambda^*\cap m\P}z^\mu \qquad\text{
and } \\
\RR(D_\S,L_\S^{-m})(z)
&=(-1)^k\sum_{\mu\in\Lambda^*\cap\inter(m\P)}z^{-\mu}
\end{align*}
as rational functions on $T^\C$.  
\qed
\end{proposition}

For every nonnegative integer $m$, let $p(m)$ denote the number of
lattice points in $m\P$.  Then by Proposition
\ref{proposition:lattice}, $\RR(D_\S,L_\S^m)(1)=p(m)$ and
$\RR(D_\S,L_\S^{-m})(1)$ is $(-1)^k$ times the number of lattice
points in the interior of $m\P$.  The following result was proved for
simple rational polytopes in \cite{me:on}.

\begin{corollary}
The counting function $p$ is a quasi-polynomial\upn, whose period is a
divisor of the smallest positive integer $l$ such that $l\P$ is a
lattice polytope.  It satisfies the Ehrhart reciprocity law
$p(-m)=(-1)^k\#\bigl(\Lambda^*\cap\inter(mP)\bigr)$.
\end{corollary}

\begin{proof}
We defined $D_\S$ as the symplectic quotient $M_0$ of
$M=T^*T\times\C^n$ and $L_\S$ as the quotient $L_0$ of the line bundle
$L=L_{T^*T}\boxtimes L_{\C^n}$.  Select a small regular value $\eta$
of $\tilde\psi_r$; then $\RR(D_\S,L_\S^m)=\RR(M_\eta,L_\eta^m)$ for
all $m\in\Z$ by Theorem \ref{theorem:independent}.  The result now
follows from Proposition \ref{proposition:lattice} and the
Kawasaki-Riemann-Roch formula for $M_\eta$ and $L_\eta^m$.
\end{proof}
 
This implies the well-known result that the counting function of a
lattice polytope is polynomial.  (See e.~g.\ \cite{fu:in}.)  Recall
from Section \ref{subsection:delzant2} that the shift
desingularization $M_\eta$ is none other than the Delzant space
associated to the labelled polytope $(\S_\eta,\P_\eta)$, where
$\S_\eta$ is as in \eqref{equation:shiftlabel}.  Guillemin pointed out
in \cite{gu:ri} that Proposition \ref{proposition:lattice} leads to an
Euler-MacLaurin type formula for the number of lattice points in $\P$,
namely
\begin{equation}\label{equation:khovanskii}
p(1)=\lim_{\eta\to0}\sum_{\F\pre\P}\lie
T_\F\biggl(\frac\partial{\partial\eta}\biggr)\vol\P_\eta.
\end{equation}
Here $\vol\P_\eta$ is the normalized Euclidean volume of $\P_\eta$ and
the $\lie T_\F$ are certain infinite-order differential operators
associated to the faces of $\P$.  These operators depend on the
component of the set of regular values containing $\eta$; see
\cite{gu:ri} for details.  There is an analogous formula involving the
moment polytope $\tilde\P_\eps$ associated to the canonical
desingularization of $D_\S$.  These identities are generalizations of
the Khovanskii-Pukhlikov formula for the number of lattice points in a
simply laced lattice polytope.  Purely combinatorial proofs were given
independently by Brion and Vergne \cite{br:res}.  See their paper for
a discussion of the relationship between \eqref{equation:khovanskii}
and similar formul\ae\ proven by Cappell and Shaneson \cite{ca:eu}.

\section{The general case}\label{section:nonabelian}

This section contains the proofs of the remaining theorems of Section
\ref{section:results} for general compact groups.  We reduce the
general case to the abelian case by means of the cross-section theorem
and local symplectic cutting with respect to certain subtori of the
maximal torus.  In this section $(M,\omega,\Phi)$ denotes a compact
connected Hamiltonian $G$-orbifold with moment polytope $\Delta$ and
$L$ denotes a $G$-equivariant line orbibundle on $M$.

\subsection{Induction and cutting}\label{subsection:induction}

In this section we prove Theorems \ref{theorem:independent} and
\ref{theorem:rigid}.  We start by observing that Theorem
\ref{theorem:rigid} is true at ``maximal'' values of $\Phi$.  Let
$\Phi_T={\pr}_{\t^*}\circ\Phi$ be the moment map for the $T$-action on
$M$.  Choose an invariant inner product on $\g^*$ and let
$\vert{\cdot}\vert$ denote the associated norm.

\begin{lemma}\label{lemma:norm} 
Let $\mu$ be a point in $\Delta$ of maximal norm.  Then
$\Phi\inv(\mu)$ is a component of the fixed-point set $M^{G_\mu}$ and
is equal to $\Phi_T\inv(\mu)$.  It follows that $\mu$ is a
quasi-regular value of $\Phi$.
\end{lemma}

\begin{proof}
Recall that the moment polytope for the $T$-action is equal to
$\Phi_T(M)=\hull\Phi(M^T)=\hull(\W\cdot\Delta)$.  This implies that
$\mu$ is a vertex of $\Phi_T(M)$ and therefore
$\Phi_T\inv(\mu)=\Phi\inv(\mu+\t^0)$ is a component of $M^T$.  Because
the norm on $\lie g^*$ is invariant we have
$\vert\Phi(x)\vert\le\vert\mu\vert$ for all $x\in M$, which for
$x\in\Phi\inv(\mu+\t^0)$ is only possible if $\Phi(x)=\mu$.  This
shows $\Phi\inv(\mu)=\Phi_T\inv(\mu)$.  Now $G_\mu$ is a connected
subgroup of maximal rank of $G$, so $\Phi\inv(\mu)$, being both
invariant under the $G_\mu$-action and fixed under the $T$-action, is
fixed under $G_\mu$.  Since all points in $\Phi\inv(\mu)$ are of the
same $G_\mu$-orbit type, $\mu$ is a quasi-regular value.
\end{proof}

Putting this together with Theorem \ref{theorem:rigid} (part
\ref{part:almostcomplex} of which we proved in Section
\ref{subsection:almostcomplex} and part \ref{part:hamiltonian} of
which we proved in the abelian case in Section
\ref{subsection:hamiltonian}) and Addendum \ref{addendum:tcrossh}, we
obtain the following.

\begin{lemma}\label{lemma:maxrigid}
Let $\mu$ be as in Lemma \ref{lemma:norm} and suppose $L$ is rigid.
Then $\RR(M,L)=\RR(M,L)^G=\RR(M_\mu,L_\mu)$.  Given another Lie group
$H$ that acts on $L$ and $M$ in such a way that the action commutes
with that of $G$ and the action on $M$ is symplectic\upn, this holds
as an equality of virtual characters of $H$.  
\qed
\end{lemma}

What is the relationship between the equivariant index of $M$ and that
of its cross-sections $Y_\sigma$?  Recall that for every open wall
$\sigma$ of the Weyl chamber $\tplus$ the induction map
$\Ind_{G_\sigma}^G$ is defined as the unique homomorphism $f\colon\Rep
G_\sigma\to\Rep G$ such that $\Ind_T^G=f\circ\Ind_T^{G_\sigma}$.  More
specifically, let $\mu\in\Lambda^*_{\sigma,+}\subset\Lambda^*$ be a
dominant weight for $G_\sigma$ and let $\chi_{\sigma,\mu}\in\Rep
G_\sigma$ be the corresponding irreducible character.  Then
$$
\Ind_{G_\sigma}^G\chi_{\sigma,\mu}=\Ind_T^G\zeta_\mu,
$$
which is also equal to $\RR(G\mu,E_\mu)$ by \eqref{equation:bwb}.

Let $Y_\sigma$ be the cross-section of $M$ over $\sigma$ as defined in
Section \ref{subsection:crosssection}.  We say that $Y_\sigma$ is a
\emph{global\/} cross-section of $M$ if $M=GY_\sigma$, or
equivalently, $\Delta$ is a subset of the open star
$\bigcup_{\tau\suc\sigma}\tau$ of $\sigma$.

\begin{proposition}\label{proposition:crosssection}
\begin{enumerate}
\item\label{part:induction}
Let $\ca Y_\sigma$ be a compact almost complex $G_\sigma$-orbifold and
$\ca L_\sigma$ a $G_\sigma$-equivariant line orbibundle on $\ca
Y_\sigma$.  Let $\ca M$ be the almost complex $G$-orbifold
$G\times^{G_\sigma}\ca Y_\sigma$ equipped with the $G$-equivariant
orbibundle $\ca L$ induced by $\ca L_\sigma$.  Then $\RR(\ca M,\ca
L)=\Ind_{G_\sigma}^G\RR(\ca Y_\sigma,\ca L_\sigma)$.
\item\label{part:globalcross}
If $Y_\sigma$ is a global cross-section of $M$\upn, then
$\RR(M,L)=\Ind_{G_\sigma}^G\RR(Y_\sigma,L_\sigma)$.
\end{enumerate}
\end{proposition}

\begin{proof}
The proof of \ref{part:induction} is closely analogous to the proof of
the quantum cross-section theorem of \cite{me:sym}.  If $Y_\sigma$ is
a global cross-section, then $M=G\times^{G_\sigma}Y_\sigma$ by Theorem
\ref{theorem:crosssection}, so \ref{part:globalcross} is evident from
\ref{part:induction}.
\end{proof}

A global cross-section $Y_\sigma$ is \emph{nontrivial\/} if
$Y_\sigma\ne M$, that is to say $\sigma\ne\lie a^*$.  Usually $M$ does
not possess nontrivial global cross-sections, but even then we can
obtain information on its equivariant index by dint of nonabelian
symplectic cutting, which was invented by Woodward \cite{wo:cl}.  It
is based on the fact that $M$ is the union
$\bigcup_{\sigma\pre\tplus}M_\sigma$ of $G$-invariant open subsets
$M_\sigma=GY_\sigma$, each of which carries a Hamiltonian action of
the torus $A_\sigma$ which commutes with the action of $G$.  This
action is defined by identifing $M_\sigma$ with
$G\times^{G_\sigma}Y_\sigma$ as in the symplectic cross-section
theorem and extending the natural $A_\sigma$-action on $Y_\sigma$ to
an action on $M_\sigma$ which commutes with $G$.  In other words, for
a $G_\sigma$-orbit $[g,y]$ in $G\times^{G_\sigma}Y_\sigma$ and $t\in
A_\sigma$ we put $t\cdot[g,y]=[g,ty]$.  This is well-defined because
$A_\sigma$ commutes with $G_\sigma$.  The restriction of $L$ to
$M_\sigma$ acquires likewise a $A_\sigma$-action that commutes with
$G$.  The moment map for $A_\sigma$ is the unique $G$-invariant
extension of the $A_\sigma$-moment map on $Y_\sigma$ and can be
described as follows.  Let
$$
\Phi_+\colon M\to\tplus
$$
be the composition of $\Phi$ with the quotient mapping $q$ defined in
\eqref{equation:q}.  The moment map of the $A_\sigma$-action on
$Y_\sigma$ is equal to $\pr_\sigma\circ\Phi|_{Y_\sigma}$, where
$\pr_\sigma$ is the canonical projection $\lie t^*\to\lie a_\sigma^*$.
Now observe that $\pr_\sigma=\pr_\sigma\circ q$ on $\lie
g_\sigma^*\subset\lie g^*$, and hence $\pr_\sigma\circ\Phi=
\pr_\sigma\circ\Phi_+$ on $Y_\sigma$.  The $A_\sigma$-moment map on
$M_\sigma$ is therefore equal to the $G$-invariant map
$\pr_\sigma\circ\Phi_+|_{M_\sigma}$.  This is a smooth map for all
$\sigma$, even though $\Phi_+$ is in general not smooth.

Now let $\S$ be a set of labels in $\lie t^*$ and $\P$ its associated
polyhedron.

\begin{definition}\label{definition:gadmissible}
The pair $(\S,\P)$ is \emph{admissible\/} or \emph{$G$-admissible\/}
with respect to $M$ if $\S$ has constant excess and the following
conditions hold for all open faces $\F$ of $\P$:
\begin{enumerate}
\item\label{part:subtorus}
for all walls $\sigma$ such that $\sigma\cap\F\cap\Delta$ is nonempty,
$T_\F$ is a subtorus of $A_\sigma$;
\item\label{part:locallyfree}
the action of $T_\F$ on $\Phi\inv(\F\cap\tplus)$ is locally free.
\end{enumerate}
\end{definition}

As in the abelian case, admissibility depends only on the polyhedron
$\P$, not on $\S$, and condition \ref{part:locallyfree} is satisfied
generically.  Condition \ref{part:subtorus} is tantamount to: for all
$\sigma$ such that $\sigma\cap\F\cap\Delta\ne\emptyset$, the tangent
space to $\F$ contains the annihilator of $\lie a_\sigma$ in $\lie
t^*$; in other words the orthogonal complement of $\sigma$ (with
respect to any invariant inner product) is contained in $\F$.  It
implies that every wall $\sigma$ has an open neighbourhood $O_\sigma$
inside $\star\sigma$ such that
$$
\P\cap O_\sigma\cap\Delta= \pr_\sigma\inv(\P_\sigma)\cap
O_\sigma\cap\Delta,
$$
where $\P_\sigma=\P\cap\a_\sigma^*$.  The symplectic cut
$(M_\sigma')_{\P_\sigma}$ of the $G\times A_\sigma$-invariant open
subset $M_\sigma'=\Phi\inv(GO_\sigma)$ of $M_\sigma$ with respect to
the polyhedron $\P_\sigma$ is then well-defined and condition
\ref{part:locallyfree} implies that it is an orbifold.  For
$\sigma\pre\tau$ there is a natural symplectic embedding of a
$G$-invariant open subset of $(M_\sigma')_{\P_\sigma}$ into
$(M_\tau)_{\P_\tau}$ and the result of gluing the
$(M_\sigma')_{\P_\sigma}$ together along these embeddings is a compact
Hamiltonian $G$-orbifold $(M_\P,\omega_\P,\Phi_\P)$, the
\emph{symplectic cut\/} of $M$ with respect to $\P$.  Its moment
polytope $\Delta_\P$ is equal to $\Delta\cap\P$.  The bundles
$(L|_{M_\sigma'})_{P_\sigma}$ are likewise well-defined and can be
pasted together to a global $G$-equivariant \emph{cut bundle\/} $L_\P$
on $M_\P$.  See \cite{me:sym} for details.  Put
$\Phi_{\P,+}=q\circ\Phi_\P$.  By analogy with \eqref{equation:cut} and
\eqref{equation:cutbundles}, for every open face $\F$ of $\P$ there
are canonical isomorphisms
\begin{align*}
\Phi_{\P,+}\inv(\F) &\cong \Phi_+\inv(\F)/T_\F,\\
L_\P\big|_{\Phi_{\P,+}\inv(\F)} &\cong
\Bigl(L\big|_{\Phi_+\inv(\F)}\Bigr)\Big/T_\F.
\end{align*}

\begin{proof}[Proof of Theorem \ref{theorem:rigid}]
The proof of \ref{part:almostcomplex} is in Section
\ref{subsection:almostcomplex}.  For the proof of
\ref{part:hamiltonian} we consider first the case that $\Delta$ is
contained in the degenerate wall $\lie a^*$ of $\tplus$, where $\lie
a$ is the centre of $\lie g$.  Then $M$ is in effect a Hamiltonian
$A$-orbifold and the theorem reduces to the abelian case, which was
covered in Section \ref{subsubsection:multiplicities}.

Now consider the case that $\Delta$ is not contained in $\lie a^*$.
Here the proof is by induction on the dimension of $M$.  We may assume
that the result holds for all compact connected groups $H$ and all
Hamiltonian $H$-orbifolds $Q$ with $\dim Q<\dim M$.  By Lemma
\ref{lemma:maxrigid}, $\RR(M,L)=\RR(M,L)^G= \RR(M_\mu,L_\mu)$ if
$\mu\in\Delta$ is of maximal norm.  It therefore suffices to check
that $\RR(M_\mu,L_\mu)$ is independent of $\mu\in\Delta$.

First we show that $\RR(M_\mu,L_\mu)$ is constant on the complement in
$\Delta$ of $\lie a^*$.  Let $\mu$ and $\nu$ be in $\Delta-\lie a^*$
and let $\sigma$ be the largest open wall of $\tplus$ such that $\mu$,
$\nu\in\star\sigma$.  Then $\sigma\ne\lie a^*$ and hence $\dim
Y_\sigma<\dim M$.  Choose an admissible polytope $\P$ such that
$\P\cap\tplus$ is a subset of the star of $\sigma$ and $\mu$ and $\nu$
are in $\P$.  Put $Y_{\sigma,\P}=(Y_\sigma)_\P$ and
$L_{\sigma,\P}=(L|_{Y_\sigma})_\P$.  By the induction hypothesis the
function that sends $\lambda$ to
$\RR\bigl((Y_{\sigma,\P})_\lambda,(L_{\sigma,\P})_\lambda\bigr)$ is
constant on $\P\cap\Delta$.  Moreover,
$M_\lambda=(Y_{\sigma,\P})_\lambda$ and
$L_\lambda=(L_{\sigma,\P})_\lambda$ for all $\lambda$ in
$\P\cap\Delta$, so the conclusion is
$\RR(M_\mu,L_\mu)=\RR(M_\nu,L_\nu)$.

It remains to show that $\RR(M_\nu,L_\nu)=\RR(M_\mu,L_\mu)$ where
$\mu\in\Delta\cap\a^*$ and $\nu$ is a point close to $\mu$ and
contained in the principal face $\Delta\gen$ of $\Delta$.  Because
$\lie a$ is the centre of $\lie g$, we can shift the moment map by
$\mu$ and may therefore assume that $\mu=0$.  

Assume that $0$ is a quasi-regular value of $\Phi$, so that
$\Phi\inv(0)=Z_\alpha$ for some $\alpha\in\eu A$.  By Proposition
\ref{proposition:varquot}, $M_\nu$ is a symplectic fibre orbibundle
over $M_0$ with general fibre $(F_\alpha)_\nu$ and $L_\nu$ is the
pullback of $L_0$.  Because $\nu$ is a generic value of $\Phi$, it is
a quasi-regular value and $(F_\alpha)_\nu$ is an orbifold.  By Theorem
\ref{theorem:product} (see Appendix \ref{section:product}) we have
$\RR(M_\nu,L_\nu)= \RR(M_0,L_0)\RR\bigl((F_\alpha)_\nu,\C\bigr)$.
Here $F_\alpha=(T^*G\times W_\alpha/\Upsilon_\alpha)\qu G_\alpha$, so
$F_\alpha\qu G=(W_\alpha/\Upsilon_\alpha)\qu G_\alpha$ is a point by
Lemma \ref{lemma:w}.  Moreover, $\dim F_\alpha\le\dim M_\nu<\dim M$,
so $\RR\bigl((F_\alpha)_\nu,\C\bigr)=
\RR\bigl((F_\alpha)_0,\C\bigr)=1$ by the induction hypothesis.  The
upshot is $\RR(M_\nu,L_\nu)=\RR(M_0,L_0)$.

If $0$ is not a quasi-regular value, consider the blowup $\bigl(\tilde
U,\tilde\omega,\tilde\Phi\bigr)$.  As noted in the proof of Theorem
\ref{theorem:shiftcanonical}, for a suitable choice of the blowup
parameters and a sufficiently small quasi-regular value $\nu$ of
$\Phi$, $\tilde U_\nu$ is symplectomorphic to $M_\nu$ and moreover
$\tilde L_\nu$ is isomorphic to the pullback of $L_\nu$.  As $0$ is a
quasi-regular value of $\tilde\Phi$, we have $\RR(M_\nu,L_\nu)=
\RR\bigl(\tilde U_0,\tilde L_0\bigr)$, which is by definition equal to
$\RR(M_0,L_0)$.
\end{proof}

\begin{proof}[Proof of Theorem \ref{theorem:independent}]
This follows from Theorem \ref{theorem:rigid} just as in the abelian
case.
\end{proof}

\subsection{Multiplicities}\label{subsection:nonabelianmultiplicities}

This section contains the proofs of Theorems \ref{theorem:moment} and
\ref{theorem:dual}.  The main ingredient is the nonabelian gluing
formula.

\begin{definition}\label{definition:gsubdivision}
An \emph{admissible\/} or \emph{$G$-admissible polyhedral
subdivision\/} of $\tplus$ is a collection $\lie P$ satisfying the
following conditions: every element of $\lie P$ is a $G$-admissible
polyhedron in $\t^*$, their union contains $\tplus$, for every element
of $\lie P$ all its closed faces are in $\lie P$, and the intersection
of any two elements of $\lie P$ is a closed face of each.
\end{definition}

\begin{theorem}[gluing formula,
\cite{me:sym}]\label{theorem:nonabglue}
Let $\lie P$ be an admissible polyhedral subdivision of $\tplus$.
Then
\begin{equation}\label{equation:nonabglue}
\RR(M,L)=\sum_{\P\in\lie P}(-1)^{\codim\P}\RR(M_\P,L_\P)
\end{equation}
as virtual characters of $G$.
\qed
\end{theorem}

\begin{figure}
\setlength{\unitlength}{0.007mm}
$$
\begin{picture}(10223,8766)
\texture{cccccccc 0 0 0 cccccccc 0 0 0 
        cccccccc 0 0 0 cccccccc 0 0 0 
        cccccccc 0 0 0 cccccccc 0 0 0 
        cccccccc 0 0 0 cccccccc 0 0 0 }
\thinlines
\shade\path(2400,7530)
        (3643,7668)
        (4824,7866)
        (6004,8074)
        (7235,8237)
        (8386,8355)
        (9000,8400)
        (7980,7929)
        (6642,7401)
        (6029,7094)
        (5400,6750)
        (4244,6998)
        (2998,7417)
        (2721,7483)
        (2400,7530)
\path(600,300)(600,7500)
\path(5400,6300)(600,300)
\path(10200,8700)(600,300)
\path(600,7500)(729,7412)(815,7372)(930,7335)(1032,7335)(1113,7343)
(1200,7350)
\path(1200,7350)(1992,7302)(2490,7212)
\path(2490,7212)(3386,6807)(4399,6539)(4875,6431.25)
\path(4875,6431.25)(5181,6420)(5400,6300)
\path(600,7500)(1200,7652)(1515,7659)(1800,7650)
\path(1800,7650)
        (2820,7769)
        (3796,7907)
        (4780,8071)
        (5812,8240)
        (6821,8361)
        (7803,8449)
        (8811,8537)
        (8900,8537.5)
\path(8900,8537.5)(9507,8640)(9818,8670)(10200,8700)
\path(5400,6300)(5695,6410)(6000,6600)
\path(6000,6600)
        (6386,6863)
        (6826,7098)
        (7182,7279)
        (7820,7564)
        (8381,7797)
        (9011,8085)
        (9174,8187)
        (9500,8350)
\path(9500,8350)(9827,8565)(10032,8650)(10200,8700)
\dottedline{100}(600,1950)(1200,1800)(1500,1425)
\dottedline{100}(1200,7350)(1200,1800)(4875,6431.25)
\dottedline{100}(1500,1425)(2700,2400)(2700,2137.5)
\dottedline{100}(6000,6600)(2700,2400)(9500,8350)
\dottedline{100}(4875,6431.25)(5200,6600)(5400,6750)(5700,6750)
(6000,6600)
\dottedline{100}(8900,8537.5)(9000,8400)(9500,8350)
\dottedline{100}(1200,7350)(1800,7497)(2400,7530)
\dottedline{100}(2400,7530)(2123,7553)(1800,7650)
\put(600,300){\circle*{100}} 
\put(3000,3300){\circle*{100}} 
\put(600,7500){\circle*{100}} 
\put(5400,4500){\circle*{100}} 
\put(300,0){\makebox(0,0)[lb]{\smash{$0$}}}
\put(0,7500){\makebox(0,0)[lb]{\smash{$\lambda_2$}}}
\put(2500,3300){\makebox(0,0)[lb]{\smash{$\lambda_1$}}}
\put(5700,4400){\makebox(0,0)[lb]{\smash{$\lambda_3$}}}
\put(5100,3000){\makebox(0,0)[lb]{\smash{$\tplus$}}}
\end{picture}
$$
\caption{Admissible subdivision for $\SU(4)$, intersected with Weyl
chamber.  Shaded area represents shifted Weyl chamber}
\label{diagram:admissible}
\end{figure}
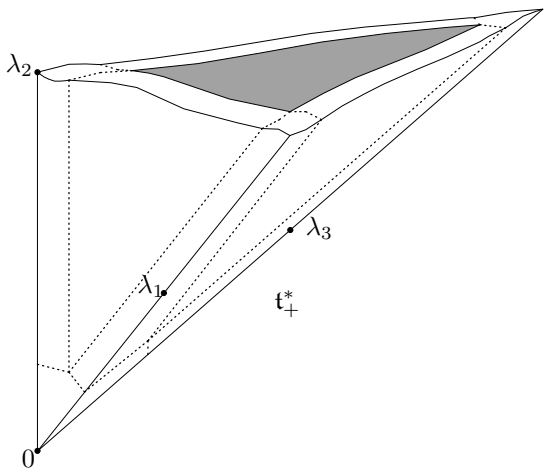

An example of an admissible polyhedral subdivision of $\tplus$ is the
subdivision that is dual to the decomposition into walls, which can be
described as follows.  For $\sigma\pre\tau$ define the polyhedral cone
$\ca C_{\sigma\tau}$ in $\lie t^*$ to be the product of $\lie a^*$ and
the cone in $[\lie g,\lie g]^*$ spanned by the vectors
\begin{equation}\label{equation:ray}
-\alpha_{j_1},-\alpha_{j_2},\dots,-\alpha_{j_r}\quad\text{and}\quad
\lambda_{i_1},\lambda_{i_2},\dots,\lambda_{i_s}.
\end{equation}
Here $r=\codim\sigma$, $s=\dim\bigl(\tau\cap[\lie g,\lie g]^*\bigr)$,
$\alpha_{j_1}$, $\alpha_{j_2},\dots$, $\alpha_{j_r}$ are the positive
simple roots perpendicular to $\sigma\cap[\lie g,\lie g]^*$, and
$\lambda_{i_1}$, $\lambda_{i_2},\dots$, $\lambda_{i_s}$ are the
fundamental weights spanning the wall $\tau\cap[\lie g,\lie g]^*$.
Hence
\begin{equation}\label{equation:dimension}
\codim\ca C_{\sigma\tau}=\dim\tau-\dim\sigma.
\end{equation}
Now choose $\lambda$ in the interior of the Weyl chamber, let
$\P_{\sigma\tau}$ be the shifted cone $\lambda+\ca C_{\sigma\tau}$,
and let $\lie P_\lambda$ be the collection of all $\P_{\sigma\tau}$
(see Diagram \ref{diagram:admissible}):
$$
\lie P_\lambda= \{\,\P_{\sigma\tau}:\sigma\pre\tau\pre\tplus\,\}.
$$
 
\begin{lemma}\label{lemma:admissible}
For generic $\lambda$ in $\inter\tplus$\upn, $\lie P_\lambda$ is an
admissible polyhedral subdivision of $\tplus$.
\end{lemma}

\begin{proof}[Outline of proof]
The closed faces of $\P_{\sigma\tau}$ are the $\P_{\rho\upsilon}$ with
$\rho\pre\sigma$ and $\tau\pre\upsilon$, so $\lie P_\lambda$ is closed
under inclusion of faces.

The intersection of $\P_{\sigma\tau}$ and $\P_{\rho\upsilon}$ is
$$
\P_{\sigma\tau}\cap\P_{\rho\upsilon}=
\P_{\sigma\wedge\rho,\tau\vee\upsilon},
$$
where $\sigma\wedge\rho$ is the largest open face contained in
$\bar\sigma\cap\bar\rho$ and $\tau\vee\upsilon$ is the smallest open
face that contains $\tau$ and $\upsilon$ in its closure.  This implies
that $\lie P_\lambda$ is closed under taking intersections.

Every wall $\sigma$ is contained in the union of all
$\P_{\upsilon\tau}$ with $\upsilon\pre\tau\pre\sigma$, so $\lie
P_\lambda$ covers $\tplus$.  In fact,
\begin{equation}\label{equation:intersection}
\sigma\cap\P_{\upsilon\tau}\ne\emptyset\quad\text{ if and only if }
\quad \upsilon\pre\tau\pre\sigma.
\end{equation}

It follows from \eqref{equation:ray} that the tangent space to
$\P_{\upsilon\tau}$ contains the annihilator of $\lie a_\upsilon$.
Therefore, if $\upsilon\pre\tau\pre\sigma$, then the tangent space to
$\P_{\upsilon\tau}$ contains $\lie a_\sigma^0$.  We conclude from
\eqref{equation:intersection} that all polyhedra in $\lie P_\lambda$
satisfy condition \ref{part:subtorus} of Definition
\ref{definition:gadmissible} for all $M$.  Furthermore,
$\upsilon\pre\tau$ implies that the sets $\{i_1,i_2,\dots,i_r\}$ and
$\{j_1,j_2,\dots,j_s\}$ are disjoint, so the vectors
\eqref{equation:ray} are linearly independent, and consequently any
(minimal) set of labels defining $\P_{\upsilon\tau}$ has constant
excess.  Because condition \ref{part:locallyfree} of Definition
\ref{definition:gadmissible} is satisfied generically, we conclude
that $\lie P_\lambda$ is an admissible polyhedral subdivision of
$\tplus$ for generic values of $\lambda$.
\end{proof}

From \eqref{equation:intersection} we obtain
$\tplus\cap\P_{\sigma\tau}\subset\star\tau$ whenever $\sigma\pre\tau$.
This means that for all $\tau\not=\a^*$ the symplectic cut
$M_{\P_{\sigma\tau}}$ possesses a nontrivial global cross-section,
namely the symplectic cut $Y_{\sigma\tau}= (Y_\tau)_{\P_{\sigma\tau}}$
of $Y_\tau$.  Let $L_{\sigma\tau}=(L|_{Y_\sigma})_{\P_{\sigma\tau}}$
denote the corresponding cut bundle.  From Proposition
\ref{proposition:crosssection}, the gluing formula and
\eqref{equation:dimension} we obtain
\begin{equation}\label{equation:subdivision} 
\RR(M,L)= \sum_{\substack{\sigma,\tau\\
\sigma\pre\tau\pre\tplus}}(-1)^{\dim\tau-\dim\sigma} 
\Ind_{G_\tau}^G\RR(Y_{\sigma\tau},L_{\sigma\tau}).
\end{equation} 
Henceforth let $\upsilon$ denote the principal wall for $M$ as defined
in Section \ref{subsection:crosssection}, so that
$\Delta\subset\bar\upsilon$.  The \emph{principal cone\/} in $\lie
P_\lambda$ is $\P_{\upsilon\upsilon}$ and the \emph{principal cut\/}
of $M$ is $M_{\P_{\upsilon\upsilon}}$.  The global cross-section
$Y_{\upsilon\upsilon}$ of $M_{\P_{\upsilon\upsilon}}$ is a Hamiltonian
$G_\upsilon$-orbifold.  Recall that $[G_\upsilon,G_\upsilon]$ acts
trivially on $Y_{\upsilon\upsilon}$, so that
$\RR(Y_{\upsilon\upsilon})$ is a trivial
$[G_\upsilon,G_\upsilon]$-character and hence
$$
\Ind_{G_\upsilon}^G\RR(Y_{\upsilon\upsilon},L_{\upsilon\upsilon})=
\Ind_T^G\RR(Y_{\upsilon\upsilon},L_{\upsilon\upsilon}).
$$
Note further that all terms in \eqref{equation:subdivision} vanish
except those for which $\sigma\pre\tau\pre\upsilon$.  For certain line
orbibundles the pair $(Y_{\upsilon\upsilon},L_{\upsilon\upsilon})$
captures all the information needed to compute $\RR(M,L)$.

\begin{theorem}[abelianization]\label{theorem:abelianization}
Suppose that $L$ has the property that for all $m\in M$ the action of
the identity component of $G_m\cap\bigl[G_{\Phi(m)},G_{\Phi(m)}\bigr]$
on the fibre $L_m$ is trivial.  Let $\upsilon$ be the principal wall
for $M$.  If $\lambda$ is generic and sufficiently small\upn, then
$$
\RR(M,L)=\Ind_T^G\RR(Y_{\upsilon\upsilon},L_{\upsilon\upsilon}).
$$
\end{theorem}

According to Lemma \ref{lemma:line} the assumption on $L$ is satisfied
for rigid, moment and dual moment bundles.

\begin{proof}
Let $\tau\pre\upsilon$.  We assert that if $\lambda$ is sufficiently
small, then for all $\sigma\pre\tau$ the $G_\tau$-equivariant
orbibundle $L_{\sigma\tau}$ on the $G_\tau$-orbifold $Y_{\sigma\tau}$
is rigid with respect to the subgroup $[G_\tau,G_\tau]$.  Note that
$Y_{\sigma\tau}$ can be regarded as the symplectic cut of
$Y_{\tau\tau}$ with respect to $\ca P_{\sigma\tau}$, so it is enough
to show this for $\sigma=\tau$.  The symplectic cross-section theorem
enables us to reduce this case to the case $\tau=\lie a^*$, where we
need to show that for $\lambda$ small $L$ induces a $[G,G]$-rigid
orbibundle on $Y_{\tau\tau}=\Phi_+\inv(\lambda+\lie a^*)/T\cap[G,G]$.
Now notice that the condition imposed on $L$ implies that $L$ is
almost equivariantly locally trivial on $\Phi\inv(\lie a^*)$ with
respect to the action of $[G,G]$.  Therefore, by Lemma
\ref{lemma:aelt}, $L$ is $[G,G]$-almost equivariantly locally trivial
on $\Phi\inv(G\lambda+\lie a^*)$ for sufficiently small $\lambda$.
The $[G,G]$-rigidity of $L_{\tau\tau}$ is now proved in the same way
as in \ref{part:cutrigid} of Lemma \ref{lemma:cutweight}.

It follows from the rigidity of $L_{\sigma\tau}$ and part
\ref{part:almostcomplex} of Theorem \ref{theorem:rigid} that
$\RR(Y_{\sigma\tau},L_{\sigma\tau})$ is constant as a
$[G_\tau,G_\tau]$-character and can therefore be regarded as a
character of $A_\tau$ or $T$.  Since $Y_{\sigma\tau}$ is equal to the
symplectic cut of $(Y_{\tau\tau})_{\ca P_{\sigma\tau}}$, over the
points in $\ca P_{\sigma\tau}\cap\Delta$ it has the same
$G_\tau$-symplectic quotients as $Y_{\tau\tau}$.  A fortiori, it has
the same $[G_\tau,G_\tau]$-symplectic quotients and therefore by the
equivariant version of \ref{part:hamiltonian} of Theorem
\ref{theorem:rigid} (where we take $G=[G_\tau,G_\tau]$ and $H=A_\tau$)
we have $\RR(Y_{\sigma\tau},L_{\sigma\tau})=
\RR(Y_{\tau\tau},L_{\tau\tau})$ as characters of $A_\tau$, and hence
as characters of $T$, for all $\sigma\pre\tau$ such that $\ca
P_{\sigma\tau}\cap\Delta$ is nonempty.  From
\eqref{equation:subdivision} we now conclude
$$
\RR(M,L)=
\sum_{\tau\pre\upsilon}\sum_{\substack{\sigma\pre\tau\\ \ca
P_{\sigma\tau}\cap\Delta\ne\emptyset}} 
(-1)^{\dim\tau-\dim\sigma}
\Ind_T^G\RR(Y_{\tau\tau},L_{\tau\tau}).
$$
The result now follows from the combinatorial identities
$$
\sum_{\substack{\sigma\pre\tau\\ \ca
P_{\sigma\tau}\cap\Delta\ne\emptyset}}
(-1)^{\dim\tau-\dim\sigma}=
\begin{cases}
1 & \mbox{if $\tau=\upsilon$},\\
0& \mbox{otherwise},
\end{cases}
$$
which derive from the fact that for any simplicial cone $\ca C$ the
sum $\sum_{\ca F\pre\ca C}(-1)^{\dim\ca F}$ is equal to $1$ if $\ca C$
is a point and $0$ otherwise.
\end{proof}

\begin{proof}[Proof of Theorem \ref{theorem:moment}]
Choose a generic $\lambda$ in the principal wall $\upsilon$.  We can
choose $\lambda$ so small that whenever $0<t\le1$ the subdivision
$\lie P_{t\lambda}$ is $G$-admissible and Theorem
\ref{theorem:abelianization} holds with $\lambda$ replaced by
$t\lambda$.  Let us denote by $Y^t$ the global cross-section of the
principal cut $M$ with respect to $\lie P_{t\lambda}$ and by $L^t$ the
corresponding orbibundle.  In view of Theorem
\ref{theorem:abelianization} and the fact that
$\chi_\mu=\Ind_T^G\zeta_\mu$ we need only show that
\begin{equation}\label{equation:abelianization}
\RR(Y^t,L^t)=
\sum_{\mu\in\Lambda^*\cap\Delta}\RR\bigl(M_\mu,L_\mu\shift\bigr)
\,\zeta_\mu.
\end{equation}
Now $L^t$ is not a moment bundle on $Y^t$, so the abelian result
proved in Section \ref{subsubsection:multiplicities} does not directly
apply.  Notice however that $\eu L=\bigcup_{0<t\le1}L^t$ is an
\emph{asymptotic\/} moment bundle on $\eu Y=\bigcup_{0<t\le1}Y^t$ as
defined in Section \ref{subsubsection:asymptotic}.  Using the notation
of that section we obtain from Theorem \ref{theorem:asymptotic} that
for $0<t\le1$
$$
\RR(Y^t,L^t)= \sum_{\mu\in\Lambda^*\cap\Delta}
\RR\bigl(Y^t_{\gamma_\mu(t)},(L^t)_{\gamma_\mu(t)}\shift\bigr)
\,\zeta_\mu.
$$
In addition, it follows from the cross-section theorem that the
quotients $Y^t_{\gamma_\mu(t)}$ and $M_{\gamma_\mu(t)}$ are isomorphic
and also $(L^t)_{\gamma_\mu(t)}\shift\cong L_{\gamma_\mu(t)}\shift$,
so
$$
\RR\bigl(Y^t_{\gamma_\mu(t)},(L^t)_{\gamma_\mu(t)}\shift\bigr)=
\RR\bigl(M_{\gamma_\mu(t)},L_{\gamma_\mu(t)}\shift\bigr)=
\RR\bigl(M_\mu,L_\mu\shift\bigr)
$$
by Theorem \ref{theorem:independent}.  This proves
\eqref{equation:abelianization}.
\end{proof}

\begin{proof}[Proof of Theorem \ref{theorem:dual}]
The proof of the multiplicity formula
\eqref{equation:dualmultiplicity} is completely analogous to the proof
of Theorem \ref{theorem:moment}.  The formula implies that the support
of the multiplicity function is contained in the orbit of
$-\inter\Delta$ under the affine action $w\odot\mu=w(\mu+\rho)-\rho$
of $\eu W$,
\begin{equation}\label{equation:necklace}
\supp N_{L\inv}\subset -\rho+\bigcup_{w\in\lie
W}w(\rho-\inter\Delta).
\end{equation}
Since $\inter\Delta$ is entirely contained in the principal open wall
$\sigma$ of $M$, it follows from Lemma \ref{lemma:reflect} below that
the intersection of the right-hand side of \eqref{equation:necklace}
with $\Lambda^*_+$ is contained in
$*\bigl(\inter\Delta-2(\rho-\rho_\sigma)\bigr)$.
\end{proof}

We thank Dan Barbasch for helping us prove the following lemma.

\begin{lemma}\label{lemma:reflect}
Let $\lambda$ be a dominant weight and let $\sigma$ be the open wall
of $\tplus$ containing $\lambda$.  Let $w_\sigma$ be the longest
element in the Weyl group $\eu W_\sigma\subset\eu W$ of the
centralizer $G_\sigma$.  Then the following conditions are equivalent.
\begin{enumerate}
\item\label{part:reflect}
There exists $w\in\eu W$ such that $w\odot(-\lambda)$ is
dominant\upn;
\item\label{part:regular}
$\lambda-\rho$ is regular\upn;
\item\label{part:dominantregular}
$w_\sigma(\lambda-\rho)$ is dominant regular\upn;
\item\label{part:dominant}
$\lambda-2(\rho-\rho_\sigma)$ is dominant.
\end{enumerate}
If \ref{part:reflect} holds\upn, then $w=w_0w_\sigma$ and
$w\odot(-\lambda)= \bigl(\lambda-2(\rho-\rho_\sigma)\bigr)^*=
\lambda^*-2(\rho-\rho_\sigma^*)$.
\end{lemma}

\begin{proof}
Let $R$ be the root system of $G$ and $R_\sigma=\{\alpha\in
R:(\lambda,\check\alpha)=0\}$ the root system of $G_\sigma$.  Let
$R^+$ and $R_\sigma^+$ denote the corresponding sets of positive
roots.  Note that $w_\sigma$ fixes $\sigma$, so
$w_\sigma\lambda=\lambda$.  Furthermore, $w_\sigma$ permutes the
elements of $R^+-R_\sigma^+$ and sends $R_\sigma^+$ to $R_\sigma^-$.
This implies
\begin{equation}\label{equation:reflect}
w_\sigma\rho=\rho-2\rho_\sigma\qquad\text{and}\qquad
w_\sigma(\lambda-\rho)=\lambda+2\rho_\sigma-\rho.
\end{equation}

If $w\odot(-\lambda)$ is dominant for some $w\in\eu W$, then
$w(-\lambda+\rho)$ is dominant regular, so $\lambda-\rho$ is regular.
This shows that \ref{part:reflect} implies \ref{part:regular}.  The
implications \ref{part:dominantregular} $\Rightarrow$
\ref{part:dominant} $\Rightarrow$ \ref{part:reflect} are obvious from
\eqref{equation:reflect}.  

Next we show that \ref{part:regular} implies
\ref{part:dominantregular}.  It suffices to show that
$(\lambda-\rho,\check\alpha)\ge0$ for all $w_\sigma$-positive roots
$\alpha$.  If $\alpha=-\beta$ with $\beta\in R^+_\sigma$, then
$(\lambda-\rho,\check\alpha)=(\rho,\check\beta)>0$. If $\alpha\in
R^+-R_\sigma^+$, then can write
$\alpha=\beta_1+\beta_2+\dots+\beta_k$, where $\beta_1$,
$\beta_2,\dots$, $\beta_k$ are simple, $\beta_1$ is not in
$R^+_\sigma$, and every partial sum
$\alpha_i=\beta_1+\beta_2+\dots+\beta_i$ is in $R^+$.  Note that
$(\lambda,\check\alpha_1)$ is a positive integer and
$(\lambda-\rho,\check\alpha_1)=(\lambda,\check\alpha_1)-1$, so
$(\lambda-\rho,\check\alpha_1)\ge0$.  Note also that
$(\lambda-\rho,\check\alpha_{i+1}-\check\alpha_i)$ is either positive
or equal to $-1$ for every $i$.  This implies that if
$(\lambda-\rho,\check\alpha)=(\lambda-\rho,\check\alpha_k)$ was
negative, then $(\lambda-\rho,\check\alpha_i)$ would be equal to $0$
for some $i$, which contradicts the regularity of $\lambda-\rho$.  We
conclude that $(\lambda-\rho,\check\alpha)\ge0$.

Finally, if \ref{part:reflect} holds, then $-w(\lambda-\rho)$ is
dominant regular, and so is $w_\sigma(\lambda-\rho)$ by
\ref{part:dominantregular}.  It follows that $w=w_0w_\sigma$ and hence
$w\odot(-\lambda)= \lambda^*-2\rho+2\rho_\sigma^*$ by
\eqref{equation:reflect}.
\end{proof}

\appendix

\section{Normal forms}\label{section:normalform}

Section \ref{subsection:minimal} contains a brief review of minimal
coupling and some observations on deformation equivalence and
equivariant blowing up.  In Section \ref{subsection:darboux} we prove
a relative version of the constant-rank embedding theorem and a number
of related embedding and deformation results.

\subsection{Minimal coupling}\label{subsection:minimal}

Let $(B,\omega_B)$ be a symplectic orbifold and
$P\stackrel{\pi}{\longleftarrow}B$ a principal $H$-orbibundle, where
$H$ is a compact Lie group.  Let $\theta\in \Omega^1(P,\lie h)$ be a
principal connection on $P$ and $(Q,\omega_Q)$ a Hamiltonian
$H$-orbifold with moment map $\Phi_Q\colon Q\to\lie h^*$.  Let
$\pr_{P,Q}$ denote the projection from $P\times Q$ onto $P$, resp.\
$Q$, and let $\langle{\cdot},{\cdot}\rangle\colon\lie h^*\times\lie
h\to\R$ be the dual pairing.  Then the two-form
$$
\pr_P^*\pi^* \omega_B + \pr_Q^*\omega_Q + \d\langle
\pr_Q^*\Phi_Q,\pr_P^*\theta\rangle
$$
on the principal orbibundle
$$ H\longrightarrow P\times Q\longrightarrow P\times^HQ $$
is basic, and therefore descends to a closed two-form on the
associated orbibundle $P\times^H Q$, which is called the \emph{minimal
coupling form}.  The following result is due to Sternberg.  A proof
for the manifold case can be found in \cite{gu:sy}.  The
generalization to orbifolds is straightforward.

\begin{theorem}\label{theorem:minimalcoupling} 
\begin{enumerate} 
\item
The minimal coupling form is nondegenerate in a neighbourhood of
$P\times^H \Phi_Q\inv(0)$.
\item\label{part:base}
A Hamiltonian $G$-action on $B$ with moment map $\Psi_B\colon B\to\lie
g^*$ that lifts to an action on $P$ by $\theta$-preserving principal
orbibundle automorphisms induces a Hamiltonian $G$-action on
$P\times^HQ$.  The moment map on $P\times^HQ$ is a sum
$\Phi_B+\Phi_\theta$\upn, where $\Phi_B$ and $\Phi_\theta$ are defined
as follows\upn: $\langle\Phi_B,\xi\rangle$ and
$\langle\Phi_\theta,\xi\rangle$ are the functions induced by the
$H$-invariant functions $\pr_P^*\pi^*\langle\Psi_B,\xi\rangle$\upn,
resp.\ $\langle \pr_Q^*\Phi_Q,\pr_P^*\iota(\xi_P)\theta\rangle$.
\item\label{part:fibre}
A Hamiltonian action of a Lie group $G'$ on $Q$ with moment map
$\Psi_Q\colon Q\to (\g')^*$ that commutes with the action of $H$
induces an $G'$-action on $P\times^HQ$\upn, which is Hamiltonian with
moment map induced by the $H$-invariant map $\pr_Q^*\Psi_Q$.
\end{enumerate} 
\end{theorem} 

Here an \emph{automorphism\/} of a fibre orbibundle refers to a
diffeomorphism of the total space that preserves the structure group
and maps fibres to fibres and hence induces a diffeomorphism from the
base onto itself.  (In the case of the principal orbibundle $P$ this
simply means an $H$-equivariant diffeomorphism of the total space.)

Weinstein observed that the associated orbibundle can also be obtained
as a symplectic quotient.  The \emph{universal phase space\/} of the
principal orbibundle $P$ is the orbifold $ P\times\lie h^*$.  It
carries a closed two-form
$\omega_\theta=\pr_P^*\pi^*\omega_B+\d\langle\pr_{\lie
h^*},\pr_P^*\theta\rangle$, which is nondegenerate in a neighbourhood
of $P\times\{0\}$.  The $H$-action on the universal phase space is
Hamiltonian with moment map given by $\pr_{\lie h^*}$.  The $H$-action
on $P\times\lie h^*\times Q$ is therefore Hamiltonian with moment map
given by
$$\Psi(p,\beta,q)=\beta+\Phi_Q(q).$$
Since the $H$-action on $P$ is locally free, the symplectic quotient
$(P\times\lie h^*\times Q)\qu H$ is a symplectic orbifold.  The map
$(p,q)\mapsto(p,-\Phi_Q(q),q)$ is an $H$-equivariant diffeomorphism
onto $\Psi\inv(0)$, and therefore descends to a diffeomorphism
$$
P\times^HQ\longrightarrow(P\times\lie h^*\times Q)\qu H,
$$
which one can easily show to be symplectic (with respect to the
minimal coupling form defined by the connection on $ P$).  Sometimes
the form $\omega_\theta$ on the universal phase space is globally
nondegenerate.

\begin{example}\label{example:model}
Suppose $H$ is a subgroup of $G$ acting on $G$ by right
multiplication.  Let $B$ be the symplectic manifold $T^*(G/H)$ and
define $P$ to be the pullback of the bundle $B\to G/H$ under the
projection $G\to G/H$.  Let us identify $T^*G$ with $G\times\g^*$ by
means of left-invariant one-forms.  Then $B\cong G\times^H\h^0$, where
$\h^0$ is the annihilator of $\h$ in $\g^*$, so $P\cong G\times\h^0$
is a principal $H$-bundle over $B$.  Choose a $G$-equivariant
connection $\theta$ on $P$, that is an $H$-equivariant splitting $\lie
g=\lie h\oplus\lie m$.  Then $\lie h^0\cong\lie m^*$ and we can
identify $\lie h^*$ with a subspace of $\lie g^*$.  We thus obtain a
diffeomorphism from $P\times\lie h^*\cong G\times\lie h^0\times\lie
h^*$ to $T^*G\cong G\times\lie g^*$, and it is easy to see that the
pullback of the standard symplectic form on $T^*G$ is equal to
$\omega_\theta$.  Thus $\omega_\theta$ is symplectic globally.  We
conclude that the associated bundle
\begin{equation}\label{equation:bundle}
F(H,Q)=G\times^H(\h^0\times Q)\cong P\times^HQ\cong (T^*G\times Q)\qu
H
\end{equation}
is a symplectic orbifold.  (In this example $F(H,Q)$ is not merely an
orbibundle but a genuine locally trivial fibre bundle, whose fibres
happen to be orbifolds.)

Now note that the action of $G$ on itself by left multiplication lifts
to a Hamiltonian $G$-action on $B=T^*(G/H)$ and to an action on $P$ by
orbibundle automorphisms that preserve the connection.  By
\ref{part:base} of Theorem \ref{theorem:minimalcoupling} we obtain a
Hamiltonian $G$-action on $F(H,Q)$ with moment map given by
\begin{equation}\label{equation:map}
[g,\beta,q]\longmapsto g(\beta+\Phi_Q(q)),
\end{equation}
where $[g,\beta,q]$ denotes the $H$-orbit through $(g,\beta,q)\in
G\times\h^0\times Q$.  The zero level set is therefore the bundle
$G\times^H \Phi_Q\inv(0)$, and $F(H,Q)\qu G\cong Q\qu H$ (reduction in
stages).

Finally, let $K$ be a compact Lie subgroup of $\Diff(Q)$ containing
the image of $H$ under the action map $\rho\colon H\to\Diff(Q)$.
Suppose that $K$ acts on $Q$ in a Hamiltonian fashion and that the
$H$-moment map $\Phi_Q$ is equal to the $K$-moment map followed by the
natural map $\rho^*\colon\lie k^*\to\lie h^*$.  Let $N_{G\times K}(H)$
be the normalizer of $H$ under the embedding $H\to G\times K$ given by
$h\mapsto(h,\rho(h))$.  Then the quotient
\begin{equation}\label{equation:group}
K(H,Q)=N_{G\times K}(H)/H
\end{equation}
acts on $F(H,Q)$ in a Hamiltonian fashion, and the action commutes
with the action of $G$.  This group has a particularly simple
interpretation.
\end{example}

\begin{lemma}\label{lemma:automorphism}
The group $K(H,Q)$ is canonically isomorphic to
$\Aut\bigl(F(H,Q)\bigr)^G$\upn, the group of $G$-equivariant
automorphisms of the fibre bundle $F(H,Q)$ that preserve the structure
group $K$.
\qed
\end{lemma}

\begin{example}\label{example:blowdeform}
As in Section \ref{subsubsection:blow} let $S$ be a locally closed
$G$-invariant symplectic suborbifold of the Hamiltonian $G$-orbifold
$M$ with normal bundle $N$ and let $\K$ be a $G$-invariant compact
subset of $M$ such that $\K\cap S$ is closed.  Choose a $G$-invariant
complex structure on $N$; then $N\cong P\times^{\U(k)}\C^k$, where $P$
is the unitary frame orbibundle of $N$.  According to \ref{part:base}
of Theorem \ref{theorem:minimalcoupling}, a neighbourhood of the zero
section in $N$ is a Hamiltonian $G$-orbifold.  Its moment map
$\Phi_N=\Phi_S+\Phi_\theta$ has the property that $\Phi_N\inv(0)\cap
S=\Phi_S\inv(0)$.

As noted in Section \ref{subsubsection:blow} there exist $\delta>0$
and a $G$-invariant open neighbourhood $U$ of $\K$ such that over
$U\cap S$ the minimal coupling form is nondegenerate on the disc
bundle $N(\delta)$ and $N(\delta)|_{U\cap S}$ embeds properly,
equivariantly and symplectically into $U$.  For all $\eps<\delta$ the
blowup $\Bl(U,S,\omega,j,\theta,\iota,\eps)$ of $U$ along $S$ is
well-defined.  Let us now assume that the projection $\pr_S$ preserves
the level set $\Phi_N\inv(0)$, in other words,
$\pr_S\Phi_N\inv(0)=\Phi_S\inv(0)$, or equivalently, $\Phi_N\inv(0)$
is equal to $\Phi_S\inv(0)\cap\Phi_\theta\inv(0)$.  This implies that
$\Phi_N\inv(0)$ is a conical subset of $N$, because $\Phi_B$ is
constant along the fibres of $\pr_S$ and $\Phi_\theta$ is quadratic.
If $\eps_0<\eps_1$ we can retract the orbifold $U-N(\eps_0)$ smoothly
and equivariantly onto $U-N(\eps_1)$ by pushing points outward along
the fibres of $N$.  This retraction leaves the complement of
$N(\delta)$ in $U$ fixed and preserves the zero level set of $\Phi_N$.
It therefore induces a deformation equivalence between the
$\eps_0$-blowup and the $\eps_1$-blowup.
\end{example}

\begin{lemma}\label{lemma:blowdeform}
If the projection $\pr_S$ preserves the zero fibre of $\Phi_N$\upn,
then for all $\eps_0$ and $\eps_1<\delta$ the blowups
$\Bl(U,S,\omega,j,\theta,\iota,\eps_0)$ and
$\Bl(U,S,\omega,j,\theta,\iota,\eps_1)$ are deformation equivalent as
Hamiltonian $G$-orbifolds.  
\qed
\end{lemma}

\subsection{Embedding theorems}\label{subsection:darboux}

In this section we present an addendum to the constant-rank embedding
theorem and some other embedding and deformation results that rely on
a straightforward extension of Moser's method for proving the Darboux
Theorem.  Moser's method shows that under certain conditions the fact
that two symplectic forms are deformation equivalent implies that they
are strongly isotopic.  It is a trivial observation that the isotopy
obtained by Moser's method ``depends smoothly on parameters''.  This
leads to a relative version of the constant-rank embedding theorem and
also enables us in certain cases to deform a path of diffeomorphisms
to a path of symplectomorphisms.

The details are as follows.  Let $\eu B$ be an orbifold and let
$\pi\colon\eu M\to\eu B$ be a fibre orbibundle over $\eu B$.  A
\emph{relative symplectic form\/} on $\eu M$ or a \emph{symplectic
form on $\eu M$ over\/} $\eu B$ is a two-form $\omega$ on the vertical
tangent orbibundle $T_\vertical\eu M=\ker\d\pi$, the restriction of
which to every fibre of $\pi$ is closed and nondegenerate.  It is
obvious how to define, in the presence of $G$-actions on $\eu B$ and
$\eu M$ such that $\pi$ is equivariant, a \emph{relative moment map\/}
$\Phi\colon\eu M\to\lie g^*$ for the action.

Consider a smooth path $\omega_t$ of relative symplectic forms on $\eu
M$ defined for $0\le t\le1$.  Denote by $\dot{\omega}_t$ the
$t$-derivative of $\omega_t$ and suppose that
$\dot{\omega}_t=\d_\pi\sigma_t$, where
$\sigma\colon[0,1]\to\Gamma\bigl(\Lambda^1T_\vertical\eu M\bigr)$ is a
smooth path of vertical one-forms on $\eu M$, and $\d_\pi$ denotes the
vertical exterior derivative.  Define a time-dependent vertical vector
field $\Xi_t$ on $M$ by
$$ 
\Xi_t=-\omega_t^\flat\sigma_t,
$$
where $\omega_t^\flat\colon(T_\vertical\eu M)^*\to T_\vertical\eu M$
denotes the lowering operator associated to $\omega$, and let $\psi_t$
be its flow.  Then clearly $\psi_t$ preserves the fibres of $\pi$ and
$\ca L(\Xi_t)\omega_t=-\dot{\omega}_t$, so
\begin{equation}\label{equation:moser}
\psi_t^*\omega_t=\omega_0
\end{equation}
for those $t$ and those points of $\eu M$ at which the flow $\psi_t$
is defined.

Sometimes there is a natural choice for the one-forms $\sigma_t$ and
an a priori estimate for the existence interval of the flow $\psi_t$.

\begin{example}\label{example:dmw}
Let $\eu Y$ be a fibre orbibundle over $\eu B$ and $\eu Y\to\eu M$ a
locally closed embedding that commutes with the projections $\eu
Y\to\eu B$ and $\eu M\to\eu B$.  Regard $\eu Y$ as a subset of $\eu M$
and suppose that $\dot{\omega}_t|_{\eu Y}=0$.  Let $\eu O\to\eu B$ be
a relative tubular neighbourhood of $\eu Y$ (i.~e.\ the fibres of
$\eu O$ are tubular neighbourhoods of the fibres of $\eu Y$) and
define a homotopy $I\colon\eu O\times[0,1]\to\eu O$ by radial
retraction along the fibres of the orbibundle projection $\pr_{\eu
Y}\colon\eu O\to\eu Y$.  Then $I_1=\id_{\eu O}$ and $I_0=\pr_{\eu Y}$.
Let $\kappa_I$ be the associated chain homotopy on the de Rham complex
of $\eu O$, which is given by
$\kappa_I\alpha=\int\iota(\partial/\partial t)I^*\alpha\,\d t$ for
$\alpha\in\Omega(\eu O)$.  Put $\sigma_t=\kappa_I\dot{\omega}_t$.
Then
$$
\d_\pi\sigma_t= \d_\pi\kappa_I\dot{\omega}_t=
\d_\pi\kappa_I\dot{\omega}_t+\kappa_I\d_\pi\dot{\omega}_t=
\dot{\omega}_t-\pr_{\eu Y}^*\bigl(\dot{\omega}_t|_{\eu Y}\bigr)=
\dot{\omega}_t.
$$
Furthermore, $\Xi_t=0$ on $\eu Y$, so $\eu Y$ is fixed under the flow
$\psi_t$.  It follows that there exists an open neighbourhood $\eu O'$
of $\eu Y$ contained in $\eu O$ such that the flow $\psi_t$ is defined
for all $t\in[0,1]$ and all initial values in $\eu O'$.
\end{example}

As a first application of the relative Moser method, we show how in
special circumstances an isotopy can be deformed to a symplectic
isotopy.

\begin{example}\label{example:straightline}
In the setting of Example \ref{example:dmw}, let $\eu B=[0,1]$, $\eu
M=M\times\eu B$, and $\eu Y=Y\times\eu B$.  Here $M$ is an orbifold
and $Y$ a locally closed suborbifold of $M$.  Let $\omega$ be a fixed
symplectic form on $M$ and let $F\colon O\times\eu B\to O$ be an
isotopy of a tubular neighbourhood $O$ of $Y$ leaving $Y$ pointwise
fixed.  Assume that $F$ is symplectic at all points of $Y$ in the
sense that $(F_b^*\omega)_y=\omega_y$ for all $b\in\eu B$ and all $y$
in $Y$.  (We do not assume that $F$ starts at the identity.)  Put $\eu
O=O\times\eu B$.  Define a path of vertical two-forms $\omega_t$ on
$\eu O$ by putting $\omega_t=(1-t)\omega+tF_b^*\omega$ on
$O\times\{b\}\subset\eu O$.  Then $(\omega_t)_y=\omega_y$ for all
$y\in\eu Y$, so $\omega_t$ is symplectic on a neighbourhood of $\eu Y$
in $\eu O$.  Furthermore, on $Y\times\{b\}$ we have
$\dot{\omega}_t=(-\omega+F_b^*\omega)=0$, so $\dot{\omega}_t=0$ on
$\eu Y$.  We are therefore in the situation of Example
\ref{example:dmw} and obtain a flow $\psi_t\colon\eu O'\to\eu O$ which
is defined for all $t$ and for all initial values in a small open $\eu
O'$ containing $\eu Y$ and satisfies \eqref{equation:moser}.  We may
assume $\eu O'$ is of the form $O'\times\eu B$ for some open subset
$O'$ of $O$ containing $Y$.  For $b\in\eu B$ let $\rho_b\colon O'\to
O'\times\{b\}$ be the diffeomorphism $\rho_b(m)=(m,b)$.  For each
$t\in[0,1]$ define an isotopy $F^{(t)}\colon O'\times\eu B\to O$ by
$F^{(t)}_b=F_b\rho_b\inv\psi_t\rho_b$.  Then $F^{(0)}=F$ and
$\bigl(F^{(1)}_b\bigr)^*\omega=
\bigl(\rho_b\inv\bigr)^*\psi_1^*\rho_b^*F_b^*\omega=\omega$, so
$F^{(1)}$ is symplectic.  This means that we have constructed a path
of isotopies joining $F=F^{(0)}$ to the symplectic isotopy $F^{(1)}$.

We note two additional properties of $\psi_t$.  Firstly, if $F_b$ is a
symplectomorphism for some $b\in\eu B$, then on $O'\times\{b\}$ we
have $\sigma_t=\kappa_I\dot{\omega}_t=0$ for all $t$ and so $\Xi_t=0$.
Therefore the flow $\psi_t$ is trivial on $O'\times\{b\}$.  Secondly,
suppose that $\omega$ is invariant under the action of a compact Lie
group $K$ which leaves $Y$ invariant and preserves the projection
$O\to Y$.  Then if $F_b$ is $K$-equivariant for some $b\in\eu B$, the
flow $\psi_t$ is equivariant on $O'\times\{b\}$.
\end{example}

In this context we also have the following elementary result.

\begin{lemma}\label{lemma:hamiltonian}
Every symplectic isotopy $O'\times[0,1]\to O$ that leaves $Y$ fixed
and starts at the identity is Hamiltonian.
\end{lemma} 

\begin{proof}
Let $F$ be such an isotopy.  Consider its infinitesimal generator
$\eta_b$ (where $b\in[0,1]$) and note that the one-form
$\beta=\iota(\eta_b)\omega$ is closed because $F$ is symplectic, and
that $\beta|_Y=0$ because $Y$ is fixed under $H$.  The function
$\kappa_I\beta$ satisfies $\d\kappa_I\beta= \d\kappa_I\beta+
\kappa_I\d\beta= \beta-\pr_Y^*(\beta|_Y)= \beta= \iota(\eta_b)\omega$.
In other words, the vector field $\eta_b$ is generated by the
time-dependent Hamiltonian $\kappa_I\beta$.
\end{proof}

The following lemma is used in Section \ref{subsubsection:blow}.  We
use the notation of that section and of Example
\ref{example:blowdeform}.

\begin{lemma}\label{lemma:circle}
For $\delta'<\delta$ let $f\colon N(\delta')\to N(\delta)$ be a
symplectic map restricting to the identity on $S$.  Then there exist
$\delta''<\delta'$ and a symplectic isotopy $H\colon
N(\delta'')\times[0,1]\to N(\delta)$ such that $H_0=f$\upn, $H_1$ is
$S^1$-equivariant\upn, $H_b|_S=\id_S$ for all $b\in[0,1]$.  If $f$ is
$G$-equivariant\upn, then $H$ can be chosen to be equivariant.
\end{lemma}

\begin{proof}
Below we construct an isotopy $F\colon N(\delta')\times[0,1]\to U$
fixing $S$ such that $F_0=f$, $F_1$ is $S^1$-equivariant and
$F_b^*\omega_m=\omega_m$ for all $m\in S$ and $0\le b\le1$.  Then we
put $H=F^{(1)}$ as in Example \ref{example:straightline} above.  As we
have seen, $H$ is a symplectic isotopy, and since $F_0=f$ preserves
the symplectic form, the flow $\psi_t$ is trivial on $
N(\delta')\times\{b\}$ for all $t$, so $H_0=F_0\circ\id=f$.
Furthermore, since $F_1$ is equivariant, so is $H_1$.

The construction of the isotopy $F$ is in two stages.  Between time
$b=0$ and $1/2$ we isotop $f$ to its fibre derivative $T_Nf$ by means
of the obvious isotopy $F(m,b)=(1-2b)\inv f\bigl((1-2b)m\bigr)$.  The
fact that $f$ is a symplectic map and leaves $S$ fixed implies that
for all $m\in S$ the derivative $T_mf$ preserves the direct sum
decomposition of (uniformized) tangent spaces $\tilde{T}_mN=
\tilde{T}_mS\oplus\tilde{T}_0N_m$.  It follows from this that $f$ has
the same derivative as $T_Nf$ at all points in the zero section $S$,
and in fact $T_mF_b=T_mf$ for all $b\in[0,1/2]$ and all $m$ in $S$.
Consequently, $(F_b^*\omega)_m= (f^*\omega)_m=\omega_m$ for all $b$
and $m$.  

The second half of the isotopy comes about as follows.  Note that
$T_Nf$ is an element of $\Aut_S(N)$, that is the group of linear
automorphisms of $N$ that preserve the symplectic forms on the fibres
and restrict to the identity on $S$.  If $P$ is the Hermitian frame
orbibundle of $N$, then $\Aut_S(N)$ can be viewed as the space of
sections of the associated orbibundle $P\times^{\U(n)}\Sp(2n,\R)$.
Using the retraction of $\Sp(2n,\R)$ onto its maximal compact subgroup
$\U(n)$ we can construct a path in $\Aut_S(N)$ defined for $1/2\le
b\le1$ starting at $T_Nf$ and ending at a Hermitian automorphism of
$N$.  Observe that symplectic orbibundle automorphisms preserve the
symplectic form at all points of the zero section and that Hermitian
automorphisms commute with the scalar $S^1$-action.  By composing the
two isotopies we obtain the requisite isotopy $F$.

It is not hard to check that each step in this proof can be made
equivariant with respect to the action of $G$.  It follows that the
isotopy $H$ can be made equivariant.
\end{proof}

Another application of the relative Moser method is the relative
Darboux-Moser-Weinstein Theorem: if $\omega_0$ and $\omega_1$ are
relative symplectic forms on a fibre orbibundle $\pi\colon\eu M\to\eu
B$ such that $\omega_{0,y}=\omega_{1,y}$ for all $y$ in a locally
closed suborbibundle $\eu Y$ of $\eu M$, then there exist open
neighbourhoods $\eu U_0$ and $\eu U_1$ of $\eu Y$ and a diffeomorphism
$f\colon\eu U_0\to\eu U_1$ commuting with $\pi$ such that $f(y)=y$,
$\d_\pi f_y=\id_{T_y\eu M}$ for all $y\in\eu Y$, and
$f^*\omega_1=\omega_0$.  The proof is word for word the same as in the
absolute case, relying on linear interpolation between $\omega_0$ and
$\omega_1$.  In turn this leads to relative versions of all the usual
embedding theorems in symplectic geometry.

As an example we state the relative constant-rank embedding theorem.
Let $\eu Z$ be a fibre orbibundle over $\eu B$ and let $\tau$ be a
vertical two-form on $\eu Z$ that is closed on every fibre.  Assume
that $\tau$ has constant rank on $\eu Z$.  Assume further that $G$
acts on $\eu B$, that $\eu Z$ is an equivariant orbibundle, and that
the action on $\eu Z$ is Hamiltonian in the sense that there exists a
$G$-equivariant map $\Phi_{\eu Z}\colon\eu Z\to\lie g^*$ satisfying
$\langle\d_\pi\Phi_{\eu Z},\xi\rangle=\iota(\xi_{\eu Z})\tau$.  Let
$\eu N$ be a $G$-equivariant symplectic vector orbibundle over $\eu Z$
with fibre symplectic form $\sigma$.  Now let $\omega$ be a relative
symplectic form on $\eu M$ and assume $G$ acts on $\eu M$ in a
Hamiltonian fashion with relative moment map $\Phi$.  An
\emph{embedding of $\eu Z$ into $\eu M$ with normal bundle\/} $\eu N$
is an embedding of fibre orbibundles $\iota\colon\eu Z\to\eu M$ such
that $\iota^*\omega=\tau$, $\iota^*\Phi=\Phi_{\eu Z}$, and the
pullback under $\iota$ of the relative symplectic normal bundle of
$\iota(\eu Z)$ in $(\eu M,\omega)$ is isomorphic to $(\eu N,\sigma)$.

The \emph{standard\/} embedding $\eu Z\hookrightarrow\eu Y$ with
normal bundle $\eu N$ is constructed as follows.  As an orbifold, $\eu
Y$ is the total space of the direct sum $\eu S\oplus\eu N$, where $\eu
S$ is the orbibundle on $\eu Z$ dual to the suborbibundle $\ker\tau$
of the vertical tangent bundle $T_\vertical\eu Z$.  The relative
symplectic form and moment map $(\omega_{\eu Y},\Phi_{\eu Y})$ on $\eu
Y$ are constructed in two stages.

At the first stage one chooses a section $s$ of the orbibundle
map $T_\vertical^*\eu Z\to\eu S$ and defines a closed two-form
$\omega_{\eu S}$ on the fibres of the projection $\eu S\to\eu B$ by
$\omega_{\eu S}=\pr_Z\tau+s^*\Omega$, where $\Omega$ is the
standard symplectic form on the fibres of $T_\vertical^*\eu Z\to\eu
B$.  Near $\eu Z$ the form $\omega_{\eu S}$ is nondegenerate in the
vertical direction, and the $G$-action on $\eu S$ is Hamiltonian with
moment map given by $\Phi_{\eu S}=\pr_{\eu Z}^*\Phi_{\eu
Z}+s^*\Phi_\vertical$, where
$\langle\Phi_\vertical(p),\xi\rangle=p(\xi_{\eu Z})$, the standard
moment map on $T_\vertical^*\eu Z$.

At the second stage one notices that as an orbifold $\eu Y=\eu
S\oplus\eu N$ is identical to the total space of the pullback of $\eu
N$ along the map $\eu S\to\eu Z$, and by means of (fibrewise) minimal
coupling constructs a relative symplectic form $\omega_{\eu Y}$ on
$\eu Y$, using the relative symplectic form $\omega_{\eu S}$ on the
base $\eu S$, an invariant $\sigma$-compatible almost complex
structure $\eu J$ on $\eu N$, and a (relative) connection one-form
$\theta$ on the orbibundle $\eu P$ of $\eu J$-unitary frames on $\eu
N$.  By Theorem \ref{theorem:minimalcoupling} the $G$-action on $\eu
Y$ is Hamiltonian with respect to $\omega_{\eu Y}$ with relative
moment map $\Phi_{\eu Y}$, and it is straightforward to check that the
zero section $\eu Z\hookrightarrow\eu Y$ is an embedding of $\eu Z$
with normal bundle $\eu N$.  The relative version of the constant-rank
embedding theorem is now proved in the same way as the absolute
version; cf.\ e.~g.\ \cite{sj:st}.

\begin{theorem}[relative constant-rank
embeddings]\label{theorem:constantrank}
For every embedding $\iota$ of $\eu Z$ into $\eu M$ with normal bundle
$\eu N$ there exist a $G$-invariant open neighbourhood $\eu U$ of $\eu
Z$ in $\eu Y$ and an isomorphism of relative Hamiltonian $G$-orbifolds
$$f\colon(\eu U,\omega_{\eu Y},\Phi_{\eu Y})\longrightarrow (\eu
M,\omega,\Phi)$$
onto an open neighbourhood of $\iota(\eu Z)$ in $\eu M$ such that the
diagram
$$
\xymatrix{{\eu Z}\ar[d]\ar[dr]^{\iota}\\{\eu U}\ar[r]_f & {\eu M}}
$$
commutes.  
\qed
\end{theorem}

The zero section $\eu Z\hookrightarrow\eu S$ is a coisotropic
embedding of $\eu Z$ and the zero section $\eu S\hookrightarrow\eu Y$
is a symplectic embedding of $\eu S$.  It is not hard to see that $\eu
S$ is a \emph{minimal\/} symplectic suborbifold of $\eu Y$ containing
$\eu Z$ in the sense that if $\eu S'$ is a locally closed symplectic
suborbifold of $\eu Y$ such that $\eu Z\subset\eu S'\subset\eu S$,
then $\eu S'$ is open in $\eu S$.  Theorem \ref{theorem:constantrank}
thus proves the existence of minimal symplectic suborbifolds
containing a given constant-rank suborbifold.  To what extent are
minimal symplectic suborbifolds unique?  We shall answer this question
in the absolute case only, though even there the proof uses the
relative constant-rank embedding theorem.

\begin{theorem}\label{theorem:centre}
Let $(M,\omega,\Phi)$ be a Hamiltonian $G$-orbifold and let $Z$ be a
$G$-invariant compact suborbifold of constant rank.  Let $S_0$ and
$S_1$ be minimal $G$-invariant locally closed symplectic suborbifolds
of $M$ containing $Z$.  Then there exist an automorphism $f\colon M\to
M$ of the Hamiltonian $G$-orbifold $M$ and $G$-invariant open
neighbourhoods $U_0$ and $U_1$ of $Z$ such that $f$ fixes $Z$ and maps
$U_0\cap S_0$ onto $U_1\cap S_1$.
\end{theorem}

\begin{proof}
Let $N$ be the symplectic normal bundle of $Z$ in $M$, let
$Z\hookrightarrow Y$ be the standard embedding of $Z$ with normal
bundle $N$, and choose an embedding $U\hookrightarrow M$ as in Theorem
\ref{theorem:constantrank}.  Below we find an equivariant symplectic
isotopy $F\colon U'\times[0,1]\to U$ of an invariant open $U'$ such
that $Z\subset U'\subset U$, $F_0$ is the identity on $U'$, and $F$
leaves $Z$ fixed.  According to Lemma \ref{lemma:hamiltonian}, $F$ is
generated by a time-dependent Hamiltonian vector field.  We extend
this vector field to $M$ by multiplying its Hamiltonian function by a
smooth cutoff function that is supported on $U'$ and identically equal
to $1$ on a smaller $U''\subset U'$.  Since $Z$ is compact, the
resulting Hamiltonian vector field is compactly supported and hence
integrates to a globally defined isotopy $\check F$ of $M$; and
$f=\check{F}_1$ is the desired automorphism.

To construct the isotopy $F$, we may without loss of generality
replace $M$ with the model space $Y$ and assume $S_0$ to be the
suborbibundle $S=(\ker\tau)^*$ of $Y$.  This means that we can
identify $S_0$ with the orbibundle $TS_0|_Z$ over $Z$.  The
construction is in two steps.  First we find a symplectic isotopy of
$U$ that fixes $Z$ and moves $S_0=TS_0|_Z$ to $TS_1|_Z$, and then we
construct $F$ in the special case where $TS_0|_Z=TS_1|_Z$.

\emph{Step}~1.  Regard $Y$ as an orbibundle over $Z$ and note that
both $TS_0|_Z$ and $TS_1|_Z$ are suborbibundles of $Y$ that are
complementary to $N$.  We can therefore select a path of
suborbibundles $N^\perp_t$ of $Y$ defined for $0\le t\le1$ such that
$N^\perp_0=TS_0|_Z$, $N^\perp_1=TS_1|_Z$, and $N^\perp_t$ is
complementary to $N$ for all $t$.  We assert that, for all $t$, near
the zero section the total space of $N^\perp_t$ is a symplectic
suborbifold of $Y$ and that $Z$ is coisotropic in $N^\perp_t$.  This
is proved by showing that $TN^\perp_t|_Z$ is a symplectic
suborbibundle of $TY|_Z$, as follows.  By construction, $TY|_Z$ is
canonically a symplectic direct sum $R\oplus K\oplus K^*\oplus N$,
where $K=\ker\tau$, $R=TZ/K$, and the orbibundle $K\oplus K^*$ carries
the canonical symplectic form on its fibres.  The fact that
$TN^\perp_t|_Z$ is symplectic now follows from the first assertion of
Lemma \ref{lemma:complement}.  We then apply the relative
\emph{coisotropic\/} embedding theorem, that is to say, we apply
Theorem \ref{theorem:constantrank} with $\eu B=[0,1]$, $\eu
Z=Z\times\eu B$, $\eu N=0$, $\eu Y=S_0\times\eu B$, and $\eu
M=\bigcup_tN^\perp_t\times\{t\}$.  As a result we obtain a
$G$-invariant open neighbourhood $\eu O$ of $\eu Z$ in $\eu Y$ and an
isomorphism of relative Hamiltonian $G$-orbifolds
$$h\colon(\eu O,\omega_{\eu Y},\Phi_{\eu Y})\longrightarrow (\eu
M,\omega,\Phi)$$
onto an open neighbourhood of $\eu Z$ in $\eu M$ such that the
relevant commutative diagram commutes.  We can choose $\eu O$ to be of
the form $O\times[0,1]$, where $Z\subset O\subset S_0$.  In other
words, $h$ is (the track of) a symplectic isotopy of $O$ which fixes
$Z$ and $h_t$ maps $S_0$ to $N^\perp_t$.  After composing $h_t$ with
the map $h_0\inv\colon S_0\to S_0$ we may also assume that $h$ starts
at the identity.  We can view $h$ as an embedding of $\eu O$ into
$Y\times[0,1]$ and as such want to extend it to an isotopy of a full
neighbourhood $U$ of $Z$ in $Y$.  This is achieved by applying the
relative \emph{symplectic\/} embedding theorem to the embedding of the
relative symplectic manifold $\eu O$ into $Y\times[0,1]$.  To this end
we need to calculate the symplectic normal bundle $\eu E$ of
$S_0\times[0,1]$ in $Y\times[0,1]$.  Let $\pi$ denote the projection
$S_0\to Z$.  The restriction of $\eu E$ to $S_0\times\{0\}$ is equal
to $\pi^*N$, which is by definition equal to $Y$, considered as a
symplectic orbibundle over $S_0$.  The unit interval being
contractible, we conclude that $\eu E$ is isomorphic to
$Y\times[0,1]$, considered as a symplectic orbibundle over
$S_0\times[0,1]$.  By the relative symplectic embedding theorem, $h$
lifts to a symplectic embedding $H$ of $U'\times[0,1]$ into
$Y\times[0,1]$, where $U'\subset U$ is an open subset of $Y$
containing $Z$, as in the following commutative diagram:
$$
\xymatrix{
{\pr}_Z^*N\ar[d] & Y\times[0,1]\ar[l]\ar[d]\ar@{ >.>}[dr]^H \\
Z\times[0,1] & S_0\times[0,1]\ar@{ >->}[r]^h \ar[l]_{\pi\times\id} &
Y\times[0,1].
}
$$
By composing $H$ with the projection $Y\times[0,1]\to Y$ we find the
desired isotopy moving $TS_0|_Z$ to $TS_1|_Z$.

\emph{Step}~2.  We may henceforth assume that $S_0=TS_0|_Z=TS_1|_Z$.
Let $O$ be an open subset of $S_0$ containing $Z$ and let $h\colon
O\to S_1$ be any diffeomorphism onto an open subset of $S_1$ that
fixes $Z$ and satisfies $T_xh=\id$ for all $x$ in $Z$.  Such a map can
be found for instance by choosing a projection map $p$ of $TY|_Z$ onto
$TS_0|_Z$; the restriction of $p$ to $S_1$ has derivative equal to the
identity at all points of $Z$ and can therefore be locally inverted.
Let $H(x,b)=b\inv h(bx)$ be the isotopy deforming $h$ to its fibre
derivative; then $H_0=\id$, $H_1=h$, and $T_xH_b=\id$ for all $x$ in
$Z$.  Consider the forms $(1-t)\omega+tH_b^*\omega$ on $S_0$.
Applying Moser's trick with parameter $b$ as in Example
\ref{example:dmw} we find an open $O'\subset O$ and an isotopy
$I\colon O'\times[0,1]\to O$ such that $I|_Z=\id$, $I_0=\id$, and
$I_b^*H_b^*\omega=\omega$.  The isotopy $\check H\colon
O'\times[0,1]\to Y$ defined by $\check H_b=H_bI_b$ therefore satisfies
$\check H_0=\id$, $\check H_1=hI_1$ maps $O'\subset S_0$ to $S_1$ and
$\check H_b^*\omega=\omega$.  This symplectic isotopy can now be
extended to a neighbourhood $U'$ of $Z$ in $Y$ by use of the relative
symplectic embedding theorem, as in Step~1 above.
\end{proof}

\begin{lemma}\label{lemma:complement}
Let $K$ be a vector space and let $R$ and $N$ symplectic vector
spaces.  Let $V$ be the symplectic direct sum $R\oplus K\oplus
K^*\oplus N$\upn, where $K\oplus K^*$ carries the canonical symplectic
form.  Let $N^\perp$ be any complementary subspace to $N$ in
$K^*\oplus N$.  Then $R\oplus K\oplus N^\perp$ is a symplectic
subspace of $V$.  Let $P\colon V\to N$ be the linear projection with
kernel $R\oplus K\oplus N^\perp$.  Then the restriction of $P$ to
$(R\oplus K\oplus N^\perp)^\omega$ is a symplectic isomorphism onto
$N$.
\end{lemma}

\begin{proof}
It clearly suffices to prove this for $R=0$.  Let $d=\dim K$, $2n=\dim
V$.  There exists a symplectic basis $e_1$, $e_2,\dots$, $e_n$, $f_1$,
$f_2,\dots$, $f_n$ of $V$ such that $e_1$, $e_2,\dots$, $e_d$ form a
basis of $K$, $f_1$, $f_2,\dots$, $f_d$ form a basis of $K^*$, and
$e_{d+1}$, $e_{d+2},\dots$, $e_n$, $f_{d+1}$, $f_{d+2},\dots$, $f_n$
form a basis of $N$.  Furthermore, if $a_i=Pf_i$ for $i=1$, $2,\dots$,
$d$, then $f_1-a_1$, $f_2-a_2,\dots$, $f_d-a_d$ form a basis of
$N^\perp$.  Put $f_i'=f_i-a_i+\sum_{j=1}^d\alpha_{ij}e_j$ for $i=1$,
$2,\dots$, $d$, where $\alpha_{ij}=\frac1{2}\omega(a_i,a_j)$.  It is
easy to check that $e_1$, $e_2,\dots$, $e_d$, $f'_1$, $f'_2,\dots$,
$f'_d$ form a symplectic basis of $K\oplus N^\perp$, so $K\oplus
N^\perp$ is symplectic.

It is also easy to see that the vectors
\begin{align*}
e_i'&=e_i+\sum_{j=1}^d\omega(e_i,a_j)e_j,\\
f_i'&=e_i+\sum_{j=1}^d\omega(f_i,a_j)e_j,\\
\end{align*}
defined for $i=d+1$, $d+2,\dots$, $n$, are a symplectic basis of
$(K\oplus N^\perp)^\omega$.  Clearly, the projection $P$ sends $e_i'$
to $e_i$ and $f_i'$ to $f_i$ and therefore maps $(K\oplus
N^\perp)^\omega$ symplectically onto $N$.
\end{proof}

\section{A product formula}\label{section:product}

The following result was proved in a holomorphic context by Borel in
Appendix II of \cite{hi:to}.

\begin{theorem}\label{theorem:product}
Let $B$ be a compact almost complex orbifold and let $X$ be an
almost complex fibre orbibundle over $B$ with compact general
fibre $Y$ and orbibundle projection $\pi$.  Assume that the
structure group of $X$ can be reduced to a compact Lie group.
Then
$$
\RR(X,\pi^*E)=\RR(B,E)\RR(Y,\C)
$$
for every complex vector orbibundle $E$ over $B$.
\end{theorem}

We establish a slightly stronger result, namely an integration
formula, Theorem \ref{theorem:integration}, for the Todd form of the
vertical tangent bundle.  For simplicity we present the proof in the
manifold category; the proof for orbifolds is analogous.

\subsection{Cartan map and equivariant curvature}

Let $K$ be a compact (but not necessarily connected) Lie group, let
$\C[\k]$ be the graded algebra of polynomials on $\k$, and let $Y$ be
a $K$-manifold.  We denote by
$$
\Omega_K^k(Y)=\bigoplus_{i+2j=k}\bigl(\Omega^i(Y)\otimes
\C[\k]_j\bigr)^K
$$
the $\Z$-graded algebra of equivariant differential forms.  We also
consider the $\Z_2$-graded algebra
$$
\hat\Omega_K(Y)=\bigl(\Omega(Y)\otimes\C[[\k]]\bigr)^K
$$
of equivariant forms with coefficients in the formal power series
$\C[[\k]]$.  These algebras carry a differential of degree $1$ defined
by
$$
(\d_K\alpha)(\xi)=\d\alpha(\xi)-\iota(\xi_Y)\alpha(\xi).
$$
The cocycles in $\Omega_K^*(Y)$ and $\hat\Omega_K(Y)$ are denoted by
$\ca Z_K^*$ and $\hat{\ca Z}_K$, respectively, and the cohomology
groups by $H_K^*$ and $\hat H_K$, respectively.

Now let $P\to B$ be a $K$-principal bundle with connection
$\theta\in\Omega^1(P,\k)^{K}$.  The curvature of $\theta$ is the basic
two-form $F^\theta= \d\theta+\frac1{2}[\theta,\theta]
\in\Omega^2(P,\k)^K$.  It can be viewed as a $K$-equivariant map
$\k^*\to\Omega^2(P)$ and as such extends uniquely to an equivariant
multiplicative map $\C[[\k]]\to\Omega^*(P)$.  In other words, given
$\alpha\in\hat\Omega_K(P)$ we can substitute the curvature in the
$\k$-slot to get a $K$-invariant differential form $\alpha(F^\theta)$
and thus we obtain a map
$j^\theta\colon\hat\Omega_K(P)\to\Omega(P)^K$.  Let
$\hor^\theta\colon\Omega^*(P)\to\Omega^*_{\hor}(P)$ be the projection
onto the horizontal forms defined by the connection $\theta$.  The
composition
$$
\Car^\theta=\hor^\theta\circ
j^\theta\colon\hat\Omega_K(P)\longrightarrow 
\Omega^*_\basic(P)\cong\Omega^*(B)
$$
is called the \emph{Cartan map}.  Neither $\hor^\theta$ nor $j^\theta$
is a cochain map, but $\Car^\theta$ is.  Its restriction to
$\C[[\k]]\subset\hat\Omega_K(P)$ is known as the \emph{Chern-Weil
map}.

Consider the product $P\times Y$ and the associated bundle
$$
X=P\times^KY.
$$
Let $\pr_P\colon P\times Y\to P$ be the projection onto the first
factor and $\Car^{\pr_P^*\theta}\colon \hat\Omega_K(P\times
Y)\to\Omega^*(X)$ the Cartan map for the principal $K$-bundle $P\times
Y\to X$.  If $Y$ is compact, we can integrate forms over the fibres
and thus obtain a diagram
\begin{equation}\label{equation:integration}
\vcenter{\xymatrix{ 
{\hat\Omega_K}(Y)\ar[r]^-{\pr_Y^*}\ar[d]_{\int_Y}
& {\hat\Omega_K}(P\times
Y)\ar[r]^-{\Car^{\pr_P^*\theta}}\ar[d]_{\int_Y}
& {\Omega}^*(X)\ar[d]_{\int_Y} \\
{\hat\Omega_K}(\pt) \ar[r]^-{\pr_\pt^*}
& {\hat\Omega_K}(P)\ar[r]^-{\Car^\theta}
& {\Omega}^*(B).
}}
\end{equation}
We assert that this diagram is commutative.  The commutativity of the
square on the left is obvious; for the commutativity of the square on
the right it suffices to show that
$$
\int_Y\circ\;\hor^{\pr_P^*\theta}=\hor^{\theta}\circ\int_Y
\qquad\text{and}\qquad
\int_Y\circ\;j^{\pr_P^*\theta}=j^{\theta}\circ\int_Y.
$$
These identities follow from the fact that the forms
$\pr_P^*\theta\in\Omega^1(P\times Y,\k)^K$ and
$F^{\pr_P^*\theta}=\pr_P^*F^\theta\in\Omega^2(P\times Y,\k)^K$ have no
components in the $Y$-direction.

Now let $H$ be another compact Lie group and let $Q$ be an
$H$-principal bundle over $Y$.  Assume that the $K$-action on $Y$
lifts to an action on $Q$ that commutes with the $H$-action.  Choose a
$K$-invariant connection $\phi\in\Omega^1(Q,\h)^{H\times K}$ on $Q$.
Its \emph{$K$-equivariant curvature\/} is the $H$-basic
$K$-equivariant form $F_K^\phi\in\Omega^2_K(Q,\h)$ defined by
$$
F_K^\phi=\d_K\phi+\frac1{2}[\phi,\phi]=F^\phi-\Psi^\phi.
$$
Here $\Psi^\phi\colon Q\to \k^*$ is the map defined by
$\langle\Psi^\phi,\eta\rangle=\iota(\eta_Q)\phi$ for all $\eta\in\k$,
which, being $H$-invariant, descends to $Y$.  The connection $\theta$
on the $K$-bundle $P$ and the connection $\phi$ on the $H$-bundle $Q$
can be combined to a connection $\phi^\theta$ on the $H$-bundle
$P\times^KQ\to P\times^KY=X$, and the curvature of $\phi^\theta$ can
be expressed in terms of the equivariant curvature of $\phi$ in the
following manner.

\begin{lemma}\label{lemma:curvature}
\begin{enumerate}
\item\label{part:connection}
The $K$-horizontal part of $\pr_Q^*\phi$\upn, which is given by
$$
\hor^\theta(\pr_Q^*\phi)= \pr_Q^*\phi-
\langle\pr_Q^*\Psi^\phi,\pr_P^*\theta\rangle \in\Omega^1(P\times
Q,\k)^K,
$$
is a $K$-basic $H$-connection one-form on $P\times Q$ and represents a
connection one-form $\phi^\theta\in\Omega^1(P\times^KQ,\h)^H$.
\item\label{part:curvature}
The curvature form $F^{\phi^\theta}\in\Omega^2(P\times^KQ,\h)^H$ of
$\phi^\theta$\upn, regarded as a $K$-basic $\h$-valued form on
$P\times Q$\upn, is equal to
$\Car^{\pr_P^*\theta}(\pr_Q^*F^\phi_K)$\upn, where $F^\phi_K$ is the
$K$-equivariant curvature of $\phi$.
\end{enumerate}
\end{lemma}

\begin{proof}[Proof of \ref{part:curvature}]
By \ref{part:connection} the pullback of $F^{\phi^\theta}$ to $P\times
Q$ is equal to the curvature of the connection
$\hor^\theta(\pr_Q^*\phi)$, which is equal to the $K$-basic form
\begin{equation}\label{equation:curvature}
F^{\hor^\theta\phi}=
F^\phi+
\frac1{2}\bigl[\langle\Psi^\phi,\theta\rangle,
\langle\Psi^\phi,\theta\rangle\bigr]
-\langle\d\Psi^\phi,\theta\rangle- \langle\Psi^\phi,\d\theta\rangle
-\bigl[\langle\Psi^\phi,\theta\rangle,\phi\bigr],
\end{equation}
where we are suppressing the pullback maps from the notation.  On the
other hand
\begin{equation}\label{equation:substitute}
j^{\theta}F^\phi= F^\phi-\langle\Psi^\phi,\d\theta\rangle-
\frac1{2}\bigl\langle\Psi^\phi,[\theta,\theta]\bigr\rangle. 
\end{equation}
It is clear that the forms \eqref{equation:curvature} and
\eqref{equation:substitute} agree on $K$-horizontal vectors, so that
\eqref{equation:curvature} is the $K$-horizontal part of
\eqref{equation:substitute}.
\end{proof}

The $H$-Cartan map $\Car^\phi\colon\hat\Omega_H(Q)\longrightarrow
\Omega^*(Y)$ for $Q$ has a $K$-equivariant analogue
$$
\Car^\phi_K\colon\hat\Omega_{H\times
K}(Q)\longrightarrow\hat{\Omega}_K(Y),
$$
which is defined by $\Car^\phi_K=\hor^\phi\circ j^\phi_K$, where
$j^\phi_K(\alpha)=\alpha(F^\phi_K)$.  (Notice that $\alpha(F^\phi_K)$
is well-defined as a formal power series on $\k$ with values in
$\Omega^*(Y)$.)  It is a cochain map and therefore induces a map $\hat
H_{H\times K}(Q)\to\hat H_K(Y)$.  Its restriction to the subalgebra
$\C[[\h]]^H\subset\hat\Omega_{H\times K}(Q)$ is the
\emph{$K$-equivariant Chern-Weil map}.

\subsection{Integration formula for the Todd form}

Let $Y$ be an almost complex $K$-manifold and choose a $K$-invariant
Hermitian inner product and connection on the tangent bundle of $Y$.
These choices give rise to a unitary frame bundle $Q$ of $Y$ and a
principal connection $\theta$ on it.  Let $H=\U(n)$ and consider the
Todd series $\Td\in\C[[\h]]^H$, which is defined by
$$
\Td(x_1,x_2,\dots,x_n)=\prod_{j=1}^n\frac{x_j}{1-\exp(-x_j)}
$$
for $(x_1,x_2,\dots,x_n)$ in $i\R^n\cong\lie t$, the Cartan subalgebra
of $H$.  The form $\Td(Y)=\Car^\phi(\Td)\in\ca Z^*(Y)$ is the Todd
form of $Y$.

Let $P\to B$ be any $K$-principal bundle with connection $\phi$ and
consider the vertical tangent bundle $V=P\times^KTY$ of $X=P\times^KY$
over $B$.  Then $V$ is the $\C^n$-bundle on $X$ associated to the
principal $H$-fibration $P\times^KQ\to X$, on which we have the
connection $\phi^\theta$, and $\Td(V)=\Car^{\phi^\theta}(\Td)\in\ca
Z^*(X)$ is the Todd form of $V$.  As before let
$\int_Y\colon\Omega^*(X)\to\Omega^*(B)$ denote integration over the
fibres.

\begin{theorem}\label{theorem:integration}
The form $\int_Y\Td(V)\in\Omega^*(B)$ is a constant function on $B$.
Its value is equal to $\int_Y\Td(V)=\int_Y\Td(Y)$.
\end{theorem}

\begin{proof}
First we reduce the general case to the case where $K$ is connected.
Let $K^0$ be the identity component of $K$ and consider the finite
covers $\tilde B=P/K^0$ of $B$, $\tilde X=P\times^{K^0}Y$ of $X$ and
$\tilde V=P\times^{K^0}TY$ of $V$.  Then $\tilde X$ is a bundle over
$\tilde B$ with fibre $Y$ and $\tilde V$ is the pullback of $V$ under
the covering map $\tilde X\to X$.  Clearly $\Td(\tilde V)$ is the
pullback to $\tilde X$ of $\Td(V)\in\Omega^*(X)$ and $\int_Y\Td(\tilde
V)$ is the pullback to $\tilde B$ of $\int_Y\Td(V)\in\Omega^*(B)$, so
it is enough to show that $\int_Y\Td(\tilde V)=\int_Y\Td(Y)$.  We
may therefore assume $K$ to be connected.

Now consider $\Td_K(Y)=\Car^\phi_K(\Td)\in\hat\ca Z_K(Y)$, the
\emph{equivariant\/} Todd form of $Y$.  It follows from
\ref{part:curvature} of Lemma \ref{lemma:curvature} that $\Td(V)$ is
the image of $\Td_K(Y)$ under the composite map
$\Car^{\pr_P^*\theta}\circ\pr_Y^*\colon\hat\Omega_K(Y)\to\Omega(X)$.
We conclude from the commutativity of diagram
\eqref{equation:integration} that $\int_Y\Td(V)$ is the image of
$\int_Y\Td_K(Y)$ under the map $\Car^\theta\circ\pr_\pt^*$.  By the
Berline-Vergne equivariant index theorem \cite{be:eq},
$\int_Y\Td_K(Y)\in\hat{\Omega}_K(\pt)= \C[[\k]]^K$ is the equivariant
arithmetic genus of $Y$ (here we use that $K$ is connected), which by
\ref{part:almostcomplex} of Theorem \ref{theorem:rigid} is constant
and equal to $\int_Y\Td(Y)$.  Hence
$$
\int_Y\Td(V)= \Car^\theta\circ\pr_\pt^*\int_Y\Td(Y)=
\int_Y\Td(Y)
$$
as an element of $\ca Z^0(B)\cong\R$.
\end{proof}

\begin{proof}[Proof of Theorem \ref{theorem:product}]
We can write $X=P\times^KY$, where $K$ is a
compact Lie group and $P$ a principal $K$-bundle over $B$.  By
the Hirzebruch-Riemann-Roch theorem
$$
\RR(X,\pi^*E)=\int_{X}\Ch(\pi^*E)\Td(X)=\int_X\pi^*\Ch(E)\Td(X).
$$
Choose a connection on $P$ such that the induced connection on
$X$ is invariant under the almost complex structure.  Then we can
write $TX=V\oplus TB$, where $V=P\times^KTY$
is the vertical tangent bundle of $X$ over $B$.  Hence
$\Td(X)=\Td(V)\,\pi^*\Td(B)$ and
$$
\RR(X,\pi^*E)= \int_{X}\pi^*\bigl(\Ch(E)\Td(B)\bigr)\Td(V)=
\int_{B}\biggl(\Ch(E)\Td(B)\int_Y\Td(V)\biggr).
$$
The result now follows from Theorem \ref{theorem:integration} and
Hirzebruch-Riemann-Roch.
\end{proof}

\section{Notation}\label{section:notation}

\begin{tabbing}
\indent \= $M_\mu$; $M_0=M\qu G$\quad \= \kill 
\> $G$; $T$ \> compact connected Lie group; maximal torus \\
\> $\eu W$; $\Lambda$ \> Weyl group; integral lattice in $\lie t$ \\
\> $w_0$; $*$ \> longest Weyl group element; involution
$\mu\mapsto\mu^*=-w_0\mu$ of $\lie t^*$ \\
\> $w\odot\mu$ \> affine action $w\odot\mu=w(\mu+\rho)-\rho$ of
$\eu W$ \\
\> $\Lambda^*$; $\Lambda^*_+$ \> weight lattice $\Hom_\Z(\Lambda,\Z)$;
monoid of dominant weights \\
\> $\zeta_\mu$ \> character of $T$ defined by $\mu\in\Lambda^*$ \\
\> $\chi_\mu$ \> irreducible character of $G$ with highest weight
$\mu\in\Lambda^*_+$ \\
\> $\Rep G$; $\Ind_H^G$ \> representation ring; induction functor \\
\> $\tplus$; $\sigma$ \> positive Weyl chamber in $\lie t^*$; open
wall of $\tplus$ \\
\> $\star\sigma$; $\eu S_\sigma$ \> open star
$\bigcup_{\tau\suc\sigma}\tau$ of $\sigma$; natural slice in $\lie
g^*$ at $\sigma$ \\
\> $(M,\omega,\Phi)$ \> Hamiltonian $G$-orbifold with moment map \\ 
\> $\xi_M$ \> vector field on $M$ induced by $\xi\in\lie g$ \\
\> $\Delta$ \> Kirwan polytope $\Phi(M)\cap\tplus$ \\
\> $\inter\Delta$ \> relative interior of $\Delta$ \\
\> $L$ \> $G$-equivariant line orbibundle on $M$ \\
\> $\RR(M,L)$ \> equivariant index of $M$ with coefficients in $L$ \\
\> $N_L=N$ \> multiplicity function of $L$ \\
\> $M_\mu$; $M_0=M\qu G$ \> symplectic quotient of $M$ at $\mu$;
resp.\ $0$ \\
\> $L_\mu$; $L_0=L\qu G$ \> quotient orbibundle at $\mu$; resp.\ $0$
\\
\> $L_\mu\shift$ \> shifted quotient orbibundle at $\mu$ \\
\> $M_{\ge0}$; $L_{\ge0}$ \> symplectic cut of $M$ w.\ r.\ t.\ circle
action; cut bundle \\
\> $Y_\sigma$; $M_\sigma$ \> cross-section $\Phi\inv(\eu S_\sigma)$;
its saturation $GY_\sigma$ \\
\> $\S$; $\P$; $\F$ \> set of labels; polyhedron; open face \\
\> $D_\S$ \> Delzant space associated to set of labels $\S$ \\
\> $\Iso(E_1,E_2)$ \> isomorphisms from fibre bundle $E_1\to B_1$ to
$E_2\to B_2$ \\
\> $\Aut(E)$ \> automorphisms of a fibre bundle $E\to B$ \\
\> $\Aut_B(E)$ \> automorphisms of $E\to B$ that map each fibre to
itself
\end{tabbing}


\providecommand{\bysame}{\leavevmode\hbox to3em{\hrulefill}\thinspace}



\begin{thebibliography}{10}

\bibitem{ar:sy}
J.~Arms, J.~Marsden, and V.~Moncrief, \emph{Symmetry and bifurcations of
  momentum mappings}, Comm. Math. Phys. \textbf{78} (1981), 455--478.

\bibitem{at:sp}
M.~F. Atiyah and F.~Hirzebruch, \emph{Spin manifolds and group actions}, Essays
  on Topology and Related Topics (Geneva, 1969) (A.~Haefliger and
  R.~Narasimhan, eds.), Springer-Verlag, Berlin-Heidelberg-New York, 1970.

\bibitem{au:to}
M.~Audin, \emph{The topology of torus actions on symplectic manifolds},
  Progress in Mathematics, vol.~93, Birkh\"auser, Boston, 1991.

\bibitem{ba:ri2}
P.~Baum, W.~Fulton, and R.~MacPherson, \emph{Riemann-{R}och for singular
  varieties}, Inst. Hautes {\'E}tudes Sci. Publ. Math. \textbf{45} (1976),
  101--167.

\bibitem{be:eq}
N.~Berline and M.~Vergne, \emph{The equivariant index and {K}irillov's
  character formula}, Amer. J. Math. \textbf{107} (1985), 1159--1190.

\bibitem{bo:on}
R.~Bott, \emph{On the fixed point formula and the rigidity theorems of
  {W}itten}, Nonperturbative Quantum Field Theory (Carg\`ese, 1987)
  (G.~'t~Hooft et~al., eds.), NATO Advanced Science Institutes Series B:
  Physics, vol. 185, Plenum Press, New York-London, 1988, pp.~13--32.

\bibitem{br:res}
M.~Brion and M.~Vergne, \emph{Residue formul\ae, vector partition functions and
  lattice points in rational polytopes}, preprint LMENS-96-22, Ecole normale
  sup\'erieure, Paris, 1996.

\bibitem{br:qu2}
J.-L. Brylinski, \emph{Quantization commutes with reduction in geometric
  invariant theory}, preprint, Pennsylvania State University, 1996.

\bibitem{ca:eu}
S.~Cappell and J.~Shaneson, \emph{{E}uler-{M}aclaurin expansions for lattices
  above dimension one}, C. R. Acad. Sci. Paris S{\'e}r. I Math. \textbf{321}
  (1995), 885--890.

\bibitem{ca:ci}
A.~Canas da~Silva, Y.~Karshon, and S.~Tolman, \emph{Quantization of
  presymplectic manifolds and circle actions}, preprint, Massachusetts
  Institute of Technology, 1997, dg-ga/9705008.

\bibitem{de:ha}
T.~Delzant, \emph{Hamiltoniens p\'eriodiques et images convexes de
  l'application moment}, Bull. Soc. Math. France \textbf{116} (1988), 315--339.

\bibitem{du:he}
J.~J. Duistermaat, \emph{The heat kernel {L}efschetz fixed point formula for
  the spin-c {D}irac operator}, Progress in Nonlinear Differential Equations
  and Their Applications, vol.~18, Birkh\"auser, Boston, 1996.

\bibitem{du:sy}
J.~J. Duistermaat, V.~Guillemin, E.~Meinrenken, and S.~Wu, \emph{Symplectic
  reduction and {R}iemann-{R}och for circle actions}, Math. Res. Letters
  \textbf{2} (1995), 259--266.

\bibitem{fu:in}
W.~Fulton, \emph{Introduction to toric varieties}, Annals of Mathematics
  Studies, vol. 131, Princeton University Press, Princeton, 1993.

\bibitem{gu:re}
V.~Guillemin, \emph{Reduced phase spaces and {R}iemann-{R}och}, Lie Groups and
  Geometry in Honor of B. Kostant (Massachusetts Institute of Technology, 1994)
  (R.~Brylinski et~al., eds.), Progress in Mathematics, vol. 123, Birkh\"auser,
  Boston, 1995, pp.~305--334.

\bibitem{gu:ri}
\bysame, \emph{Riemann-{R}och for toric orbifolds}, J. Differential Geom.
  \textbf{45} (1995), no.~1, 53--73.

\bibitem{gu:co1}
V.~Guillemin and S.~Sternberg, \emph{Convexity properties of the moment
  mapping}, Invent. Math. \textbf{67} (1982), 491--513.

\bibitem{gu:ge}
\bysame, \emph{Geometric quantization and multiplicities of group
  representations}, Invent. Math. \textbf{67} (1982), 515--538.

\bibitem{gu:bi}
\bysame, \emph{Birational equivalence in the symplectic category}, Invent.
  Math. \textbf{97} (1989), no.~3, 485--522.

\bibitem{gu:sy}
\bysame, \emph{Symplectic techniques in physics}, Cambridge Univ. Press,
  Cambridge, 1990, second reprint with corrections.

\bibitem{hi:to}
F.~Hirzebruch, \emph{Topological methods in algebraic geometry}, third ed.,
  Grundlehren der mathematischen Wissenschaften, vol. 131, Springer-Verlag,
  Berlin-Heidelberg-New York, 1966.

\bibitem{je:lo2}
L.~C. Jeffrey and F.~C. Kirwan, \emph{Localization and the quantization
  conjecture}, Topology \textbf{36} (1997), no.~3, 647--693.

\bibitem{ka:ri}
T.~Kawasaki, \emph{The {R}iemann-{R}och theorem for complex {V}-manifolds},
  Osaka J. Math. \textbf{16} (1979), 151--159.

\bibitem{ki:con}
F.~C. Kirwan, \emph{Convexity properties of the moment mapping, {III}}, Invent.
  Math. \textbf{77} (1984), 547--552.

\bibitem{ki:pa}
\bysame, \emph{Partial desingularisations of quotients of nonsingular varieties
  and their {B}etti numbers}, Ann. of Math. (2) \textbf{122} (1985), 41--85.

\bibitem{ko:ap}
C.~Kosniowski, \emph{Applications of the holomorphic {L}efschetz formula},
  Bull. London Math. Soc. \textbf{2} (1970), 43--48.

\bibitem{le:sy2}
E.~Lerman, \emph{Symplectic cuts}, Math. Res. Letters \textbf{2} (1995),
  247--258.

\bibitem{le:co}
E.~Lerman, E.~Meinrenken, S.~Tolman, and C.~Woodward, \emph{Convexity by
  symplectic cuts}, Topology, to appear.

\bibitem{le:ha}
E.~Lerman and S.~Tolman, \emph{Hamiltonian torus actions on symplectic
  orbifolds and toric varieties}, preprint, Massachusetts Institute of
  Technology, 1995, dg-ga/9511008.

\bibitem{le:rr}
R.~Levy, \emph{The {R}iemann-{R}och theorem for complex spaces}, Acta Math.
  \textbf{158} (1987), no.~3--4, 149--188.

\bibitem{mc:re}
D.~McDuff, \emph{Remarks on the uniqueness of symplectic blowing up},
  Symplectic Geometry (Warwick, 1990) (D.~Salamon, ed.), London Mathematical
  Society Lecture Note Series, vol. 192, Cambridge University Press, Cambridge,
  1993, pp.~157--168.

\bibitem{mc:fr}
\bysame, \emph{From symplectic deformation to isotopy}, preprint, State
  University of New York, Stony Brook, 1996, dg-ga/9606004.

\bibitem{mc:in}
D.~McDuff and D.~Salamon, \emph{Introduction to symplectic topology}, Oxford
  Mathematical Monographs, Oxford University Press, New York, 1995.

\bibitem{me:sym}
E.~Meinrenken, \emph{Symplectic surgery and the {S}pin$^{\rm c}$-{D}irac
  operator}, Adv. in Math., to appear, dg-ga/9504002.

\bibitem{me:on}
\bysame, \emph{On {R}iemann-{R}och formulas for multiplicities}, J. Amer. Math.
  Soc. \textbf{9} (1996), 373--389.

\bibitem{me:lo}
E.~Meinrenken and C.~Woodward, \emph{A symplectic proof of {V}erlinde
  factorization}, preprint, Massachusetts Institute of Technology, 1996,
  dg-ga/9612018.

\bibitem{ot:re}
M.~Otto, \emph{A reduction scheme for phase spaces with almost {K}\"ahler
  symmetry. {R}egularity results for momentum level sets}, J. Geom. Phys.
  \textbf{4} (1987), 101--118.

\bibitem{sj:co}
R.~Sjamaar, \emph{Convexity properties of the moment mapping re-examined}, Adv.
  in Math., to appear, dg-ga/9408001.

\bibitem{sj:ho}
\bysame, \emph{Holomorphic slices, symplectic reduction and multiplicities of
  representations}, Ann. of Math. (2) \textbf{141} (1995), 87--129.

\bibitem{sj:sy}
\bysame, \emph{Symplectic reduction and {R}iemann-{R}och formulas for
  multiplicities}, Bull. Amer. Math. Soc. (N.S.) \textbf{33} (1996), 327--338.

\bibitem{sj:st}
R.~Sjamaar and E.~Lerman, \emph{Stratified symplectic spaces and reduction},
  Ann. of Math. (2) \textbf{134} (1991), 375--422.

\bibitem{ti:sy}
Y.~Tian and W.~Zhang, \emph{Symplectic reduction and analytic localization},
  preprint, Courant Institute, New York, 1996.

\bibitem{ve:qu}
M.~Vergne, \emph{Quantification g\'eom\'etrique et multiplicit\'es}, C. R.
  Acad. Sci. Paris S{\'e}r. I Math. \textbf{319} (1994), 327--332.

\bibitem{ve:eq}
\bysame, \emph{Equivariant index formula for orbifolds}, Duke Math. J.
  \textbf{82} (1996), 637--652.

\bibitem{wo:cl}
C.~Woodward, \emph{The classification of transversal multiplicity-free group
  actions}, Ann. Global Anal. Geom. \textbf{14} (1996), 3--42.

\end{thebibliography}
\end{document}